\documentclass[english]{article}
\usepackage{amssymb}
\usepackage{amsmath,color}
\usepackage{amsmath}
\usepackage{babel}

\def\thebibliography#1{\section*{References}\list
  {[\arabic{enumi}]}{\settowidth\labelwidth{[#1]}\leftmargin\labelwidth
    \advance\leftmargin\labelsep
    \usecounter{enumi}}
    \def\newblock{\hskip .11em plus .33em minus -.07em}
    \sloppy
    \sfcode`\.=1000\relax}
\newcommand{\refbook}[3]{{\sc #1}{\em\ #2}{\ #3}}
\newcommand{\refer}[5]{{\sc #1}{\ #2}{\em\ #3}{\bf\ #4}{\ #5}}

\newtheorem{lem}{Lemma}[section]

\newtheorem{teo}[lem]{Theorem}
\newtheorem{os}[lem]{Remark}
\newtheorem{defi}[lem]{Definition}
\newtheorem{prop}[lem]{Proposition}
\newtheorem{esem}[lem]{Example}

\newcommand{\qed}{\thinspace\null\nobreak\hfill\hbox{\vbox{\kern-.2pt\hrule
 height.2pt depth.2pt\kern-.2pt\kern-.2pt \hbox to2.5mm{\kern-.2pt\vrule
 width.4pt \kern-.2pt\raise2.5mm\vbox to.2pt{}\lower0pt\vtop
 to.2pt{}\hfil\kern-.2pt \vrule
 width.4pt \kern-.2pt}\kern-.2pt\kern-.2pt\hrule height.2pt depth.2pt
 \kern-.2pt}}\par\medbreak}

\newcommand{\R}{\mathbb{R}}

\newcommand{\C}{\mathbb{C}}

\newcommand{\N}{\mathbb{N}}

\newcommand{\eps}{\varepsilon}

\newcommand{\ov}{\overline}

\newcommand{\supp}{\emph{supp\,}}

\date{}
\begin{document}

\title{Scale invariant elliptic operators with singular coefficients}
\author{G. Metafune \thanks{Dipartimento di Matematica ``Ennio De Giorgi'', Universit\`a del Salento, C.P.193, 73100, Lecce, Italy.
e-mail:  giorgio.metafune@unisalento.it}\qquad N. Okazawa \thanks{Department of Mathematics, Tokyo University of Science, Japan. email: okazawa@ma.kagu.tus.ac.jp}\qquad M. Sobajima \thanks{Department of Mathematics, Tokyo University of Science, Japan and Dipartimento di Matematica ``Ennio De
Giorgi'', Universit\`a del Salento, C.P.193, 73100, Lecce, Italy. email: msobajima1984@gmail.com} \qquad C. Spina \thanks{Dipartimento di Matematica ``Ennio De Giorgi'', Universit\`a del Salento, C.P.193, 73100, Lecce, Italy.
e-mail:  chiara.spina@unisalento.it}}

\maketitle
\begin{abstract}
\noindent We show that a realization of the operator $L=|x|^\alpha\Delta +c|x|^{\alpha-1}\frac{x}{|x|}\cdot\nabla -b|x|^{\alpha-2}$ generates a semigroup in $L^p(\R^N)$ if and only if $D_c=b+(N-2+c)^2/4 > 0$ and  $s_1+\min\{0,2-\alpha\}<N/p<s_2+\max\{0,2-\alpha\}$, where $s_i$ are the roots of the equation $b+s(N-2+c-s)=0$, or $D_c=0$ and  $s_0+\min\{0,2-\alpha\} \le N/p \le s_0+\max\{0,2-\alpha\}$, where $s_0$ is the unique root of the above equation. The domain of the generator is also characterized.

\par

\end{abstract}

%%%%% Section 1 %%%%%
\section{Introduction}

In this paper we make a systematic investigation of the operator 
\begin{equation} \label{defL}
L=|x|^\alpha\Delta +c|x|^{\alpha-1}\frac{x}{|x|}\cdot\nabla -b|x|^{\alpha-2}
\end{equation}
in $L^p(\R^N)$, $N \ge 1$, $1<p<\infty$. Here $\alpha, b,c$ 
are unrestricted real numbers. 
%%%%%%==red letters==%%%%%%%%%%%%%%%%%%%%%%%%%%%%
Operators of the form  
%%%%%%==red letters==%%%%%%%%%%%%%%%%%%%%%%%%%%%%
$L(s)$
$=(s+|x|^\alpha)\Delta 
+ c|x|^{\alpha-1}\frac{x}{|x|}\cdot\nabla,\ s=0,1$, 
or operators containing a more general diffusion matrix 
in the second order part have been already studied 
in literature. See for example 
\cite{for-lor}, 
\cite{met-spi1}, 
\cite{met-spi2}, 
\cite{met-spi-tac}, 
\cite{met-spi-tac1}, 
\cite{soba-spin}, 
where generation results, domain characterization 
and spectral properties have been proved and 
\cite{met-spi3}, 
\cite{spina}, 
where kernel estimates have been deduced 
via weighted Nash inequalities. 
%%%%%%==red letters==%%%%%%%%%%%%%%%%%%%%%%%%%%%%
Operators of the form 
\eqref{defL} with $\alpha = 0$ 
have been studied in $L^{p}$-spaces with weight $|x|^{-\beta}$ 
for real $\beta$ (see 
\cite{AGG06}, 
\cite{OSY}). 
\\
In order to treat the singularity at zero we introduce $\Omega=\R^N\setminus\{0\}$ and define
$C_c^\infty(\Omega)$ as the space of infinitely continuously differentiable functions with compact support in $\Omega$.  
We define $L_{min}$ as the closure in $L^p(\R^N)$ of $(L,C_c^\infty(\Omega))$ and $L_{max}=(L,D_{max}(L))$ where

\begin{equation} \label{defDmax}
D_{max}(L)=\{u\in W^{2,p}_{loc}(\Omega)\cap L^p(\R^N):\ Lu\in L^p(\R^N)\}.
\end{equation}
The domain of $L_{min}$ will be denoted by $D_{min}(L)$. Note that if $u\in D_{max}(L)$ and $f=Lu$ the equation $Lu=f$ is satisfied in the sense of distributions in $\Omega$ rather than in $\R^N$.
\noindent
We study when suitable realizations of $L$ between $L_{min}$ and $L_{max}$ generate a semigroup in $L^p(\R^N)$.
The introduction of  $C_c^\infty(\Omega)$ instead of $C_c^\infty(\R^N)$ is unavoidable to treat the singularity at $0$ but sometimes leads 
%%%%%%%%%%%%%%%%%%==page 2==%%%%%%%%%%%%%%%%%%%%%
to unnatural difficulties. For example, 
if $\alpha=b=c=0$ and $N \ge 3$, then the Laplacian with domain 
$W^{2,p}(\R^N)$ coincides with $\Delta_{min}$ 
if and only if $p \le N/2$ and with $\Delta_{max}$ 
if and only if $p \ge N/(N-2)$. Similar problems happen 
when $C_c^\infty(\R^N) \subset D_{max}(L)$ 
(depending on $\alpha,b,c,p$) and this explains 
why we need also intermediate operators 
between $L_{min}$ and $L_{max}$. 
\noindent
When $\alpha=c=0$, $L$ becomes the Schr\"odinger operator with inverse square potential which is widely studied in the literature. A famous result in \cite{BG} shows that the parabolic equation $u_t=Lu$ presents instantaneous blow-up for positive solutions when 
%%%%%%==red letters==%%%%%%%%%%%%%%%%%%%%%%%%%%%%
$D_{0} := b +(N-2)^2/4 <0$, 
where $4D_{0}$ is the discriminant of the quadratic equation 
$$
      f_{0}(s) := -s^2 + (N-2)s + b = 0. 
$$
%$b +(N-2)^2/4 <0$. 
In the general case we show that the elliptic equation 
$\lambda u-Lu=f$, 
%%%%%%==red letters==%%%%%%%%%%%%%%%%%%%%%%%%%%%%
with
$\lambda, f \ge 0$,  has no positive solution if 
$\alpha \neq 2$ and 
%%%%%%==red letters==%%%%%%%%%%%%%%%%%%%%%%%%%%%%
$D_{c} := b+(N-2+c)^2/4 <0$.
The case $\alpha=2$ is 
special in the whole paper and the above restriction 
%%%%%%==red letters==%%%%%%%%%%%%%%%%%%%%%%%%%%%%
%{\color{red} 
%$D_{c}<0$% 
%}
is not necessary.
\noindent
We obtain positive results under the assumption 
%%%%%%==red letters==%%%%%%%%%%%%%%%%%%%%%%%%%%%%
$D_{c} \ge 0$.
%$b+(N-2+c)^2/4 > 0$. 

\noindent
In order to formulate our main results we introduce the 
%%%%%%==red letters==%%%%%%%%%%%%%%%%%%%%%%%%%%%%
quadratic function
\begin{equation} \label{deff}
f(s)=b+s(N-2+c-s) 
%%%%%%==red letters==%%%%%%%%%%%%%%%%%%%%%%%%%%%%
\;=-s^{2}+(N-2+c)s+b
\end{equation}
whose discriminant is $4D_c$. Its roots  are
%%%%%%==red letters==%%%%%%%%%%%%%%%%%%%%%%%%%%%%
$s_{1}$, $s_{2}$ 
$(s_{1} < s_{2})$ given by 
\begin{equation} \label{defs}
%%%%%%==red letters==%%%%%%%%%%%%%%%%%%%%%%%%%%%%
s_{1}=\frac{N-2+c}{2}-\sqrt{b+\left(\frac{N-2+c}{2}\right)^2}, 
\quad
s_{2}=\frac{N-2+c}{2}+\sqrt{b+\left(\frac{N-2+c}{2}\right)^2}. 
\end{equation}
\noindent
Note that $f$ has the maximum at $s_0= (N-2+c)/2$ with $f(s_0)=D_c$.

\noindent Our main result in the case $D_c>0$ is the following which summarizes Theorems  \ref{nosemi}, \ref{gen-int}.

\begin{teo} \label{main}
Let $1<p<\infty$, $\alpha \neq 2$, 
%%%%%%==red letters==%%%%%%%%%%%%%%%%%%%%%%%%%%%%
$D_{c}=$
$b+(N-2+c)^2/4>0$. 
%$b+(N-2+c)^2/4 > 0$. 
%%%%%%==red letters==%%%%%%%%%%%%%%%%%%%%%%%%%%%%
Then
a suitable realization of 
$L_{min} \subset L_{int} \subset L_{max}$ 
generates a semigroup in $L^p(\R^N)$ 
if and only if 
%%%%%%==red letters==%%%%%%%%%%%%%%%%%%%%%%%%%%%%
$$
s_1+\min\{0,2-\alpha\}<N/p<s_2+\max\{0,2-\alpha\}. 
$$
In this case the 
%%%%%%==red letters==%%%%%%%%%%%%%%%%%%%%%%%%%%%%
generated
semigroup is bounded analytic and positive. 
The domain of $L_{int}$ is given by equation (\ref{Dreg}).
\end{teo}
In general the semigroup is not contractive. The case $\alpha=2$ is special and much simpler: no restriction on $N/p$ is needed, see Proposition \ref{L2}.

\noindent
We observe that $L$ generates a semigroup in some $L^p(\R^N)$ if and only if the open intervals 
$(s_1+\min\{0,2-\alpha\},s_2+\max\{0,2-\alpha\})$  and $(0,N)$ intersect. This is always the case when $b > 0$ since $s_1$ and $s_2$ have opposite signs but easy examples show that the contrary can happen if $b \le 0$, see the last section of this paper. In such cases no realization of $L$ between $L_{min}$ and $L_{max}$ is a generator but it can happen that $L$ endowed with a suitable domain is a generator. We refer the reader to \cite{SW} where it is shown that for every $b \in \R$ a suitable realization of $\Delta-b|x|^{-2}$ is self-adjoint and non-positive in $L^2(\R^N)$. \\

\noindent In the critical case $D_c=0$ we prove the following result in Section 6.
\begin{teo} \label{main1}
Let $1<p<\infty$, $\alpha \neq 2$, 
%%%%%%==red letters==%%%%%%%%%%%%%%%%%%%%%%%%%%%%
$D_{c}=$
$b+(N-2+c)^2/4=0$ and $s_0=\frac{N-2+c}{2}$. 
Then
a suitable realization of 
$L_{min} \subset L_{int} \subset L_{max}$ 
generates a semigroup in $L^p(\R^N)$ 
if and only if 
%%%%%%==red letters==%%%%%%%%%%%%%%%%%%%%%%%%%%%%
$$
s_0+\min\{0,2-\alpha\} \le N/p \le s_0+\max\{0,2-\alpha\}. 
$$
In this case the 
%%%%%%==red letters==%%%%%%%%%%%%%%%%%%%%%%%%%%%%
generated
semigroup is bounded analytic and positive. 
The domain of $L_{int}$ is given by equations (\ref{Lint-critical}), (\ref{Lint-criticalMaggiore}). 
\end{teo}
 Note that the endpoints are included in Theorem \ref{main1} but excluded in Theorem \ref{main}. 
We also point out that 
the validity of the equalities $L_{int}=L_{min}$ and $L_{int}=L_{max}$ is also characterized through the paper.\\

\noindent
The paper is organized as follows. In Section 2 we prove and recall some preliminary results. In Section 3 we partially generalize the results in \cite{BG} by showing that if $b+(N-2+c)^2/4 <0$  the equation $ u-Lu=f$ has no positive distributional solutions for certain positive $f$ with compact 
%%%%%%%%%%%%%%%%%%%%%%%%%%==page 3==%%%%%%%%%%%%%
support. In Section 4 we show that $L_{min}$ generates an analytic semigroup
when $s_1+2-\alpha <N/p <s_2+2-\alpha$ and characterize 
its domain, using Rellich inequalities from 
\cite {met-soba-spi}. The proof is done first for very 
large $b>0$ showing sectoriality and then extended to 
the precise range above using a perturbation argument 
%%%%%%==red letters==%%%%%%%%%%%%%%%%%%%%%%%%%%%%
in \cite{Sobajima-per}, as
%of the third author, 
stated in the  Appendix. Generation results for $L_{max}$ are deduced by duality. The sharpness of the above intervals is then shown using the  asymptotics of special radial solutions: in particular the "only if" part of Theorem \ref{main} is proved in Theorem \ref{nosemi}. The operator $L_{int}$ is introduced in Section 5. Using the results of Section 4 for $L_{min}$ 
we give a proof of the "if" part of 
Theorem \ref{main}, see Theorem \ref{gen-int} for a more precise formulation. The critical case $D_c=0$ is studied in Section 6, using the methods of Section 5 but adding a logarithmic term in the weighted estimates. In contrast with Section 5, we do not prove directly the resolvent estimates in $\R^N$ but first show a weaker form in the unit ball and then improve them in the whole space by scaling. In Section 7 we present some examples. It is worth mentioning that our main results, specialized to the case of Schr\"odinger operators with inverse square potentials, yields more precise results than those already known. In particular we show that the semigroup exists in the same range of $p$ as in \cite{Li-sobol} when $D_c>0$ but we are able to characterize the domain of the generator  in addition to
the domain of the form. The precise range of existence  of the semigroup is also given in the critical case and seems to be new.

\smallskip

\noindent Our result are valid when $N=1$ with $[0,\infty[$ instead of $\R$. In the statements, however,  we keep the notation $\R^N$ even when $N=1$. Accordingly $\Omega=]0,\infty[$ and all balls $B_r$ should be replaced by the intervals $]0,r[$. With these (formal) changes all proofs hold in the one-dimensional case with, at most, some simplifications.

\bigskip
\noindent\textbf{Notation.} We use $\Omega$ for $\R^N \setminus \{0\}$ and for $]0, \infty[$ when $N=1$.
$C_c^\infty(V)$ denotes the space of infinitely continuously differentiable functions with compact support in $V$. We adopt standard notation for $L^p$ and Sobolev spaces. The unit sphere in $\R^N$ is denoted by $S^{N-1}$ and $B_r$ stands for the ball with center at $0$ and radius $r$.

%%%%% Section 2 %%%%%
\section{Preliminary results}
Here we collect some known or simple fact necessary to our analysis.
Observe that if $I_\lambda u(x)=u(\lambda x)$ for $\lambda >0$, then $(I_\lambda)^{-1}LI_\lambda=\lambda^{2-\alpha}L$. Note that $L$ is scale invariant when $\alpha=2$. Other symmetry properties follow from the use of the Kelvin transform. Let $Tu(x)=|x|^{2-N}u(x|x|^{-2})$. A straightforward but tedious computation shows that
\begin{equation} \label{kelvin}
T^{-1}LT=|x|^{4-\alpha}\Delta-c|x|^{3-\alpha}\frac{x}{|x|}\cdot \nabla +\left (c(2-N)-b\right )|x|^{2-\alpha}.
\end{equation}
In particular the power $\alpha$ is changed into $4-\alpha$. Many proofs will be subdivided according to $\alpha<2$ and $\alpha>2$. If $\alpha <2$ the degeneracy at infinity is easy to treat but that at the origin is the real source of the difficulties. Conversely when $\alpha >2$, using the Kelvin transformation and noticing that it maps the unit ball into its complement, one can study only case, e.g., $\alpha <2$ and reduce the other to it. Observe however that the Kelvin transform is an isomorphism in $L^p$ if and only if $p=2N/(N-2)$.
\noindent
Let us show the closedness of $L_{min}$ and $L_{max}$.

\begin{prop} 
The operator $L_{max}$ is closed and $(L,C_c^\infty(\Omega))$ is  closable.
\end{prop}
{\sc Proof.} The closedness of $L_{max}$ is an immediate consequence of local elliptic regularity, since $L$ has regular coefficients outside the origin. Since $C_c^\infty(\Omega) \subset D_{max}(L)$, the closability of $(L,C_c^\infty(\Omega))$ follows from the closedness of $L_{max}$.
\qed

\noindent
Next we introduce the formal adjoint   

\begin{equation} \label{defLtilde}
\tilde{L}=|x|^\alpha\Delta +\tilde{c}|x|^{\alpha-1}\frac{x}{|x|}\cdot\nabla -\tilde{b}|x|^{\alpha-2}
\end{equation}
with 
\begin{equation} \label{tilde}
\left\{
\begin{array}{ll}
 \tilde{c}= 2\alpha-c\\ 
  \tilde{b}=b+(c-\alpha)(\alpha-2+N)
\end{array}\right.
\end{equation}
acting on $L^{p'}(\R^N)$. Observe that the function 
\begin{equation} \label{ftilde}
\tilde f (s)=\tilde b+s(N-2+\tilde c-s)=f(N+\alpha-2-x)
\end{equation}
 defined as in (\ref{deff}) and relative to $\tilde L$, has roots $\tilde{ s}_i=s_i+\alpha-c$, $i=1,2$, where $s_1, s_2$ are defined in (\ref{defs}) and that its discriminant $\tilde{b}+(N-2+\tilde c)^2/4$ coincides with that of $f$ (that is with  $b+(N-2+c)^2/4$).

\noindent Then we have 
%%%%%%%%%%%%%%%%%%%%%%==page 4==%%%%%%%%%%%%%%%%%
\begin{prop}\label{adjoint}
\[
(\tilde{L}_{max})=(L_{min})^*
\quad\textrm{and}\quad
(\tilde{L}_{min})=(L_{max})^*
\]
\end{prop}
{\sc Proof.} The first identity is immediate consequence of the definitions and of interior elliptic regularity. Taking the adjoiont in the equality $(\tilde{L}_{max})=(L_{min})^*$ one obtains $(\tilde{L}_{max})^*=(L_{min})$, by the closedness of $L_{min}$, which is the second one (with the roles of $L$ and $\tilde L$ interchanged).
\qed

\noindent
As pointed out in the Introduction, the case $\alpha=2$ is quite special. Let us state the result in the next proposition (see \cite[Section 6]{met-soba-spi} for the proof).

\begin{prop} \label{L2}
Consider the operator $L$ defined in (\ref{defL}) with $\alpha=2$ and let $1< p< \infty$. Then $L_{max}=L_{min}$ generates an analytic semigroup of positive operators $(T(t))_{t \ge 0}$ in $L^p$ satisfying $\|T(t)\|_p \le e^{(b-\omega_p) t}$, 
%%%%%%==red letters==%%%%%%%%%%%%%%%%%%%%%%%%%%%%
$\omega_p =f(N/p)-b 
= \frac{N}{p}\left(\frac{N}{p'}-2+c\right)$.
%$\omega_p=\frac{N}{p^2}\left[p(N-2+c)-N\right]$. 
Finally
 \begin{align*}
D_{max}(L)=& \{u \in L^p(\R^N)\cap W^{2,p}(\Omega), \  
 |x|\nabla u, |x|^2 D^2u \in L^p(\R^N)\},
\end{align*}
\end{prop}
\noindent
When $\alpha \in \R$ we introduce the domain
\begin{equation}\label{domain}
D_{p,\alpha}=\{u\in L^p(\R^N) \cap  W^{2,p}_{loc}(\Omega),\ |x|^\alpha D^2u,\ |x|^{\alpha-1} \nabla u,\ |x|^{\alpha-2}u\in L^p(\R^N)\}
\end{equation}
endowed with its canonical norm and note that it 
to that in the above proposition when $\alpha=2$. Note that extra integrability condition for $u$ is relevant near $0$ when $\alpha <2$ and near infinity when $\alpha >2$.

\begin{lem}  \label{densita}
The space $C_c^\infty(\Omega)$ is dense in $D_{p,\alpha}$. The following interpolation property holds in $D_{p,\alpha}$: there exist $C, \eps_0$ depending on $N,p,\alpha$ such that for every $u \in D_{p,\alpha}$ and $\eps \le \eps_0$
\begin{equation} \label{int}
\||x|^{\alpha-1}\nabla u\|_p \le \eps \|Lu\|_p +\frac{C}{\eps}\||x|^{\alpha-2}u\|_p.
\end{equation}
\end{lem}
{\sc Proof.} 
 Let us first observe that a function $u\in
W^{2,p}(\Omega)$ with compact support in $\Omega$ can be approximated by a
sequence of $C^\infty$ functions with compact support in $\Omega$ in the
$D_{p,\alpha}$ norm, by using standard mollifiers. Let $u$ in $D_{p,\alpha}$  and $\eta_n$ be 
smooth functions such that $\eta_n=1$ in $B_{n}\setminus B_{1/n}$, $\eta_n=0$ in
$\R^N\setminus (B_{2n} \cup B_{1/2n})$, $0\leq \eta_n \leq 1$ and $|\nabla \eta_n(x)| \le C|x|^{-1}$, $|D^2\eta_n(x)| \le C|x|^{-2}$. If $u\in D_{p,\alpha}$, then
$u_n=\eta_n u$ are compactly supported functions in $W^{2,p}(\Omega)$,
$u_n\to u$ in $L^p(\R^N)$, $|x|^{\alpha-2}u_n\to
|x|^{\alpha-2}u$ in $L^p(\R^N)$ by dominated convergence.
Concerning the convergence of the derivatives we have
$$|x|^{\alpha-1}\nabla u_n=|x|^{\alpha-1}\nabla\eta_n(x) u+|x|^{\alpha-1}\eta_n(x)\nabla u.$$
As before $|x|^{\alpha-1}\eta_n(x)\nabla u
\to |x|^{\alpha-1}\nabla u$ in $L^p(\R^N)$. For the left term,
since $\nabla \eta_n (x)$ can be different from zero only for $1/2n \le
|x| \le 1/n$ or $n \le
|x| \le 2n$ we have
$$|x|^{\alpha-1}|\nabla\eta_n (x) ||u|\leq C|x|^{\alpha-2}|u|(\chi_{\{(2n)^{-1}\leq |x|\leq n^{-1}\}}+\chi_{\{n\leq |x|\leq 2n\}}),$$
and the right hand side tends to $0$ as $n\to\infty$. A similar
argument shows the convergence of the second order derivatives in
the weighted  $L^p$ norm and the proof of the density is complete. Concerning (\ref{int}) we observe that the weaker inequality
\begin{equation} \label{int1}
\||x|^{\alpha-1}\nabla u\|_p \le \eps \||x|^\alpha D^2 u\|_p +\frac{C}{\eps}\||x|^{\alpha-2}u\|_p
\end{equation}
holds in $C_c^\infty(\Omega)$  by \cite[Lemma 4.4]{met-soba-spi}. By applying the classical Calder\'on-Zygmund $\|D^2v\|_p \le C\|\Delta v\|_p$ to $v=|x|^\alpha u$ and using (\ref{int1}) to interpolate the gradient terms 
we get 
\begin{align*}
&\||x|^\alpha D^2 u\|_p \leq
  C\left( \|D^2(|x|^\alpha u)\|_p+\||x|^{\alpha-1}\nabla u\|_p +\||x|^{\alpha-2}u\|_p\right) \\
& \leq
C\left( \|\Delta(|x|^\alpha
u)\|_p+\||x|^{\alpha-1}\nabla u\|_p +\||x|^{\alpha-2}u\|_p\right)\leq
C\left( \||x|^\alpha \Delta
u\|_p+\||x|^{\alpha-1}\nabla u\|_p +\||x|^{\alpha-2}u\|_p\right)\\&\leq
C\left( \||x|^\alpha L
u\|_p+\||x|^{\alpha-1}\nabla u\|_p +\||x|^{\alpha-2}u\|_p\right)\leq
C\left( \||x|^\alpha L
u\|_p+\eps \||x|^{\alpha}D^2 u\|_p +C_\eps\||x|^{\alpha-2}u\|_p\right).
\end{align*}
Taking $\eps$ small, (\ref{int}) follows. By the density of $C_c^\infty (\Omega)$ in $D_{p,\alpha}$ the proof is complete.
\qed

\noindent The following lemma is useful to study the equality $L_{min}=L_{max}$.

\begin{lem}\label{aux}For every $\alpha\neq 2$, 
\[
D_{max}(L)\cap D(|x|^{\alpha-2})=D_{p, \alpha}\subset D_{min}(L).
\]
\end{lem}
{\sc Proof.} The inclusion $D_{p,\alpha} \subset D_{max}(L)\cap D(|x|^{\alpha-2})$ is evident and  $D_{p,\alpha} \subset D_{min}(L)$ follows from the density of $C_c^\infty(\Omega)$ in $D_{p,\alpha}$. Let $u\in D_{max}(L)\cap D(|x|^{\alpha-2})$, we  define $v=|x|^{\alpha-2}u \in L^p(\R^N)$ and note that $Lu=\tilde{L}v-bv$
where 
$$
\tilde L=|x|^2\Delta +(4-2\alpha+c)x\cdot\nabla +(2-\alpha)(N-\alpha+c).
$$
Then  $v \in D_{max}(\tilde L)$ and  therefore, by Proposition \ref{L2}, $|x| \nabla v, |x|^2 D^2 v \in L^p(\R^N)$. This yields $u \in D_{p,\alpha}$ and concludes the proof. \qed

\noindent We need also the asymptotic behavior of the solutions of certain singular ordinary differential equations  related to Bessel equations. We recall that the numbers $s_1, s_2$ are defined in (\ref{defs}).

\begin{lem}\label{bhv}
Let $\alpha \neq 2$, $b,c\in\R$ and $\lambda>0$ and assume that $k:=b+(\frac{N-2+c}{2})^2 \ge 0.$ The differential equation
\begin{equation}
\lambda u-r^\alpha\left(u''+\frac{N-1+c}{r}u'-\frac{b}{r^2}u\right)=0
\end{equation}
has  two positive solutions 
$u_1$ and $u_2$ 
with the following behavior: if $\alpha<2$ and $k>0$, then
\begin{align}
\label{u1-behave}
u_1(r)\approx r^{-s_1}\quad {\rm near }\ 0, 
\qquad 
u_1(r)\approx r^{-\frac{N-2+c}{2}+\frac{\alpha-2}{4}}e^{\frac{2}{2-\alpha}\lambda^{\frac{1}{2}}r^{\frac{2-\alpha}{2}}} \quad {\rm near }\ \infty, 
\\
\label{u2-behave}
u_2(r)\approx r^{-s_2}\quad {\rm near }\ 0, 
\qquad 
u_2(r)\approx r^{-\frac{N-2+c}{2}+\frac{\alpha-2}{4}}e^{-\frac{2}{2-\alpha}\lambda^{\frac{1}{2}}r^{\frac{2-\alpha}{2}}} \quad {\rm near }\ \infty;
\end{align}
if $\alpha <2$ and $k=0$, then
\begin{align}
\label{u1-behave-res}
u_1(r)\approx r^{-\frac{N-2+c}{2}} \quad {\rm near }\ 0, 
\qquad 
u_2(r)\approx- r^{-\frac{N-2+c}{2}}\log r\quad {\rm near }\ 0, 
\end{align}
and the behavior at $\infty$ is as above. When $\alpha>2$, (\ref{u1-behave}), (\ref{u2-behave}) and (\ref{u1-behave-res}) hold with $0$ and $\infty$ interchanged in each of them.
\end{lem}
{\sc Proof.} Defining  $\tilde{u}(r)=r^{\frac{N-2+c}{2}}u(r)$ we obtain
\begin{align*}
r^2\tilde{u}''(r)+r\tilde{u}'(r)
&=r^{\frac{N-2+c}{2}}\left(r^2u''(r)+r(N-1+c)u'(r)+\left(\frac{N-2+c}{2}\right)^2u(r)\right)
\\
&=r^{\frac{N-2+c}{2}}\left(\lambda r^{2-\alpha}+b+\left(\frac{N-2+c}{2}\right)^2 \right)u(r)
=\left(\lambda r^{2-\alpha}+k \right)\tilde{u}(r).
\end{align*}
Setting $v(r)=\tilde{u}\bigl(c r^{\gamma}\bigr)$,
we have 
\begin{align*}
r^2v''(r)+rv'(r)
&=
\gamma^2
\left[
c^2r^{2\gamma}
\tilde{u}''(cr^{2\gamma})
+
cr^\gamma
r^{\frac{2}{2-\alpha}}
\tilde{u}'(cr^\gamma)\right]=
\gamma^2 \left(\lambda c^{2-\alpha}r^{\gamma(2-\alpha)}+k \right )
v(r).
\end{align*}
Choosing 
\begin{equation} \label{constants}
\gamma=\frac{2}{2-\alpha} \qquad c=\left (\frac{(2-\alpha)^2}{4\lambda}\right )^{\frac{1}{2-\alpha}}
\end{equation}
it follows that $v$ satisfies the Bessel equation
\begin{equation} \label{mbessel}
r^2v''(r)+rv'(r)=
(\nu^2 +r^2)v(r)
\end{equation}
with $\nu^2=\bigl(\frac{2}{2-\alpha}\bigr)^2k$. 
%\begin{equation}\label{mbessel}
%r^2v''+rv'-(\nu^2+r^2)v=0.
%\end{equation}
for which 
the modified Bessel functions $ I_\nu$ and $K_\nu$ 
constitute a basis. We 
note that both $I_\nu$ and $K_\nu$ are positive, 
$I_\nu$ is monotone increasing and $K_\nu$ is monotone decreasing. Moreover, by \cite[Section 7.5]{BW},
\begin{align*}
&I_\nu(r) \approx r^\nu \quad {\rm near}\  0, \qquad I_\nu(r) \approx \frac{e^r}{\sqrt r} \quad {\rm near}\ \infty \\
&K_\nu(r) \approx r^{-\nu}\  (\nu >0), \quad K_\nu \approx -\log r\  (\nu=0) \quad {\rm near}\  0, \qquad K_\nu(r) \approx \frac{e^{-r}}{\sqrt r} \quad {\rm near}\ \infty .
\end{align*}
Since
\[
u_1(r):=
r^{-\frac{N-2+c}{2}}I_{\nu}\left (\bigl |1-\frac{\alpha}{2}\bigr|\lambda^{\frac{1}{2}}r^{\frac{2-\alpha}{2}}\right), \qquad 
u_2(r):=
r^{-\frac{N-2+c}{2}}K_{\nu}\left(\bigl |1-\frac{\alpha}{2}\bigr|\lambda^{\frac{1}{2}}r^{\frac{2-\alpha}{2}}\right). 
\]
all the assertions readily follow.
\qed

\noindent The following elementary consequence of H\"older inequality will be used several times;  we state it here to fix the parameters. 
\begin{lem} \label{intgamma}
Assume that $\mu$ is a measure and that all powers are integrable with respect to $\nu$. If  $\gamma_1 \le \gamma_2 \le \gamma_3$, then
$$
\||x|^{\gamma_2}\|_p \le \||x|^{\gamma_1}\|_p^\tau \||x|^{\gamma_3}\|_p^{1-\tau}
$$ with $\tau=\frac{\gamma_3-\gamma_2}{\gamma_3-\gamma_1}$ and the norms are taken in $L^p$ with respect to $\mu$.
\end{lem}

%%%%% Section 3 %%%%%
\section{Non existence of positive solutions for  $b+\left(\frac{N-2+c}{2}\right)^2<0$}
 A famous result in \cite{BG}, see also \cite{CM}, \cite{BV} for different proofs,  states that the equation $u_t=\Delta u-b|x|^{-2}u$ does not admit  positive solution if $b+(N-2)^2/4 <0$. Note that 
%%%%%%==red letters==%%%%%%%%%%%%%%%%%%%%%%%%%%%%
$b_0 : =(N-2)^2/4$
is the best constant in Hardy inequality in $L^2(\R^N)$. A detailed analysis of the solution for $b \ge -b_0$ is done in \cite{VZ}, including an investigation of oscillating solutions for $b <-b_0$.
The above result does not say that the symmetric operator $\Delta-b|x|^{-2}$ does not generate a semigroup for $b <-b_0$. In fact in \cite{SW} it is proved that for every $b \in \R$ the operator above, endowed with a suitable domain, generates a self-adjoint semigroup of positivity preserving operators. However the semigroup solution so produced, satisfies the parabolic equation in a weaker sense than in \cite{BG}, namely it is a distributional solution in a set $\R^N \setminus F$ where $F$ is a closed set of measure zero.
In this section we show that a phenomenon similar to that of \cite{BG} occurs, independently of $\alpha$. We prove it for the elliptic problem rather than for the parabolic one.

%%%%%%%%%%%%%%%%%%==page 5==%%%%%%%%%%%%%%%%%%%%%
\begin{teo}
Let $\alpha \neq 2$, $b+\left(\frac{N-2+c}{2}\right)^2<0$. 
Then, for every $\lambda>0$, there exists a radial function 
$0 \leq \phi \in C_c^\infty(\Omega)$, $\phi\not\equiv 0$, 
such that the problem 
\begin{equation}  \label{oscillante}
\lambda u- Lu =\phi
\end{equation}
does not admit any positive distributional solution 
in $\Omega$.
\end{teo}
{\sc Proof.} Assume that $\alpha <2$ and 
that there exists $u\geq 0$ satisfying  (\ref{oscillante}) as a distribution in $\Omega$. By local elliptic regularity, $u \in C^\infty(\Omega)$.
Set 
$$
v(r)=\int_{S^{N-1}}u(r\omega)d\omega.
$$
Since $u \ge 0$, then $v \ge 0$ and,  by the divergence theorem, we have for $r > \delta >0$
\begin{align*}
v'(r)=&\int_{S^{N-1}}\nabla u(r\omega)\cdot \omega\, d\omega=r^{1-N}\int_{|\eta|=r}\nabla u(\eta) \cdot \frac{\eta}{r} \, d\eta=r^{1-N}\int_{B_r\setminus B_\delta}\Delta u(x)\, dx\\
&+r^{1-N}\int_{|\eta|=\delta}\nabla u(\eta) \cdot \frac{\eta }{\delta}\, d\eta 
\end{align*}
hence 
$$
\frac{d}{dr} \left (r^{N-1}v'(r)\right )=\int_{|\eta|=r} \Delta u(\eta)\, d\eta=r^{N-1}\int_{S^{N-1}}\Delta u(r\omega)\, d\omega
$$
and therefore
$$
v''+ (N-1+c)\frac{v'}{r}-\frac{b}{r^2}v=\int_{S^{N-1}} \left (\Delta u(r \omega)+\frac{c}{r}\nabla u(r\omega)\cdot r\omega -\frac{b}{r^2}u(r\omega)\right )\, d\omega=\int_{S^{N-1}}r^{-\alpha }Lu(r\omega)\, d\omega
$$
Then 
%%%%%%==red letters==%%%%%%%%%%%%%%%%%%%%%%%%%%%%
it follows from 
(\ref{oscillante}) %%%%%==(18)==%%%%%%%%%%%%%%%%%
that 
$v$ satisfies
$$\lambda v- r^\alpha\left[v''+ (N-1+c)\frac{v'}{r}-\frac{b}{r^2}v\right]=\phi(r).$$
Setting 
$w(s)=e^{(\frac{N-2+c}{2})s}v(e^s)$ we get
%$$z''(r)=\lambda r^{-\alpha}z(r)+\left[b+\left(\frac{N-1+c}{2}\right)\left(\frac{N-3+c}{2}\right)\right]\frac{z(r)}{r^2}-r^{-\alpha}\phi(r)$$
%and defining $w(s)=e^{-\frac{s}{2}}z(e^s)$ the equation above becomes
\begin{equation}  \label{eq-osc}
w''(s)=(k+\lambda e^{(2-\alpha)s})w(s)-e^{(\frac{3}{2}-\alpha) s}\phi (e^s),\ \ s\in \R
\end{equation}
where 
$$ 
k=%\frac{1}{4}+ b+\left(\frac{N-1+c}{2}\right)\left(\frac{N-3+c}{2}\right)=
b+\frac{(N-2+c)^2}{4} <0.
$$
We choose $m \in \R$ such that $(k+\lambda e^{(2-\alpha)s}) \le k/2 <0$ for $s \le m$. 
By the Sturm Comparison Theorem all non-zero solutions of the homogeneous equation 
\begin{equation} \label{eq-osc1}
 \zeta''(s)=(k+\lambda e^{(2-\alpha)s})\zeta(s)
\end{equation} are oscillating for $s \le m$.
By variation of parameters we write 
$$w(s)=u_2(s)\int_{-\infty}^su_1(t)g(t)dt + u_1(s)\int_{s}^\infty u_2(t)g(t)dt+c_1u_1(s)+c_2u_2(s),$$
where 
$c_1, c_2\in\C$, $g(s)=e^{(\frac{3}{2}-\alpha) s}\phi (e^s)$ and 
$u_i$, $i=1,2$ are linearly independent solutions of (\ref{eq-osc1}) with Wronskian equal to 1.
Since $g$ is compactly supported  we have  for $s$ near $-\infty$
$$w(s)= u_1(s)\int_{\supp g} u_2(t)g(t)dt+c_1u_1(s)+c_2u_2(s).$$
However $w$ is non-negative, because $v \ge 0$, and also oscillating near $-\infty$ since solves (\ref{eq-osc1}). Hence $w=0$ near $-\infty$ and therefore 
$$c_1=-\int_{\supp g} u_2(t)g(t)dt, \quad c_2=0.$$
%%%%%%%%%%%%%%%%%%%%%==page 6==%%%%%%%%%%%%%%%%%%
This gives  
\begin{align*}
w(s)&=u_2(s)\int_{-\infty}^su_1(t)g(t)dt + u_1(s)\int_{s}^\infty u_2(t)g(t)dt-u_1(s)\int_{\supp g} u_2(t)g(t)dt\\
&=u_2(s)\int_{-\infty}^su_1(t)g(t)dt -u_1(s)\int_{-^\infty}^s u_2(t)g(t)dt=
\int_{-\infty}^s (u_1(t)u_2(s)- u_1(s)u_2(t))g(t)dt.
\end{align*}
For fixed $s$ the function $t \mapsto G(s,t)=u_1(t)u_2(s)- u_1(s)u_2(t)$ is also oscillating near $t=-\infty$. Therefore, 
if we choose $g \neq 0$ such that $G(s,t)<0$ on $\supp g$, we get $w(s)<0$ and this contradicts $v\geq 0$.
The case $\alpha >2$ is similar arguing near $+\infty$ instead of $0$.
\qed

%%%%% Section 4 %%%%%
\section{The case  $b+\left(\frac{N-2+c}{2}\right)^2>0$: $L_{min}$ and $L_{max}$}
We always assume that $\alpha \neq 2$.
We recall the function $f$ defined in Introduction 
$$f(s)=b+s(N-2+c-s)$$ and note that
$$\max_{s\in\R}f(s)=b+\left(\frac{N-2+c}{2}\right)^2>0$$
Its roots are $s_1, s_2$ defined in (\ref{defs})
and $f(s)>0$ if and only if $s_1 <s<s_2$. Observe that the equation $Lu=0$ has the two radial solutions $|x|^{-s_1}, |x|^{-s_2}$.

\begin{defi} \label{Leps} In order to approximate $L$ with uniformly elliptic operators we set for 
%%%%%%==red letters==%%%%%%%%%%%%%%%%%%%%%%%%%%%%
 $\eps$ $(0<\eps<1<\eps^{-1})$,
%$\eps >0$, 
$\Omega_\eps=B_{\eps^{-1}} \setminus B_\eps$ and $L_\eps=L$ in $\Omega_\eps$ with Dirichlet boundary conditions. Since $L_\eps$ is uniformly elliptic it follows that $D(L_\eps)=$ 
%%%%%%==red letters==%%%%%%%%%%%%%%%%%%%%%%%%%%%%
$W^{2,p}(\Omega_{\eps}) \cap W^{1,p}_0 (\Omega_{\eps})$.
%W^{2,p}(D_{\eps}) \cap W^{1,p}_0 (D_{\eps})$. 
To shorten the notation we also write $L_0$ for $L_{min}$.
\end{defi}

\subsection{Positive results for $L_{min}$}

We first prove necessary and sufficient conditions under which $L_\eps$ and $L_{min}$ are sectorial 
%%%%%%==red letters==%%%%%%%%%%%%%%%%%%%%%%%%%%%%
in the sense of 
\cite[Definitions 1.5.8]{Goldstein}; 
note that the sectoriality (or more precisely, 
sectorial-valuedness) in a Hilbert space 
was originally introduced in 
\cite[Section V.3.10]{Kato}.

\begin{prop}  \label{gen-pot}
Let $1<p<\infty$. 
If $$f\left(\frac{N-2+\alpha}{p}\right)=b+\frac{(N+\alpha-2)^2}{p'p}+\frac{(N+\alpha-2)(c-\alpha)}{p}>0,$$ or equivalently, 
\begin{equation*} 
%b> -\frac{(N+\alpha-2)^2}{p'p}-\frac{(N+\alpha-2)(c-\alpha)}{p},
s_1+\frac{2-\alpha}{p}<\frac{N}{p}<s_2+\frac{2-\alpha}{p},
\end{equation*}
then the  operators $L_\eps$, $L_{min}$ are sectorial in $L^p(\Omega_\eps)$, $L^p(\R^N)$, respectively, with sectoriality constants independent of $\eps$. 
Moreover $L_{min}$ is dissipative in $L^p(\R^N)$ if and only if 
\begin{equation} \label{condition-sec} 
b+\frac{(N+\alpha-2)^2}{p'p}+\frac{(N+\alpha-2)(c-\alpha)}{p} \ge 0.
\end{equation}
Furthermore, if $f(\frac{N-2+\alpha}{p})=0$ and $\frac{N-2+\alpha}{p}=\frac{N-2+c}{2}$, then $L_\eps$ and $L_{min}$ are sectorial.
\end{prop}
The dissipativity of $L_{\min}$ with $\alpha = 0$ 
is independently proved in 
\cite{Oka-Soba}
with constant $f((N-2)/p)$.% 
\medskip

\noindent
{\sc Proof.} Let $u\in D(L_\eps)$ for $\eps >0$ or $u \in C_c^\infty(\Omega)$ when $\eps=0$. 
%Set $u^\star=\ov{u}|u|^{p-2}$. 
Multiply  $Lu$  by $\ov{u}|u|^{p-2}$ and integrate it over $\R^N$. 
The integration by parts is straightforward when $p\geq 2$. 
For $1<p<2$, $|u|^{p-2}$ becomes singular near the zeros of $u$. 
It is possible to prove that the integration by parts is allowed 
also in this case (see \cite{met-spi}).  
Put $v=|x|^{\frac{N-2+\alpha}{p}}u$. Then 
\begin{align}
Lu
&=
|x|^{\alpha}{\rm div}\,
\left(|x|^{-\frac{N+\alpha-2}{p}}\nabla v-\frac{N+\alpha-2}{p}|x|^{-\frac{N+\alpha-2}{p}-2}xv
\right) 
\\
&\quad+
c|x|^{\alpha-\frac{N+\alpha-2}{p}-2}x\cdot\nabla v
-
\left(b+\frac{c(N+\alpha-2)}{p}\right)|x|^{\alpha-\frac{N+\alpha-2}{p}-2}v. 
%\left(|x|^{-(N-2-\frac{N+\alpha-2}{p})}\overline{v}|v|^{p-2}\right)
\end{align}
Setting $v^\star=\overline{v}|v|^{p-2}$, by integration by parts we have 
\begin{align*}
&\int_{\Omega_\eps}(-Lu)\overline{u}|u|^{p-2}\,dx \\
&=
\int_{\Omega_\eps}
|x|^{2-N}\left(\nabla v-\frac{N+\alpha-2}{p}\frac{xv}{|x|^2}\right) 
\cdot\left(\nabla v^\star-\left(N-2-\frac{N+\alpha-2}{p}\right)\frac{xv^\star}{|x|^2}\right) 
\,dx
\\
&\quad-c
\int_{\Omega_\eps}
|x|^{-N}(x\cdot\nabla v) v^\star\, dx
+
\left(b+\frac{c(N+\alpha-2)}{p}\right)\int_{\R^N}|x|^{-N}vv^\star\,dx
\\
&=
\int_{\Omega_\eps}|x|^{2-N}\left(\nabla v\cdot\nabla v^*\right)\,dx 
-
2\left(\frac{N-2+c}{2}-\frac{N+\alpha-2}{p}\right)
\int_{\Omega_\eps}|x|^{-N}(x\cdot\nabla v) v^*\, dx
\\
&\quad+
\left[b+\frac{N+\alpha-2}{p}\left(N-2+c-\frac{N+\alpha-2}{p}\right)\right]
\int_{\Omega_\eps}|x|^{-N}vv^*\,dx.
\end{align*}
By taking real and imaginary parts of both sides of the equality,
and since ${\rm div } ( x|x|^{-N})=0$ we have
\begin{align} \label{zeroreal}
 \nonumber &Re\left ( \int_{\Omega_\eps}(-Lu)\ov{u}|u|^{p-2}\,dx\right) \\
 \nonumber &=(p-1)\int_{\Omega_\eps}|x|^{2-N}|u|^{p-4}|Re(\ov{v}\nabla v)|^2\,dx
+
\int_{\Omega_\eps}|x|^{2-N}|v|^{p-4}|Im(\ov{v}\nabla v)|^2\,dx\\
&-2\left(\frac{N-2+c}{2}-\frac{N+\alpha-2}{p}\right)
\int_{\Omega_\eps}|x|^{-N}(x\cdot Re(\ov{v}\nabla v))|v|^{p-2}\, dx \\
\nonumber &+
f\left(\frac{N+\alpha-2}{p}\right)\int_{\Omega_\eps}|x|^{-N}|v|^p\,dx
\\
\nonumber &=(p-1)\int_{\Omega_\eps}|x|^{2-N}|u|^{p-4}|Re(\ov{v}\nabla v)|^2\,dx
+
\int_{\Omega_\eps}|x|^{2-N}|v|^{p-4}|Im(\ov{v}\nabla v)|^2\,dx
\\
\nonumber&\quad+
f\left(\frac{N+\alpha-2}{p}\right)\int_{\Omega_\eps}|x|^{-N}|v|^p\,dx \\
\nonumber &=\int_{\Omega_\eps}|x|^{2-N}(\nabla v\cdot \nabla v^*)\, dx+f\left(\frac{N+\alpha-2}{p}\right)\int_{\Omega_\eps}|x|^{-N}|v|^p\,dx,
\end{align}
\begin{align*}
Im\left (\int_{\Omega_\eps}(-Lu)\ov{u}|u|^{p-2}\,dx\right )
&=(p-2)\int_{\Omega_\eps}|x|^{2-N}|u|^{p-4}Re(\ov{u}\nabla u) \cdot Im(\ov{u}\nabla u)\,dx
\\
&\quad 
-2\left(\frac{N-2+c}{2}-\frac{N+\alpha-2}{p}\right)
\int_{\Omega_\eps}|x|^{-N}(x\cdot Im(\ov{v}\nabla v))|v|^{p-2}\,dx.
\end{align*}
Therefore setting
\begin{align*}
B^2&=\int_{\Omega_\eps}|v|^{p-4}|x|^{2-N}|Re(\ov{v}\nabla v)|^2 dx,\\
C^2&=\int_{\Omega_\eps}|v|^{p-4}|x|^{2-N}|Im(\ov{v}\nabla v)|^2 dx,\\
D^2&=\int_{\Omega_\eps}|x|^{(\alpha-2)}|u|^p dx=\int_{\Omega_\eps}|x|^{-N}|v|^p dx,
\end{align*}
we see that
\begin{align}  \label{reale}
Re\bigg(\int_{\Omega\eps}(-Lu)\ov{u}|u|^{p-2}\,dx\bigg)
= (p-1)B^2+C^2+f\left(\frac{N+\alpha-2}{p}\right)D^2
%\left(p-1+\frac{(c-\alpha)
%p}{N-2+\alpha}\right)\int_{\Omega_\eps}|x|^\alpha|u|^{p-4}|Re(\ov{u}\nabla
%u)|^2dx\\\nonumber&+
%\int_{\Omega_\eps}|x|^\alpha|u|^{p-4}|Im(\ov{u}\nabla u)|^2dx
\end{align}
%%%%%%%%%%%%%%%%%%%%%==page 8==%%%%%%%%%%%%%%%%%%
and
\begin{align}%\label{immaginaria}
\nonumber
&\bigg|Im\bigg(\int_{\Omega_\eps}(-Lu)\ov{u}|u|^{p-2}\,dx\bigg)\bigg|
\\ \nonumber
&\leq|p-2|\left(\int_{\Omega_\eps}|v|^{p-4}|x|^{2-N}|Re(\ov{v}\nabla v)|^2\,dx\right)^\frac{1}{2}
\left(\int_{\Omega_\eps}|v|^{p-4}|x|^{2-N}|Im(\ov{v}\nabla v)|^2\, dx\right)^\frac{1}{2}
\\ \nonumber
&
+2\left|\frac{N-2+c}{2}-\frac{N+\alpha-2}{p}\right|
\int_{\Omega_\eps}|v|^{p-2}|x|^{1-N}|Im(\ov{v}\nabla v)|\,dx
\\ \nonumber& \leq|p-2|BC+2\left|\frac{N-2+c}{2}-\frac{N+\alpha-2}{p}\right|CD.\nonumber
\end{align}
By condition (\ref{condition-sec}) or condition $\frac{N-2+\alpha}{p}=\frac{N-2+c}{2}$, 
we see that %there exists a positive constant $l_\alpha$  such that
%$$
%     (p-1)B^2+C^2 + 
%f\left(\frac{N+\alpha-2}{p}\right)D^2 
%\geq l_\alpha\left(|p-2|BC+2\left|\frac{N-2+c}{2}-\frac{N+\alpha-2}{p}\right|CD\right)
%$$
%and, consequently,
\begin{equation} \label{sectconstants}
\left|Im\bigg(\int_{\Omega_\eps}(-Lu)u|u|^{p-2}\,dx\bigg)\right|
\leq l_\alpha\left\{Re\bigg(\int_{\Omega_\eps}(-Lu)u|u|^{p-2}\,dx\bigg)\right\},
\end{equation} 
where 
\[l_\alpha=\sqrt{\frac{(p-2)^2}{4(p-1)}+
\left|\frac{N-2+c}{2}-\frac{N+\alpha-2}{p}\right|^2f\left(\frac{N+\alpha-2}{p}\right)^{-1}}
\]
($0/0=0)$.
This shows the sectoriality of $L_\eps$ and $L_{min}$, with sectoriality constants independent of $\eps$.
%If $\tan\theta_\alpha=l_\alpha$, then $e^{\pm i\theta}L$ is
%dissipative in $\R^N$ for $0 \le \theta \le \theta_\alpha$. 
Assume now that $L_{min}$ is dissipative. Then, by (\ref{reale}) for real-valued functions, the inequality
$$(p-1)\int_{\R^N}|x|^{2-N}|\nabla v|^2|v|^{p-2}\, dx +f\left(\frac{N+\alpha-2}{p}\right)\int_{\R^N}|x|^{-N}|v|^{p}\, dx\geq 0$$ holds. By \cite[Corollary 2.3 (ii)]{costa} (with $b=(N-2)/2$ and $u=|v|^{p/2}$)
we obtain 
$f\left(\frac{N+\alpha-2}{p}\right)\geq 0.$
%$$b\geq -\frac{(N+\alpha-2)^2}{p'p}-\frac{(N+\alpha-2)(c-\alpha)}{p}.$$ 
\qed

\begin{os} \label{gen-pot-crit}
{\rm We remark that the above proposition holds also when $b+(N-2+c)^2/4=0$ and $(N-2+\alpha)/p=(N-2+c)/2$. In this case $s_1=s_2=(N-2+c)/2$, hence the condition $f(\frac{N-2+\alpha}{p})=0$ is satisfied. We also remark that the choice of the power in the substitution $v=|x|^{\frac{N-2+\alpha}{p}}u$ is the only one which leads to the term $x|x|^{-N}$ in (\ref{zeroreal}) which has zero divergence.}
\end{os}

\noindent
We can state the main result of this subsection.

\begin{teo}\label{gen-min}
Let $p,\ \alpha,\ c,\ b$ satisfy 
$$f\left(\frac{N}{p}-2+\alpha\right) 
 = \left(\frac{N}{p}-2+\alpha\right) 
   \left(\frac{N}{p'}-\alpha+c\right)+b >0,$$ 
%%%%%%%%%%%%%%%%==page 9==%%%%%%%%%%%%%%%%%%%%%%%
or equivalently, 
$$s_1+2-\alpha<\frac{N}{p}<s_2+2-\alpha.$$
Then the operator $L_{min}$ generates a bounded positive analytic semigroup in $L^p(\R^N)$, coherent with respect to all $p$ satisfying the above inequalities. Moreover,  $D_{min}(L)$ coincides with 
$$
D_{p,\alpha}=\{u\in W^{2,p}_{loc}(\Omega),\ |x|^\alpha D^2u,\ |x|^{\alpha-1} \nabla u,\ |x|^{\alpha-2}u\in L^p(\R^N)\}.\\
$$
\end{teo}

\noindent
Note that the generation interval $]s_1+2-\alpha, s_2+2-\alpha[$ differs from the conctractivity interval  $]s_1+(2-\alpha)/p, s_2+(2-\alpha)/p[$. In particular, for certain values of $p$, $L_{min}$ generates a non-contractive semigroup. 
\noindent
The proof of the theorem above  is based upon the perturbation result stated in Theorem \ref{perturbation}. 
The next lemma provides the validity of its assumptions for the operators $L_\eps$ introduced in Definition \ref{Leps}, with constants independent of $\eps$.

%%%%%==Lemma 4.4==%%%%%%%%%%%%%%%%%%%%%%%%%%%%%%%
\begin{lem}\label{dual}
Let $p,\ \alpha,\ c,\ b$ as in Theorem \ref{gen-min}. 
%%%%%%==red letters==%%%%%%%%%%%%%%%%%%%%%%%%%%%%
Put 
$$
M := f\left(\frac{N}{p}-2+\alpha\right)
  = b+\left(\frac{N}{p}-2+\alpha\right) 
      \left(\frac{N}{p'}-\alpha+c\right) > 0 
$$
as in Theorem 4.3. Then for 
%There exists a positive constant $M$ independent of 
%$\eps \ge 0$ such that with 
$V(x)=|x|^{\alpha-2}$, 
$$
     -Re\int_{\Omega_\eps}(Lu)V^{p-1}\ov{u}|u|^{p-2}\,dx 
\geq M\|Vu\|_p^p
$$ 
for every $u\in C_c^\infty(\Omega)$. 
%%%%%%==red letters==%%%%%%%%%%%%%%%%%%%%%%%%%%%%
%{\color{red} 
%where 
%
%}
\end{lem}

\noindent
{\sc Proof.} 
%%%%%%==red letters==%%%%%%%%%%%%%%%%%%%%%%%%%%%%
%Let $M_{\delta}(\alpha)$ be the constant as in 
%\eqref{sectoriality}. 
Put $\beta := \alpha+(p-1)(\alpha-2)$. 
Then since $\frac{N-2+\beta}{p}=\frac{N}{p}+\alpha-2$, we have
\[
f\left(\frac{N-2+\beta}{p}\right)=f\left(\frac{N}{p}+\alpha-2\right)>0.
\]
%Then $M_{p-1}(\beta)>0$. 
%In fact, we see that 
%$$
%  M_{p-1}(\beta) = b+\frac{(c-\beta)(N+\beta-2)}{p} 
%                 + \frac{(N+\beta-2)^2}{p\,'p} 
%$$ 
%can be written (in terms of $\alpha$) as $M$ 
%in the statement. 
Now let $u\in D(L_\eps)$ for $\eps >0$ or $u \in C_c^\infty(\Omega)$ when $\eps=0$. Then
\begin{align*}
&\quad-Re\int_{\Omega_\eps}(Lu)V^{p-1}\ov{u}|u|^{p-2}\, dx
\\
&=- Re \int_{\Omega_\eps}\left(|x|^\alpha\Delta u+c|x|^{\alpha-1}\frac{x}{|x|}\nabla u-b|x|^{\alpha-2}u\right)|x|^{(\alpha-2)(p-1)}\ov{u}|u|^{p-2}\, dx
\\
&=- Re \int_{\Omega_\eps} \left(|x|^\beta\Delta u+c|x|^{\beta-1}\frac{x}{|x|}\nabla u-b|x|^{\beta-2}u\right)\ov{u}|u|^{p-2}\, dx, 
\end{align*}
where $\beta=\alpha+(p-1)(\alpha-2)$ 
%%%%%%==red letters==%%%%%%%%%%%%%%%%%%%%%%%%%%%%
is as defined above. 
%%%%%%==red letters==%%%%%%%%%%%%%%%%%%%%%%%%%%%%
By applying 
\eqref{reale} 
with $\alpha=\beta$, 
%By choosing $\delta=p-1$ in  (\ref{sectoriality}) 
%and changing $\alpha$ with $\beta$ 
%%in the definition of $M(\delta) = M(p-1)$, 
we get
$$
   - Re\int_{\Omega_\eps}(Lu)V^{p-1}\ov{u}|u|^{p-2}\,dx 
\geq 
%%%%%%==red letters==%%%%%%%%%%%%%%%%%%%%%%%%%%%%
f\left(\frac{N-2+\beta}{p}\right)
\int_{\Omega_\eps}|x|^{\beta-2}|u|^p\, dx. 
$$
%%%%%%==red letters==%%%%%%%%%%%%%%%%%%%%%%%%%%%%
Since $\beta-2 = p(\alpha-2)$, 
this is nothing but the desired inequality with $M=f(\frac{N}{p}+\alpha-2)$. 
\qed

%%%%%%%%%%%%%%%%%%==page 10==%%%%%%%%%%%%%%%%%%%%

\noindent
{\sc Proof of Theorem \ref{gen-min}}.
{\bf Step 1.} First assume that $b$ is sufficiently large so that the conditions of Proposition \ref{gen-pot} and Lemma \ref{dual} are satisfied. Then, by (\ref{sectconstants}),   there exists $0 <\theta < \pi/2$ such that $\lambda- L_\eps$ is injective for $\lambda \in \Sigma_{\pi/2+\theta}$ for every $\eps \ge 0$ (take such a $\theta$ with $\tan \theta < l_\alpha$). For $\eps>0$, $L_\eps$ is uniformly elliptic, hence generates an analytic semigroup. By (\ref{sectconstants}) again, $\lambda-L_\eps  $ is invertible and satisfies $\|(\lambda-L_\eps)^{-1}\| \le C|\lambda|^{-1}$ for  $\lambda \in \Sigma_{\pi/2+\theta}$ with $C$ independent of $\eps$ (actually, $C=1$ in $\Sigma_\theta$). Let $f \in L^p(\R^N)$ and $u_\eps=(\lambda-L)^{-1}f$ for $\lambda \in \Sigma_{\pi/2+\theta}$. Since $\|u_\eps\| \le C|\lambda|^{-1}$ and $\lambda u_\eps -Lu_\eps=f$, by 
local elliptic regularity we can find a sequence $u_n$ such that $u_{\eps_n} \to u$ weakly  in $W^{2,p}_{loc}(\Omega)$, strongly 
in $L^p_{loc}(\Omega)$ and 
pointwise. By Lemma \ref{dual}, $M\|Vu_\eps\| \le \|Lu_\eps\| \le (C+1)\|f\|$ ($V(x)=|x|^{\alpha-2}$) and therefore $|\lambda|\|u\|_p \le C\|f\|_p$, $\| Vu\|_p \le (C+1)\|f\|_p$ and $\lambda u-Lu=f$. Since $u \in D_{max}(L) \cap D(V)$, by Lemma \ref{aux} $u \in D_{min}(L)$ and this shows that $\lambda-L_{min}$ is invertible for $\lambda \in \Sigma_{\pi/2+\theta}$ and that $\|(\lambda-L_{min})^{-1}\| \le C|\lambda|^{-1}$. If $\lambda>0$ and $f \ge 0$, then $u_\eps \ge 0 $ by the classical maximum principle, hence $u \ge 0$. Finally, if $f \in L^p(\R^N) \cap L^q(\R^N)$, then the solutions $u_\eps$ do not depend on $p,q$ and we can select the same sequence $(u_{\eps_n})$ convergent both in $L^p$ and in $ L^q$. Therefore $u$ is the same in $L^p,L^q$ and this shows the coherence of the resolvents. 

\noindent
{\bf Step 2.} Assume now that $b$ satisfies only the condition in the statement and let $A=-L+kV$ with $k$ large enough to satisfy also the conditions of Step 1. Then  $(-A)_{min}$ generates an analytic semigroup of positive contractions in $L^p(\R^N)$, coherent with respect to $1<p<\infty$. By Lemma \ref{dual}, we have 
$$\int_{\R^N} (-Lu)\, V^{p-1}\ov{u}|u|^{p-2}\, dx\geq M\|Vu\|_p^p.$$
Therefore 
 $$\int_{\R^N} (-Lu+kV)\, V^{p-1}\ov{u}|u|^{p-2}\, dx\geq (M+k)\|Vu\|_p^p.$$
By Theorem \ref{perturbation}, the operator $-(A+tV)=L-kV-tV$, with the same domain as $(-A)_{min}$, generates a bounded analytic semigroup of positive operators for every $t> -(M+k)$. 
In particular, choosing $t=-k$, we deduce that the operator $(L,D_{min}(A))$ generates a bounded analytic semigroup of positive operators in $L^p(\R^N)$. Since $C_c^\infty(\Omega)$ is a core for $A_{min}$ it follows that $(L,D_{min}(A))=L_{min}$.
Finally if also $q$ satisfies the inequalities in the statement, the coherence of the resolvents, hence of the semigroups, follows from the perturbation argument, since the unperturbed semigroups are coherent, by Step 1.

\noindent
{\bf Step 3.} Finally we prove equality (\ref{domain}). The inclusion $D_{p,\alpha} \subset D_{min}(L)$ follows since $C_c^\infty (\Omega)$ is dense in $D_{p,\alpha}$, see Lemma \ref{densita}. Conversely, let $u\in D_{min}(L)$. By \cite[Theorem 3.1]{met-soba-spi}, we have the estimate $\||x|^{\alpha-2}u\|_p \le C\|Lu\|_p$ and therefore $u \in D_{max}(L) \cap D(|x|^{\alpha-2})$. By Lemma \ref{aux}, $u \in D_{p,\alpha}$. \qed

%%%%% Section 4 %%%%%
\subsection{Positive results for $L_{max}$}
We consider the adjoint operator  
$$\tilde{L}=|x|^\alpha\Delta +\tilde{c}|x|^{\alpha-1}\frac{x}{|x|}\cdot\nabla -\tilde{b}|x|^{\alpha-2}$$ 
where $\tilde{c}= 2\alpha-c$, $ \tilde{b}=b+(c-\alpha)(\alpha-2+N)$
on $L^{p'}(\R^N)$, see (\ref{defLtilde}), (\ref{tilde}) .
Since, by Proposition \ref{adjoint}, $L_{max}=(\tilde L_{min})^*$, we deduce generation results for $L_{max}$ by duality.

\begin{teo}  \label{gen-max}
Let $p,\ \alpha,\ c,\ b$ such that
$$f\left(\frac{N}{p}\right)
%%%%%%==red letters==%%%%%%%%%%%%%%%%%%%%%%%%%%%%
=b+\omega_{p} 
=b+\frac{N}{p}\left(\frac{N}{p'}-2+c\right)
>0,$$ 
or equivalently,
$$s_1<\frac{N}{p}<s_2.$$
Then the operator $L_{max}$ generates a bounded positive analytic semigroup in $L^p(\R^N)$, coherent with respect to all $p$ satisfying the above inequalities. 
\end{teo}
{\sc Proof.}
It is sufficient to write the conditions  of Theorem \ref{gen-min}  for the operator $\tilde{L}$ in $L^{p'}(\R^N)$, to recall that $\tilde{s}_i=s_i+\alpha-c$, $i=1,2$, $s_1+s_2=N-2+c$,  and then to argue by duality.
\qed
\noindent
Observe that the condition in the above theorem is independent of $\alpha$. Observe also that if $p$ satisfies both the conditions of Theorems \ref{gen-min}, \ref{gen-max}, that is if $s_1<N/p<s_2$ and $s_1+2-\alpha<N/p <s_2 +2-\alpha$, then $L_{min}=L_{max}$.

\begin{os} \label{endpoints}
{\rm
In the next sections we shall see what happens in Theorems \ref{gen-min}, \ref{gen-max} when $N/p$ coincides with one of the endpoints. For example, if  $\alpha <2$, then $L_{min}$ generates if and only if $s_1+2-\alpha \le N/p <s_2 +2-\alpha$, see Propositions \ref{counter1}, \ref{nonex}, \ref{endpoint-Lmin}  but for the  equality $D(L_{min})=D_p$ one needs  $s_1+2-\alpha < N/p <s_2 +2-\alpha$, see Proposition \ref{lasttheta}. By duality $L_{max}$ is a generator if and only if $s_1 <N/p \le s_2$. The case $\alpha >2$ is similar with the roles of $s_1, s_2$ interchanged. 
}
\end{os}

%%%%%%%%%%%%%%%%%%%%%%%%%==page 12==%%%%%%%%%%%%%
%%%%%%%%%%%%==Section 4.3==%%%%%%%%%%%%%%%%%%%%%%
\subsection{Negative results for $L_{min}$}
We  prove that the generation conditions for $L_{min}$ given in Theorem \ref{gen-min} are sharp. 

\begin{prop} \label{counter1}
If $\alpha <2$ and  $N/p<s_1+2-\alpha$ or $\alpha>2$ and $N/p>s_2+2-\alpha$  then, for every $\lambda>0$, $\ov{Rg(\lambda-L_{min})}\neq L^p(\R^N)$. Therefore $L_{min}$ does not generate a semigroup in $L^p(\R^N)$.
\end{prop}
{\sc Proof.} We focus on the case $\alpha<2$, the other being similar. We consider the adjoint operator $\tilde L$ defined in (\ref{defLtilde}), see  Proposition \ref{adjoint},  and we prove that 
$N(\lambda-\tilde{L}_{max})\neq\{0\}$ in $L^{p'}(\R^N)$ 
by exhibiting a radial function $u\in D_{max}(\tilde{L})$ in $L^{p'}(\R^N)$ satisfying
\begin{equation}  \label{eigenf}
\lambda u-\tilde{L}u=0.
\end{equation}
By Lemma \ref{bhv} with $b$ and $c$ respectively replaced with $\tilde{b}$ and $\tilde{c}$ defined in (\ref{tilde}), $u$ can be written 
by $u=c_1u_1+c_2u_2$, where $u_j$ is defined in Lemma \ref{bhv}. 
In order to have integrability of $u(\rho)$ in $L^{p'}(\R^N)$ for large $\rho$, 
we consider the solution $u_2$ ($c_1=0$ and $c_2=1$). 
This choice will lead to an additional assumption to insure also the integrability 
near the origin: the solution $u_2$ is in $L^{p'}(B_1)$ if and only if
\begin{equation}\label{no-sg}\left(\frac{N-2+2\alpha-c}{2}+\sqrt k\right)<\frac{N}{p'}, \quad {\rm or\  equivalently, \quad }
\frac{N}{p}<s_1+2-\alpha.
\end{equation}
\qed

%%%%%%%%%%%%==Proposition 4.8==%%%%%%%%%%%%%%%%%%
\begin{prop} \label{nonex}
If $\alpha <2$ and  $N/p \ge s_2+2-\alpha$ or $\alpha>2$ and $N/p \le s_1+2-\alpha$ then
for every $\lambda>0$, $N(\lambda-L_{min})\neq \{0\}$. 
Therefore no extension of $L_{min}$ generates a semigroup in $L^p(\R^N)$.
\end{prop}
{\sc Proof. } As in the proof of Proposition \ref{counter1}, we focus on the case $\alpha<2$ and we prove the existence of a 
radial function $u\in D_{min}(L)\setminus\{0\}$ satisfying 
$$
u-Lu%=u-\left[\rho^\alpha(u''+\frac{N-1+c}{\rho}u')-b\rho^{\alpha-2}u\right]
=0.
$$
%By the same computation as in Proposition \ref{counter1} with $\tilde{c}$ and $\tilde{b}$ 
%replaced with $c$ and $b$, the above equation has two independent solutions ($k=b+ (N-2+c)^2/4>0 $)
%\[
%u_{1}\approx k^{-\frac{1}{4}}\rho^{-s_{2}}, \quad
%u_{2}\approx k^{-\frac{1}{4}}\rho^{-s_{1}} \quad{\rm for\ small\ }\rho, 
%\]
%and 
%\[
%u_1 \approx \rho^{-\frac{N-1+c}{2}+\frac{\alpha}{4}}\exp\left\{-\frac{2}{2-\alpha}\rho^\frac{2-\alpha}{2}\right\},
%\quad
%u_2 \approx \rho^{-\frac{N-1+c}{2}+\frac{\alpha}{4}}\exp\left\{\frac{2}{2-\alpha}\rho^\frac{2-\alpha}{2}\right\}
%\quad{\rm for\ large\ }\rho.
%\]
We write $u=c_1u_1+c_2u_2$, where $u_j$ is defined in Lemma \ref{bhv} for $j=1, 2$. 
The integrability of $u$ near $\infty$ implies that $u=c_2u_2$ with $c_2\neq 0$.

\noindent
We prove that $u_2\in D_{min}(L)$. We first assume that $s_2<N/p+\alpha-2$. 
In this case from \eqref{u2-behave} we have $u_1, |x|^{\alpha-2}u_1\in L^p(\R^N)$. By Lemma \ref{aux}, 
we obtain $u_2\in D_{min}(L)$. 

\noindent
Next, we assume that $s_2=N/p+\alpha-2$. 
Let $\eps >0$ with $\alpha+\eps<2$. Then using \eqref{u2-behave}, we have 
\begin{align*}
&\||x|^{\alpha-2+\eps}u_2\|^p= 
\int_{B_1}|x|^{(\alpha-2+\eps)p}|u_2|^{p}\,dx+\int_{\R^N\setminus B_1}|x|^{(\alpha-2+\eps)p}|u_2|^{p}\,dx\\
&\leq 
C_1\int_{B_1}|x|^{(\alpha-2+\eps)p-ps_2}\,dx+
C_2\int_{\R^N\setminus B_1}
\exp\left\{-C_3 |x|^\frac{2-\alpha}{2}\right\}
\,dx
\\
&\leq 
C_1 \omega_N\int_0^1r^{\eps p-1}\,dr+C_2'
\leq 
\frac{C_1\omega_N}{\eps p}+C_2'.
\end{align*}
We apply Lemma \ref{aux}  to $|x|^\eps L$ (with $\alpha+\eps$ instead of $\alpha$) to deduce that $u_2 \in D_{p, \alpha+\eps}$. Moreover the interpolation inequality (\ref{int}) yields
\begin{align*}
\||x|^{\alpha-1+\eps}\nabla u_2\|
\leq 
C'(\|(|x|^{\eps}L) u_2\|+\||x|^{\alpha-2+\eps}u_2\|)
\leq 
C'(\||x|^{\eps}u_2\|+\||x|^{\alpha-2+\eps}u_2\|)
\leq 
C''(1+\eps^{-\frac{1}{p}}).
\end{align*}
%%%%%%%%%%%%%%%%==page 14==%%%%%%%%%%%%%%%%%%%%%%
Hence we have
\[
L(|x|^{\eps}u_2)=|x|^{\eps}u_2+2\eps |x|^{\alpha-2+\eps}x\cdot\nabla u_2+
\eps(N-2+c+\eps)|x|^{\alpha-2}u_2\in L^p(\R^N).
\]
This implies that $|x|^{\eps}u_2\in D_{max}(L)\cap D(|x|^{\alpha-2})\subset |x|^{\eps}u_2\in D_{min}(L)$, by Lemma \ref{aux}. Moreover, by the above estimates
\[
\|L_{min}(|x|^\eps u_2)-u_2\|\leq C'''(\eps+\eps^{1-\frac{1}{p}}).
\]
The closedness of $L_{min}$ yields $u_2\in D_{min}(L)$ and $L_{min}u_2=u_2$.
\qed

\subsection{Negative results for $L_{min} \subset L \subset L_{max}$}

Since $L_{max}=(\tilde{L}_{min})^*$, from Propositions \ref{counter1}, \ref{nonex}
we obtain the following result.

%%%%%%%%%%%%==Proposition 4.9==%%%%%%%%%%%%%%%%%%
\begin{prop} \label{nonex2}
\begin{itemize}
\item[(i)] If $\alpha<2$ and $N/p\leq s_1 $ or  $\alpha>2$ and $N/p\geq s_2$, then
for every $\lambda>0$, $\ov{R(\lambda-L_{max})}\neq L^p(\R^N)$. 
Therefore no restriction of $L_{max}$ generates a semigroup in $L^p(\R^N)$. 
\item[(ii)] If $\alpha<2$ and $N/p >s_2 $ or  $\alpha>2$ and $N/p < s_1$, then
for every $\lambda>0$, $\ov{N(\lambda-L_{max})}\neq L^p(\R^N)$. 
Therefore  $L_{max}$ does not generate a semigroup in $L^p(\R^N)$.
\end{itemize}
\end{prop}
{\sc Proof.} We prove (i) and assume $\alpha<2$. By Proposition \ref{nonex} $N(\lambda-\tilde{L}_{min})\neq \{0\}$  in $L^{p'}(\R^N)$, that is $\ov{R(\lambda-L_{max})}\neq L^p(\R^N)$,  if 
$N/p' \ge \tilde{s}_2 +2-\alpha$. Here $\tilde{s_i}$ are the roots of the function $f$ defined in (\ref{deff}) and relative to $\tilde L$. Since $\tilde{s}_i=s_i+\alpha-c$ the above condition reads $N/p \le N+c-2-s_2=s_1$. The case $\alpha>2$ is similar. The proof of (ii) follows similarly from Proposition \ref{counter1}.
\qed

\noindent
Finally we state the following negative result for any realization $(L,D)$ such that $D_{min}(L) \subset D_{max}(L)$, $L_{min} \subset L \subset L_{max}$, in short. This proves the "only if " part in Theorem \ref{main}.

\begin{teo} \label{nosemi}
Let $\alpha \neq 2$ and assume that $N/p \le s_1+\min\{0,2-\alpha\}$ or $N/p \ge s_2+\max\{0,2-\alpha\}$. Then no realization of the operator $L$ between $L_{min}$ and $L_{max}$ generates a semigroup in $L^p(\R^N)$.
\end{teo}
{\sc Proof.} This follows immediately from Propositions \ref{nonex}, \ref{nonex2} (i). In fact, if $\alpha<2$ and $N/p \le s_1$ no restriction of $\lambda -L_{max}$ can be surjective whereas if $N/p \ge s_2+2-\alpha$ no extension of $\lambda-L_{min}$ can be injective.
\qed

%%%%%==Section 5==%%%%%
\section{The case  $b+\left(\frac{N-2+c}{2}\right)^2>0$: $L_{min} \subset L \subset  L_{max}$}
We always assume $\alpha \neq 2$ and show that a suitable realization of 
$L_{min} \subset L \subset  L_{max}$ generates a semigroup in $L^p(\R^N)$ if and only if 

\begin{equation}\label{conjecture}
s_1+\min\{0,2-\alpha\}<\frac{N}{p}<s_2+\max\{0,2-\alpha\}.
\end{equation}
To explain the meaning of the above condition let us fix $\alpha<2$. By the results of the previous section $L_{max}$ generates if $s_1<N/p <s_2$ and $L_{min}$ when $s_1 +2-\alpha <N/p <s_2+2-\alpha$ and $L_{min}=L_{max}$ if both conditions are satisfied. Therefore we have generation under (\ref{conjecture}) if $s_2 <s_1+2-\alpha$. However this last condition is not always verified: this is the case when  when $\alpha$ is very negative but also for $N=3,4$ and $\alpha=b=c=0$: as already pointed out in the Introduction 
%%%%%%%%%%%%%%%%==page 15==%%%%%%%%%%%%%%%%%%%%%%
$\Delta_{min}$ generates for $p \le N/2$ and $\Delta_{max}$ for $p \ge N/(N-2)$.
We also remark that under under the condition $s_1+(2-\alpha)/p  \le N/p \le s_2+(2-\alpha)/p$, see  Proposition \ref{gen-pot},  $L$ is dissipative  in the annulus $B_{\eps{-1}} \setminus B_\eps$ when endowed with Dirichlet boundary conditions.  A semigroup can therefore be constructed via approximation as in Step 1 of the proof of Theorem \ref{gen-min}. We do not follow this approach since it does not cover all cases considered in (\ref{conjecture}).

\subsection{The operator $L_{int}$}

We define an intermediate operator $L_{int}$ between $L_{min}$ and 
$L_{max}$ as 
\begin{equation} \label{Lint}
D_{int}(L)=\{u\in D_{max}(L)\subset L^p(\R^N)\;;\; |x|^{\theta(\alpha-2)}u\in L^p(\R^N)\, {\rm for\ every\ } \theta \in I\}, 
\end{equation}
where $I$ is the interval of all $\theta\in [0,1]$ such that 
\begin{equation}\label{conjecture-other}
f\left(\frac{N}{p}+\theta(\alpha-2)\right)>0.
\end{equation}
Note that \eqref{conjecture} is equivalent to the existence of some $\theta \in [0,1]$
satisfying \eqref{conjecture-other}. First we show the injectivity of $\lambda-L_{int}$ for ${\rm Re}\, \lambda>0$.

\begin{lem}\label{inj-int}For every $\lambda$ such that ${\rm Re}\, \lambda>0$, $\lambda-L_{int}$ is injective.
\end{lem}
{\sc Proof.} We fix $\theta \in I$, ${\rm Re}\, \lambda>0$ and
suppose that $u\in D_{int}(L)$ satisfies $(\lambda u-L_{int})u=0$. 
Set 
\begin{align*}
A
=|x|^{\theta(\alpha-2)}L|x|^{-\theta(\alpha-2)}
=|x|^{\alpha}\Delta +c_A|x|^{\alpha-2}x\cdot\nabla 
-b_A|x|^{\alpha-2}, 
\end{align*} 
where $c_A=c-2\theta(\alpha-2)$ and $b_A:=b+\theta(\alpha-2)(N-2+c-\theta(\alpha-2))$. Then 
\begin{align*}
(\lambda-A)|x|^{\theta(\alpha-2)}u=|x|^{\theta(\alpha-2)}(\lambda-L)u=0.
\end{align*} 
Setting $v:=|x|^{\theta(\alpha-2)}u\in L^p(\R^N)$, 
we have $v\in D(A_{max})$ and $(\lambda -A_{max})v=0$. On the other hand, 
 $A$ satisfies the hypothesis of Theorem \ref{gen-max}
\begin{align*}
&b_A+\frac{N}{p}\left(\frac{N}{p'}-2+c_A\right)
=b+\theta(\alpha-2)(N-2+c-\theta(\alpha-2))
+\frac{N}{p}\left(\frac{N}{p'}-2+c-2\theta(\alpha-2)\right)
\\
&=b+\left(\frac{N}{p}+\theta(\alpha-2)\right)\left(\frac{N}{p'}-2+c-\theta(\alpha-2)\right)
=f\left(\frac{N}{p}+\theta(\alpha-2)\right)
>0.
\end{align*}
Then $A_{max}$ generates a bounded analytic semigroup on $L^p(\R^N)$ and, in particular, 
$\lambda\in \rho(A_{max})$.  Hence $v=(\lambda -A_{max})^{-1}(\lambda-A_{max})v=0$ and,
by the definition of $v$, $u=0$, too. \qed

\noindent
We approximate $L_{int}$ through the  operators
\[
\begin{cases}
L_t := |x|^{\alpha}\Delta + c|x|^{\alpha-2}x\cdot\nabla -(b+k)|x|^{\alpha-2}
+k \min\{t, |x|^{\alpha-2}\},
\\
%%%%%%%%%%%%==red color==%%%%%%%%%%%%%%%%%%%%%%%%
D(L_{t}):=D_{p,\alpha}
\end{cases}
\]
where $t >0$, $D_{p,\alpha}$ is defined in (\ref{domain}) and  $k$ is a large fixed nonnegative constant for which the conditions of Proposition \ref{gen-pot} and Theorems \ref{gen-min}, \ref{gen-max} are satisfied for every $p >1$.
Observe that, in particular $D_{p,\alpha}=D_{min}(L)=D_{max}(L)$.
Then we have 
\begin{lem} \label{g0}
For every $1<p<\infty$, $L_t$ generates an analytic semigroup of positive operators in $L^p(\R^N)$, coherent with respect to $p$. Moreover
$(kt,\infty)\subset \rho(L_t)$ and $C_c^\infty(\Omega)$ is a core for $L_t$. 
\end{lem}
{\sc Proof.} 
Because of  the assumption on $k$ we see from Theorem \ref{gen-min} that 
$L-k|x|^{\alpha-2}$ with domain $D_{p}(L)$ 
generates an analytic semigroup of positive operators in $L^p(\R^N)$ for every $1<p<\infty$, coherent with respect to $p$. 
Since $k\min\{t, |x|^{\alpha-2}\}$ is bounded
the same is true for $L_t$ and moreover
$(kt,\infty)\in \rho(L_t)$. 
\qed

\noindent
We show weighted and unweighted resolvent  estimates for $L_t$ with constants independent of $t$.
\begin{lem}\label{re-int}
Let $\theta$ satisfy \eqref{conjecture-other}. 
Then there exist constants $C,C'>0$ such that for every $\lambda\in \C_+$, $t>0$  
and $u\in C_c^\infty(\Omega)$, 
\begin{gather}
\label{u-Ln}
\|u\|_p\leq \frac{C}{|\lambda|}\|\lambda u-L_tu\|_p.
\\
\label{xu-Ln}
\||x|^{\theta(\alpha-2)}u\|_p\leq \frac{C'}{|\lambda|^{1-\theta}}\|\lambda u-L_tu\|_p.
\end{gather}
Therefore $\C_+\subset \rho (L_t)$ and 
$
\|(\lambda -L_t)^{-1}\|\leq C|\lambda|^{-1}.
$
\end{lem}
{\sc Proof.} 
First we prove \eqref{xu-Ln} when $\theta\in [\frac{1}{p},1]$. We observe that the assumptions of Proposition  \ref{gen-pot} are satisfied if we replace $\alpha$ with $\beta=\alpha+(p\theta-1)(\alpha-2)$, hence $\beta-2=p\theta(\alpha-2)$. 
Therefore we consider the operator   $ |x|^{(p\theta-1)(\alpha-2)} L$. For  $u\in C_c^\infty (\Omega)$
we have 
\[
%\left|{\rm Im}
%\int_{\R^N}
%  (-L_{n}u)|x|^{(p\theta-1)(\alpha-2)}\overline{u}|u|^{p-2}
%\,dx
%\right|
%\leq
{\rm Re}\,\left[
e^{i\omega}\int_{\R^N}
  (-L u)|x|^{(p\theta-1)(\alpha-2)}\overline{u}|u|^{p-2}
\,dx\right]\geq 0, \quad 
\omega\in[-\frac{\pi}{2}+\omega_1, \frac{\pi}{2}-\omega_1], 
\]
where $\pi/2 -\omega_1>0$ is the angle of sectoriality of $ |x|^{(p\theta-1)(\alpha-2)}L$. Since $L_t=L-V_t$ with $V_t \ge 0$, the same inequality holds for $L_t$, thus for ${\rm  Re}\, \lambda >0$
\begin{equation} \label{g1}
{\rm Re \, }(\lambda e^{i \omega}) \int_{\R^N} |x|^{(p\theta-1)(\alpha-2)}|u|^{p} \,dx \le {\rm Re}\,\left[
e^{i\omega}\int_{\R^N}
  (\lambda-L_tu)|x|^{(p\theta-1)(\alpha-2)}\overline{u}|u|^{p-2}
\,dx\right]
\end{equation}
and, by choosing $\omega\in[-\frac{\pi}{2}+\omega_1, \frac{\pi}{2}-\omega_1]$ such that ${\rm Re \, }(\lambda e^{i \omega})=|\lambda| \cos \omega_1$, H\"older inequality yields
\begin{align} \label{weighted}
|\lambda|\cos\omega_1
\left\||x|^{(\theta-\frac{1}{p})(\alpha-2)}u\right\|_p^p
\leq \left\||x|^{\frac{p\theta-1}{p-1}(\alpha-2)}u\right\|_p^{p-1} 
\|\lambda u-L_tu\|_p.
\end{align}
Noting that 
\[
0 \le \theta-\frac{1}{p}<\frac{p\theta-1}{p-1}\leq \theta, 
\]
we apply Lemma \ref{intgamma} with respect to the measure $|u|^p\, dx$ to get 
\begin{equation}\label{6.3-aux1}
  \left\||x|^{\frac{p\theta-1}{p-1}(\alpha-2)}u\right\|_p
\leq
  \left\||x|^{\theta(\alpha-2)}u\right\|_p^{\frac{p\theta-1}{p-1}}
  \left\||x|^{(\theta-\frac{1}{p})(\alpha-2)}u\right\|_p^{1-\frac{p\theta-1}{p-1}} 
\end{equation}
and  
\begin{align}\label{6.3-aux2}
|\lambda|\cos\omega_1 \left\||x|^{(\theta-\frac{1}{p})(\alpha-2)}u\right\|_p^{p\theta}
\leq 
   \left\||x|^{\theta(\alpha-2)}u\right\|_p^{p\theta-1}
  \|\lambda u-L_tu\|_p.
\end{align}
On the other hand, \eqref{reale} applied again to $ |x|^{(p\theta-1)(\alpha-2)} L$ implies that 
\begin{align} \label{g2}
%&\,
M\left\||x|^{\theta(\alpha-2)}u\right\|_p^p
%\\
 &\,\leq
{\rm Re}\,
\int_{\R^N}
  (\lambda u-Lu)|x|^{(p\theta-1)(\alpha-2)}\overline{u}|u|^{p-2}
\,dx
\\
%&\,\quad
%+k
%\int_{\R^N}
%  (|x|^{\alpha}-\min\{n,|x|^{\alpha}\})|x|^{(p\theta-1)(\alpha-2)}|u|^{p}
%\,dx
%\\
\nonumber &\,\leq
{\rm Re}\,
\int_{\R^N}
  (\lambda u -L_tu)|x|^{(p\theta-1)(\alpha-2)}\overline{u}|u|^{p-2}
\,dx
\\
\nonumber &\,\leq
\left\||x|^{\frac{p\theta-1}{p-1}(\alpha-2)}u\right\|_p^{p-1} 
\|\lambda u-L_tu\|_p.
\end{align}
where $M=f(\frac{N}{p}+\theta(\alpha-2))>0$. 
Combining \eqref{6.3-aux1} and \eqref{6.3-aux2} with the above estimate, we have
\begin{align*}
\left\||x|^{\theta(\alpha-2)u}\right\|_p^p
&\,\leq
  \frac{1}{M}\left\||x|^{\theta(\alpha-2)}u\right\|_p^{p\theta-1}
  \left\||x|^{(\theta-\frac{1}{p})(\alpha-2)}u\right\|_p^{p-p\theta}
\|\lambda u-L_tu\|_p
\\
&\,\leq
  \frac{1}{M}
  \left(\frac{1}{|\lambda|\cos\omega_1}\right)^{\frac{1-\theta}{\theta}}
  \left\||x|^{\theta(\alpha-2)}u\right\|_p^{\frac{p\theta-1}{\theta}}
  \|\lambda u-L_tu\|_p^{\frac{1}{\theta}}.
\end{align*}
Therefore we obtain
\begin{align}\label{est-theta}
\left\||x|^{\theta(\alpha-2)}u\right\|_p
\leq 
\frac{1}{M^{\theta}}
  \left(  \frac{1}{|\lambda|\cos\omega_1}\right)^{1-\theta}
  \|\lambda u-L_tu\|_p.
\end{align}
Next we prove \eqref{u-Ln}. 
From Proposition \ref{gen-pot}, we have 
\begin{align*}
{\rm Re}\left[e^{i\omega}
\int_{\R^N}
  (-L u+k|x|^{\alpha-2}u)\overline{u}|u|^{p-2}
\,dx\right]
\geq 0, \quad \omega\in[-\frac{\pi}{2}+\omega_2, \frac{\pi}{2}-\omega_2]
\end{align*}
where $\pi/2-\omega_2$ is the angle of sectoriality of $L-k|x|^{\alpha-2}$.
Since $L_t=L-V_t$ with $V_t \ge 0$, the same inequality holds for $L_t$, thus for ${\rm  Re}\, \lambda >0$ arguing as for (\ref{g1})
\begin{align*}
|\lambda|\cos\omega_2 \|u\|_p^p&\,\leq
k\left\||x|^{\frac{\alpha-2}{p}}u\right\|_p^p
+
  \|u\|^{p-1}\|\lambda u-L_tu\|_p.
\end{align*}
Since $p\theta \ge 1$ we may apply H\"older inequality to obtain the estimate
\[
  \left\||x|^{\frac{\alpha-2}{p}}u\right\|_p
\leq
  \left\|u\right\|_p^{1-\frac{1}{p\theta}}  
  \left\||x|^{\theta(\alpha-2)}u\right\|_p^{\frac{1}{p\theta}}.
\]
Then we have
\begin{align*}
\|u\|_p^{\frac{1}{\theta}}
&\,\leq
\frac{k}{|\lambda|\cos\omega_2}
  \left\||x|^{\theta(\alpha-2)}u\right\|_p^{\frac{1}{\theta}}
+
\frac{1}{|\lambda|\cos\omega_2}
  \|u\|_p^{\frac{1-\theta}{\theta}}\|\lambda u-L_tu\|_p
\\
&\,\leq
\frac{k}{|\lambda|\cos\omega_2}
  \left\||x|^{\theta(\alpha-2)}u\right\|_p^{\frac{1}{\theta}}
+
  (1-\theta)\|u\|_p^{\frac{1}{\theta}}
+
  \theta
  \left(\frac{1}{|\lambda|\cos\omega_2}\|\lambda u-L_tu\|_p\right)^{\frac{1}{\theta}}, 
\end{align*}
and hence using \eqref{est-theta}, we have
\[
\|u\|_p
\leq
\left[
\frac{k}{\theta\cos\omega_2}
\frac{1}{M}
\left(\frac{1}{\cos\omega_1}\right)^{\frac{1-\theta}{\theta}}
+
\left(\frac{1}{\cos\omega_2}\right)^{\frac{1}{\theta}}
\right]^{\theta}
\frac{1}{|\lambda|}
\|\lambda u-L_u\|_p.
\]
This a-priori estimate implies 
that $\C_+\subset \rho(L_t))$ and that  \eqref{u-Ln},  \eqref{xu-Ln} hold
for every $\lambda\in \C_+$.

\noindent
To deal with the case $\theta\in (0,\frac{1}{p})$ we consider the 
adjoint operator in $L^{p'}(\R^N)$
\[
(L_t)^{*}v=
|x|^{\alpha}\Delta + \tilde{c}|x|^{\alpha-2}x\cdot\nabla -(\tilde{b}+k)|x|^{\alpha-2}
+k \min\{t, |x|^{\alpha-2}\}, 
\]
where $\tilde{c}=2\alpha-c$ and $\tilde{b}=b+(c-\alpha)(N+\alpha-2)$, see \ref{defLtilde}) and (\ref{tilde}). 
Then, taking $\tilde{\theta}:=1-\theta\in (\frac{1}{p'},1)$, we see from (\ref{ftilde}) that
\[
\tilde{f}\left(\frac{N}{p'}+\tilde{\theta}(\alpha-2)\right)
=f\left(\frac{N}{p}+\theta(\alpha-2)\right)>0. 
\]
Thus applying \eqref{u-Ln} to $L_t^*$, we obtain
that $\C_+\subset \rho(L_t^*)$ and 
\[
\|(\lambda-L_t^*)^{-1}\|
\leq
\frac{C}{|\lambda|}, \quad \lambda\in \C_+.
\]
By duality we have \eqref{u-Ln} for $L_t$. 
Finally, let $\chi\in C^\infty(\R^N)$ satisfy 
\[
\begin{cases}
\chi\equiv 1\textrm{ in }B_{1\phantom{/2}}\textrm{ and }\chi\equiv 0\textrm{ and }B_2^c
&\textrm{ if }\alpha<2,  
\\
\chi\equiv 0\textrm{ in }B_{1/2}\textrm{ and }\chi\equiv 1\textrm{ and }B_1^c
&\textrm{ if }\alpha>2.
\end{cases}
\]
Then noting that $p\theta<1$, we obtain from (\ref{g2})
\begin{align} \label{g3}
%&\,
M
\int_{\R^N}
  |x|^{p\theta(\alpha-2)}|\chi u|^{p}
\,dx
&\,\leq
{\rm Re}\,
\int_{\R^N}
  (\lambda \chi u-L_t(\chi u))|x|^{(p\theta-1)(\alpha-2)}\chi^{p-1}\overline{u}|u|^{p-2}
\,dx
\\
\nonumber &\,\leq
\||x|^{\frac{p\theta-1}{p-1}(\alpha-2)}\chi u\|^{p-1}\|\lambda \chi u-L_t(\chi u)\|
\\
\nonumber &\,\leq
2^{\frac{(1-p\theta)|\alpha-2|}{p-1}}
\|u\|_p^{p-1}(\|\lambda u-L_tu\|_p+C_1\|u\|_p+C_2\|\nabla u\|_{L^p({\rm supp}\nabla \chi)}).
\end{align}
We note that  ${\rm supp}\nabla \chi \subset \ov{B_2}\setminus B_{1/2}$, that first and second order coefficients of $L_t$ are independent of $t$ and that the zero-order coefficients of $L_t$ are uniformly bounded with respect to $t$ in the annulus $D_4=B_4\setminus B_{1/4}$.
Therefore the interior gradient estimates 
\[
\|\nabla u\|_{L^p(B_2\setminus B_{1/2})}
\leq 
C_3(\|L_tu\|_{L^p(D_4)}+\|u\|_{L^p(D_4)}) \le C_3(\|\lambda u-L_tu\|_p+(1+|\lambda|)\|u\|_p)
\]
hold with $C$ independent of $t>0$. Using these estimates, (\ref{g2}) and  \eqref{u-Ln} we obtain for $\lambda \in C_+$, $|\lambda| \ge 1$
\[
\begin{cases}
\||x|^{\theta(\alpha-2)}u\|_{L^p(B_1)}\leq C_4\|\lambda u-L_tu\|_p.
&
{\rm if}\ \alpha<2, 
\\
\||x|^{\theta(\alpha-2)}u\|_{L^p(B_1^c)}\leq C_4\|\lambda u-L_tu\|_p.
&
{\rm if}\ \alpha>2
\end{cases}
\]
wih $C_4$ independent of $\lambda,t$.
Combining the above estimate with \eqref{u-Ln} we obtain 

\begin{equation} \label{g4}
\||x|^{\theta(\alpha-2)}u\|_p\leq C_5\|\lambda u-L_tu\|_p
\end{equation}
for 
$\lambda \in C_+$, $|\lambda| \ge 1$ and $C_5$ independent of $\lambda,t$. Finally, applying (\ref{g4}) with $\lambda=e^{i \omega}$ to $u(x)=v(sx)$ we get 
$$
\||x|^{\theta(\alpha-2)}u\|_p\leq C_5 s^{(1-\theta)(2-\alpha)}\|s^{\alpha-2}e^{i \omega}u-L_{ts^{\alpha-2}}u\|_p
$$
or, with $\eta=s^{\alpha-2} e^{i \omega}$ and $\tau=ts^{\alpha-2}$,
\[
\||x|^{\theta(\alpha-2)}(\eta -L_\tau)^{-1}\|\leq \frac{C_5}{|\eta|^{1-\theta}}.
\]
\qed

\noindent
We are now in a position to state and proof the main result 
of this section, that is the "if " part of Theorem \ref{main}.  
We recall that $I$ is the interval of all 
$\theta \in [0,1]$ such that (\ref{conjecture-other}) 
is satisfied. For every $\theta \in I$ we set 
$\alpha'=\alpha'(\theta)=\theta(\alpha-2)+2$ and define
\begin{equation} \label{Dreg}
D_{reg}(L)=\left\{
\begin{array}{l}
\Bigl \{u\in D_{max}(L)\;; 
|x|^{\alpha'}D^2u, |x|^{\alpha'-1}\nabla u, |x|^{\alpha'-2}u\in L^p(B)\ \  {\rm for\ every\ } \theta \in I\, \\ 
\qquad \qquad \qquad \quad |x|^{\alpha}D^2u, |x|^{\alpha-1}\nabla u \in L^p(B^c) \Bigr \}\quad  {\rm if\ } \alpha <2;\\ \\
\Bigr \{u\in D_{max}(L)\;; 
|x|^{\alpha'}D^2u, |x|^{\alpha'-1}\nabla u, |x|^{\alpha'-2}u\in L^p(B^c)\ {\rm for\ every\ } \theta \in I\, \\
\qquad \qquad \qquad \quad  |x|^{\alpha}D^2u, |x|^{\alpha-1}\nabla u \in L^p(B)\Bigr \}\ \   {\rm if\ }\quad  \alpha >2.
\end{array}
\right.
\end{equation}
 where $B=B_1$. Note that the maximum of regularity is achieved when $1 \in I$, that is when Theorem \ref{gen-min} applies.
\begin{teo}\label{gen-int}
If \eqref{conjecture} is satisfied, 
then $L_{int}$ generates a positive analytic semigroup in $L^p(\R^N)$ 
which is coherent with respect to all $p$ satisfying \eqref{conjecture}. Moreover,  $D_{int}(L)$ defined in (\ref{Lint}) coincides with $D_{reg}(L)$ defined above.
\end{teo}
{\sc Proof.} Fix $\lambda$ with ${\rm Re}\, \lambda>0$ and recall that, by Lemma \ref{g0}, $\lambda-L_{int}$ is injective. To show the surjectivity we fix  $f\in L^p(\R^N) $ and define $u_n=(\lambda-L_n)^{-1}f \in D_{p,\alpha}$. By Lemma \ref{re-int}
$|\lambda|\|u_n\|_p \le C\|f\|_p$ and $|\lambda|^{\theta-1}\||x|^{\theta (\alpha-2)}u_n\|_p \le C\|f\|_p$ with $C$ independent of $\lambda, n$. Note that the operators $L_n$ differ only for the zero-order coefficients which are uniformly bounded on every compact subset of $\Omega$. By local elliptic regularity, the sequence $(u_n)$ is therefore bounded in $W^{2,p}_{loc}(\Omega)$ and,  passing a subsequence, we may assume that $(u_n) \to u$ weakly in $W^{2,p}_{loc}(\Omega)$, strongly in $L^p_{loc}(\Omega)$ and pointwise. Then $\lambda u-Lu=f$ and $|\lambda|\|u\|_p \le C\|f\|_p$. Moreover  $|\lambda|^{\theta-1}\||x|^{\theta (\alpha-2)}u\|_p \le C\|f\|_p$, hence $u \in D_{int}(L)$. Note that the injectivity of $\lambda-L_{int}$ actually implies that the whole sequence $(u_n)$ converges to $u$, that is $(\lambda-L_n)^{-1}f \to (\lambda-L_{int})^{-1}f$.  If $\lambda >0$, $f \ge 0$, then $u_n \ge 0$ by Lemma \ref{g0} hence $u \ge 0$. Moreover,  if $f \in L^p(\R^N) \cap L^q(\R^N)$ the solution $u$ is independent of 
$p,q$ since so are the $u_n$, by Lemma \ref{g0}, again.

\noindent
Finally we prove the equality $D_{int}(L)=D_{reg}(L)$ and focus, as usual, on the case $\alpha <2$, the other being similar.
The inclusion $D_{reg}(L) \subset D_{int}(L)$ is obvious. Let now $u \in D_{int}(L)$ and write $u=u_1+u_2$ where $u_1=u\phi$, $u_2=u(1-\phi)$ and $\phi \in C_c^\infty (\R^N)$ with support in $B_2$ and equal to 1 in $B_1$. We introduce the operator $L_2$ on $\R^N$ in this way: the coefficients of $L_2$ coincide with those of $L$ in $B_1^c$ whereas in $B_1$ they take the (constant) value that they have on $\partial B_1$. $L_2$ is therefore uniformly elliptic with Lipschitz coefficients in $B_1$ and satisfies Hypothesis 2.1 of \cite{for-lor}. By construction the function $u_2$ belongs to the maximal domain of $L_2$ and, by \cite[Proposition 2.9]{for-lor}, $|x|^\alpha D^2 u_2, |x|^{\alpha-1}\nabla u_2 \in L^p(B^c)$, that is  $|x|^\alpha D^2 u, |x|^{\alpha-1}\nabla u \in L^p(B^c)$. To treat $u_1$ we consider the operator $L_1=|x|^{\alpha'-\alpha}L$. Since $\alpha<2$ then $\alpha' \ge \alpha$ and then $u_1 \in D_{max}(L_1)$ and, by the definition of $L_{int}$, $|x|^{\alpha'
-2}u_1 \in L^p(\R^N)$. By Lemma \ref{aux}, $u_1 \in D_{p,\alpha'}$. It follows that 
$|x|^{\alpha'}D^2u_1, |x|^{\alpha'-1}\nabla u_1, |x|^{\alpha'-2}u_1\in L^p(B)$, hence the same holds for $u$.
\qed

\noindent
We observe that $L_{int}=L_{min}$  if the conditions of Theorem \ref{gen-min} are satisfied and $L_{int}=L_{max}$ if the conditions of Theorem \ref{gen-max} hold. In both cases the equality $D_{int}(L)=D_{reg}(L)$ yields a better description of $D_{min}(L)$ and $D_{max}(L)$, respectively.

\subsection{Some consequences}

In the next proposition we show that $L_{min}$ is a generator when $N/p$ coincides with one of the endpoints of the interval $(s_1+2-\alpha, s_2+2-\alpha)$ of Theorem \ref{gen-min}.

\begin{prop}\label{endpoint-Lmin}
Let $\alpha<2$ and $N/p=s_1+2-\alpha$ or $\alpha >2$ and $N/p=s_2+2-\alpha$. Then $L_{int}$ coincides with $L_{min}$.
\end{prop}
{\sc Proof.} We only treat the case $\alpha <2$ and we first show that $\overline{Rg(I-L_{min})}=L^p(\R^N)$. 
Let $v\in L^{p'}(\R^N)$ and suppose that for every $\varphi\in C_c^\infty(\Omega)$, 
\begin{equation}\label{density}
\int_{\R^N}\overline{v}(\varphi-L\varphi)\,dx=0.
\end{equation}
Fix $\eps>0$. Since $b+\eps+(N/p+\alpha-2)(N/p'-\alpha+c)=\eps$  
it follows from Theorem \ref{gen-pot} that 
the minimal realization of $L-\eps|x|^{\alpha-2}$ generates analytic semigroup in $L^p(\R^N)$ 
and its domain is $D_{p,\alpha}$.
Hence \eqref{density} holds for in $D_{p,\alpha}$ and we deduce from the invertibility of $I-L+\eps |x|^{\alpha-2}$ that for every $f\in L^p(\R^N)$, 
\begin{equation*}
\int_{\R^N}\overline{v}f\,dx=
\int_{\R^N}\overline{v}\left(\eps |x|^{\alpha-2}(I-L+\eps|x|^{\alpha-2})^{-1}f\right)\,dx.
\end{equation*}
Choosing $f=|v|^{p'-2}v$ and 
setting $w_\eps=(I-L+\eps|x|^{\alpha-2})^{-1}f$,  
%and letting $m\to \infty$, 
we have
\begin{align}\label{density2}
\int_{\R^N}|v|^{p'}\,dx&=
\int_{\R^N}\overline{v}\left(\eps|x|^{\alpha-2}w_\eps\right)\,dx.
\end{align}
Observe that Lemma \ref{dual}  implies that 
\begin{align*}
\left\||x|^{\frac{\alpha-2}{p'}}w_{\eps}\right\|_p^p
+
\eps\left\||x|^{\alpha-2}w_{\eps}\right\|_p^p
&\leq 
\int_{\R^N} |x|^{(p-1)(\alpha-2)}\overline{w_\eps}|w_\eps|^{p-2}(I-L+\eps|x|^{\alpha})w_\eps\,dx
\\
&\leq
\|f\|_p \left\||x|^{\alpha-2}w_\eps\right\|_p^{p-1}.
\end{align*}
This yields 
\begin{align*}
\left\|\eps|x|^{\alpha-2}w_\eps\right\|_p\leq \|f\|_p
\text{\quad and \quad}
\left\|\left(\eps|x|^{\alpha-2}\right)^{\frac{1}{p'}}w_{\eps}\right\|_p\leq \|f\|_p.
\end{align*}
Using the above estimates, we obtain 
\begin{equation}\label{conc}
\eps|x|^{\alpha-2}w_\eps\to 0\quad\text{weakly in }L^p(\R^N).
\end{equation}
In fact, for every $\phi\in C_c^\infty(\Omega)$, we have
\begin{align*}
\left|\int_{\R^N}\left(\eps|x|^{\alpha-2}w_\eps\right)\phi\,dx\right|
&=
\eps^{\frac{1}{p}}
\left|\int_{\R^N}\left(\eps|x|^{\alpha-2}\right)^{\frac{1}{p'}}w_\eps|x|^{\frac{\alpha-2}{p}}\phi\,dx\right|
\leq \eps^{\frac{1}{p}}\|f\|_p\left\||x|^{\frac{\alpha-2}{p}}\phi\right\|_{p'}
\to 0
\end{align*}
as $\eps \to 0$. Since $\{\eps|x|^{\alpha-2}w_\eps\}$ is bounded in $L^p(\R^N)$, 
a density argument implies \eqref{conc}. Consequently, 
combining \eqref{conc} with \eqref{density2}, we obtain $v=0$.
This means that $\overline{Rg(I-L_{min})}=L^p(\R^N)$.

\noindent Finally, we prove $L_{int}=L_{min}$. The inclusion $L_{int}\supset L_{min}$ is obvious. 
Conversely, let $u\in D_{int}(L)$. 
Since $\overline{Rg(I-L_{min})}=L^p(\R^N)$, 
we can choose $u_n\in D_{min}(L)\subset D_{int}(L) $ such that 
$(I-L_{min})u_n\to (I-L_{int})u$. Since $I-L_{int}$ is invertible 
we have 
\[
\|u_n-u\|_p\leq C\|(I-L_{int})(u_n-u)\|_p\to 0 
\]
 and the closedness of $L_{min}$ implies $u\in D(L_{min})$. 
\qed

\begin{os}
The equality $\overline{Rg(I-L_{min})}=L^p(\R^N)$ is true even when $\alpha<2$ and  $N/p=s_2+2-\alpha$ or $\alpha >2$ and $N/p=s_1+2-\alpha$, by the same proof as above. 
However, in these cases, the injectivity of $I-L_{min}$ breaks down,  see Proposition \ref{nonex}.
\end{os}
By duality one obtains a similar result for $L_{max}$, see Remark \ref{endpoints}.
\begin{prop}\label{endpoint-Lmin1}
Let $\alpha<2$ and $N/p=s_2$ or $\alpha >2$ and $N/p=s_1+2-\alpha$. Then $L_{int}$ coincides with $L_{max}$.
\end{prop}
We end this section with some remarks on $L_{int}$. We fix $\theta \in I$, that is satisfying (\ref{conjecture-other}), and define $L_\theta$ through the domain
\begin{equation*} \label{Ltheta}
D_{\theta}(L)=\{u\in D_{max}(L)\subset L^p(\R^N)\;;\; |x|^{\theta(\alpha-2)}u\in L^p(\R^N)\}.
\end{equation*}
Clearly $L_{int} \subset L_\theta$. However, since $I-L_{int}$ is invertible and $I-L_\theta$ is injective, by Lemma \ref{inj-int} (whose proof works for any fixed $\theta$), then both operator coincide and $L_{int}=L_\theta$. This means that  the extra integrability condition $|x|^{\theta (\alpha-2)}u \in L^p(\R^N)$, $u \in D_{max}(L)$, for a fixed $\theta \in I$ extends automatically to  every $\theta \in I$. 

\noindent
In the next proposition we show that, unless $1 \in I$, this integrability condition does not hold for $\theta_0=\sup I$. Note that $1 \in I$ is equivalent to say that Theorem \ref{gen-min} applies and is more restrictive than requiring that $L_{min}$ generates. Note also that $\theta_0$ can be equal to 1 even though $1\not \in I$.

\begin{prop} \label{lasttheta} Assume that (\ref{conjecture}) holds  and that $1\notin I$,  
Set $\theta_0=\sup I$ and $\alpha_0'=\theta_0(\alpha-2)+2$. 
Then there exists $u\in D_{int}(L)$ such that 
$|x|^{\alpha_0'-2}u\notin L^p(\R^N)$.
\end{prop}
{\sc Proof.} We give a proof only $\alpha<2$. In this case (\ref{conjecture}) reads $s_1 <N/p <s_2+2-\alpha$. Since $1 \not \in I$, then $f(N/p+\theta_0(\alpha-2))=0$ and then 
 $s_1=N/p+\theta_0(\alpha-2)$.
We set $u(x)=|x|^{-s_1}\zeta (x)$,
where $\zeta \in C_c^\infty(\R^N)$ is one in the unit ball $B_1$  and zero outside the ball $B_2$.
Then for every $\theta\in[0,\theta_0)$, we have 
\[
|x|^{\theta(\alpha-2)}u= |x|^{\theta(\alpha-2)-s_1}\zeta(x)=
|x|^{-\frac{N}{p}+(\theta_0-\theta)(2-\alpha)}\zeta(x)\in L^p(\R^N).
\]
 Since $L|x|^{-s_1}=0$ then $Lu\in C_c^\infty(\Omega)$ and therefore  $u\in D_{int}(L)$. 
However  $|x|^{\alpha_0'-2}u\notin L^p(\R^N)$. 
\qed

\noindent Finally let us show that for $\lambda >0$, $f \ge 0$, $(\lambda-L_{int})^{-1}f$  is the minimal  among the positive solutions $u\in D_{max}(L)$ of the equation $\lambda u-Lu=f$. This characterizes the generated semigroup as the minimal one and is important when $L_{int}$ differs both from $L_{min}$ and $L_{max}$.
First prove a maximum principle for the operator $L$ restricted to the annulus $\Omega_\eps$. Note that the classical maximum principle does not hold when $b <0$.
\begin{lem} \label{max-princ}
Let $\lambda>0$, $g\leq 0$ and let $u\in W^{2,p}(\Omega_\eps)$ solve $\lambda u-Lu=g$ in $\Omega_\eps$ with $u\leq 0$ at the boundary. Then $u\leq 0$ in $\Omega_\eps$.
\end{lem}
{\sc Proof.}
Let $\theta$ be such that $f\left(\frac{N}{p}+\theta (\alpha-2)\right)>0$. We multiply the equation $\lambda u-Lu=g$ by $|x|^{(p\theta-1)(\alpha-2)}(u^{+})^{p-1}$ and integrate over $\Omega_\eps$.  We proceed as in  Proposition \ref{gen-pot} whith $\alpha$ replaced by $\beta=\alpha+(p\theta-1)(\alpha-2)$ and observe that, since $u\leq 0$ on the boundary, no boundary terms appear after integration by parts. Setting $v=|x|^{\frac{N-2+\beta}{p}}u$ we obtain the analogous of (\ref{zeroreal}) 
\begin{align*}
\lambda \int_{\Omega_{\eps}}|x|^{(p\theta-1)(\alpha-2)}(u^+)^{p}&+(p-1)\int_{\Omega_\eps}|x|^{2-N}|\nabla v^+|^2  (v^+)^{p-2}+f\left(\frac{N}{p}+\theta (\alpha-2)\right)\int_{\Omega_{\eps}}|x|^{-N}(v^+)^p\\&=\int_{\Omega_{\eps}}|x|^{(p\theta-1)(\alpha-2)}g(u^+)^{p-1}\leq 0.
\end{align*}
It follows that 
$\int_{\Omega_{\eps}}|x|^{(p\theta-1)(\alpha-2)}(u^+)^{p}\leq 0$ and therefore $u^+=0$ in $\Omega_\eps$.
\qed

\begin{prop} \label{minimal-sem}
Let $\lambda >0$, $f \ge 0$ and let $0 \le u \in D_{max}(L)$ satisfy $\lambda u -Lu=f$. Then $ (\lambda-L_{int})^{-1}f\le u$.
\end{prop}
{\sc Proof} 
Let $u_\eps \in W^{2,p}(\Omega_\eps) \cap W^{1,p}_0(\Omega_\eps)$ be such that $\lambda u_\eps-Lu_\eps=f$ in $\Omega_\eps$. Then $u_\eps \ge 0$,   $\lambda (u_\eps-u)-L(u_\eps-u)=0$ and $u_\eps-u \le 0$ at the boundary.  By Lemma \ref{max-princ}, $u_\eps\leq u$ in $\Omega_\eps$. If $v=(\lambda-L_{int})^{-1}f$ the same argument shows that $u_\eps \le v$ in $\Omega_\eps$. Moreover, if $\eps_1 <\eps_2$, then $u_{\eps_2} \le u_{\eps_1}$ in $\Omega_{\eps_2}$, by the above lemma again.
Then  $(u_\eps)$ converges pointwise as $\eps \to 0$ to some $u_0 \le v$. By dominated convergence $u_\eps \to u_0$ in $L^p(\R^N)$ and, by elliptic interior estimates, also in $W^{2,p}_{loc}(\R^N)$. Then $\lambda u_0-Lu_0=f$ and $u_0 \in D_{max}(L)$. Since $0\le u_0 \le v$, then $u_0 \in D_{int}(L)$, hence $u_0=v$, by uniqueness.
Since $u_\eps\leq u$, letting $\eps \to 0$, it follows that $v \le u$.
\qed

%%%%% Section 6 %%%%%
\section{The critical case: $b+\left(\frac{N-2+c}{2}\right)^2=0$}
We always assume $\alpha \neq 2$ and show that a suitable realization of 
$L_{min} \subset L \subset  L_{max}$ generates a semigroup in $L^p(\R^N)$ if and only if 
\begin{equation}\label{conjecture1}
s_1+\min\{0,2-\alpha\} \le \frac{N}{p} \le s_2+\max\{0,2-\alpha\}.
\end{equation}
Note that the endpoints above are included, whereas they are excluded in (\ref{conjecture}).
In this case the function $f$ defined in (\ref{deff}) is negative except for $s=s_0=(N-2+c)/2$ where it vanishes and both $s_1$ and $s_2$ coincide with $s_0$.
We first consider the case $\alpha<2$ and we give full proofs following the method of Section 5, but adding logarithmic weights in the resolvent estimates. Moreover, we consider the operator first in the unit ball $B_1$ and then we use a gluing procedure to treat the case of the whole space. The case $\alpha>2$ will be shortly  considered in Subsection \ref{maggiore2}.

\subsection{Positive results for $\alpha<2$}
We always assume $\alpha <2$ in this subsection and fix $\theta_0\in [0,1]$ such that $\frac{N}{p}=s_0+\theta_0(2-\alpha)$,  $s_0=\frac{N-2+c}{2}$.
\begin{teo}\label{8.1}
 Assume that  
\begin{equation}\label{conj-critical}
s_0\leq\frac{N}{p}\leq s_0+2-\alpha
\end{equation}
 and define $L_{int}$ through the domain
\begin{equation} \label{Lint-critical}
D_{int}(L)=\{u\in D_{max}(L)\;;\;|x|^{\theta_0(\alpha-2)}\bigl|\log|x|\bigr|^{-\frac{2}{p}}u\in L^p(B_{1/2})\}.
\end{equation}
Then
$L_{int}$
generates a positive analytic semigroup in $L^p(\R^N)$  which is coherent with respect to all $p$ satisfying \eqref{conj-critical}. 
\end{teo}
For technical reasons we need also the the operator $L$ in the unit ball $B_1$ (with Dirichlet boundary conditions) defined on the domain
\begin{align} \label{Lint-critical1}
\nonumber D^1_{int}(L)=&\{u\in L^p(B_1) \cap W^{2,p}(B_1 \setminus B_\eps)\ \forall \ 0<\eps<1, Lu \in L^p(B_1); \\&|x|^{\theta_0(\alpha-2)}\bigl|\log|x|\bigr|^{-\frac{2}{p}}u\in L^p(B_{1/2}),\  u=0\  {\rm on\ } \partial B_1\}, 
\end{align}
where $\theta $ is as before. We denote this operator by $L^1_{int}$. Most computations will be performed on the set
\[
D_1:=\{u\in C_c^\infty(\overline{B}_1\setminus\{0\}):\;u=0\text{ on }\partial B_1\}.
\]
In the next proposition we prove the injectivity of $\lambda- L_{int} $ and $\lambda-L^1_{int}$ for positive $\lambda$.

\begin{prop}\label{inj} The operators  $\lambda-L_{int}$ and $\lambda-L^1_{int}$  are injective for $\lambda >0$.
\end{prop}
{\sc Proof.} We start with $L_{int}$. We denote by $\Delta_{S^{N-1}}$ the Laplace Beltrami on the unit sphere $S^{N-1}$. If $Q$ is a spherical harmonic of order $n \ge 0$, then $-\Delta_{S^{N-1}}Q=\lambda_n Q$ with $\lambda_n=n(n+N-2)$.
If $v\in N(\lambda-L_{int})$
we set 
\[
v_Q(r)=\int_{S^{N-1}}v(r,\omega)Q(\omega)\,d\omega, 
\]
where $Q$ is a spherical harmonic of order $n$. Then
\begin{equation}\label{v_q-profile1}
%r^{N-1+\theta_0(\alpha-2)p}
r^{ps_0-1}\bigl|\log|x| \bigr|^{-2}|v_Q(r)|^p\in L^1(0,1/2).
\end{equation}
Observe that
\begin{align*}
&v_Q''+\frac{N-1}{r}v_Q'=\int_{S^{N-1}}\left (v_{rr}+\frac{N-1}{r}v_r\right )Q(\omega) \, d\omega=\int_{S^{N-1}}\left (\Delta v -\frac{\Delta_{S^{N-1}} v}{r^2}\right )Q \, d\omega \\
&\int_{S^{N-1}}\left (Q\Delta v -v\frac{\Delta_{S^{N-1}} Q}{r^2}\right ) \, d\omega=\int_{S^{N-1}}\left (\Delta v +\frac{\lambda_n v} {r^2}\right )Q \, d\omega.
\end{align*}
This implies that $v_Q$ satisfies 
\[
\lambda v_Q-r^\alpha
          \left(
            v_Q''+\left(N-1+c\right)\frac{v_Q'}{r}-
            \left(
            b+\lambda_n
            \right)\frac{v_Q}{r^2}
          \right)=0.
\]
We use Lemma \ref{bhv} to show that $v_Q=0$. The integrabilty of $v_Q$ at $r=\infty$ and\eqref{u1-behave} imply that $v_Q=cu_2$. If $n>0$, by 
\eqref{u2-behave} 
and \eqref{v_q-profile1} we see that $c=0$, that is $v_Q=0$. 
If $n=0$,  $u_2$ behaves like $r^{-s_0}\log r$ near $0$, see \eqref{u1-behave-res}, and hence 
\eqref{v_q-profile1} is not satisfied unless $v_Q=0$. The density of spherical harmonics in $L^{p'}(S^{N-1})$ yields $v(r, \cdot)=0$ for every $r$, hence $v=0$ and this concludes the proof for $L_{int}$. In the case of $L^1_{int}$ the proof is similar: if $v_Q \not \equiv 0$ then $v_Q=c_1u_1+c_2u_2$ with $c_1 \not =0$, $c_2 \not =0$ since $v_Q(1)=0$ and $u_1, u_2$ are positive. Hence $v_Q$ behaves like $u_2$ near 0 and  \eqref{v_q-profile1} is not satisfied.
\qed

\noindent The following weighted estimates will be crucial in what follows.
\begin{prop}\label{log-est}
For every 
$v\in D_1$
\begin{align}\label{log}
Re\int_{B_1}
  |x|^{2-N}\nabla v \cdot\nabla (\overline{v}|v|^{p-2})\,dx
\geq 
\frac{p-1}{p^2}\int_{B_1}
|x|^{-N}\bigl|\log |x|\bigr|^{-2}|v|^{p}
\,dx.
\end{align}
In particular, 
if $u\in D_1$, 
and $v=|x|^{\frac{N-2+c}{2}}u$, then 
\begin{align} \label{log1}
Re\int_{B_1}
  (-Lu)|x|^{(p\theta_0-1)(\alpha-2)}\overline{u}|u|^{p-2}
\,dx
&\geq 
\frac{p-1}{p^2}\int_{B_1}
|x|^{p\theta_0(\alpha-2)}\bigl|\log |x|\bigr|^{-2}|u|^{p}
\,dx.
\end{align}
\end{prop}
Noting that $r^{\eps}|\log r|\leq (\eps e)^{-1}$ if $r\in(0,1)$ and $\eps>0$, we obtain from (\ref{log1})
%%%%%%%%%%%%%%%%%%%%%%%%%%%%%%%%%%%%%%%%%%%%%%%%%%%%%%%%%%%%%%%%%%%%%
\begin{lem}\label{power-est}
For every $\delta>0, \theta_0>0$ and 
$u\in D_1$, 
%C_c^\infty(\overline{B}\setminus\{0\})$  satisfying $u=0$ on $\partial B$, 
\begin{align}
Re\int_{B_1}
  (-Lu)|x|^{(p\theta_0-1)(\alpha-2)}\overline{u}|u|^{p-2}
\,dx
\geq 
\frac{(p-1)(2-\alpha)^2e^2\delta^2}{4}\int_{B_1}
  |x|^{p(\theta_0-\delta)(\alpha-2)}|u|^{p}
\,dx.
\end{align}
\end{lem}
{\sc Proof.} 
If $v\in D_1$ integrating by parts we obtain
\begin{align*}
\int_{B_1}
  |x|^{2-N} \nabla v \cdot \nabla (\overline{v}|v|^{p-2})
\,dx=-\int_{B_1}
  |x|^{-N}(|x|^2\Delta v+(2-N)x\cdot\nabla v)\overline{v}|v|^{p-2}
\,dx.
\end{align*}
Observe that in spherical coordinates 
\[
|x|^2\Delta v+(2-N)x\cdot\nabla v
=r^2\frac{\partial^2 v}{\partial r^2}+r\frac{\partial v}{\partial r}+\Delta_{S^{N-1}}v,
\]
and
\begin{align}
\nonumber
&
\int_{B_1}
  |x|^{2-N}\nabla v \cdot \nabla (\overline{v}|v|^{p-2})\rangle
\,dx
%\int_{B_1}
%  |x|^{-N}(|x|^2\Delta v+(2-N)x\cdot\nabla v)\overline{v}|v|^{p-2}
%\,dx
=
-\int_{B_1}
r^{-N}\left(r^2\frac{\partial^2 v}{\partial r^2}+r\frac{\partial v}{\partial r}+\Delta_{S^{N-1}}v\right)
\overline{v}|v|^{p-2}
\,dx
\\
\label{lem3.5-eq1}
&=
-\int_{S^{N-1}}
\int_{0}^1
\left(r^2\frac{\partial^2 v}{\partial r^2}+r\frac{\partial v}{\partial r}\right)
\overline{v}|v|^{p-2}
\,\frac{dr}{r}\,d\omega
+
\int_{0}^1
\left(\int_{S^{N-1}}
(-\Delta_{S^{N-1}}v)
\overline{v}|v|^{p-2}
\,d\omega\right)\,\frac{dr}{r}.
\end{align}
Since $\Delta_{S^{N-1}}$ is dissipative in $L^p(S^{N-1})$ we see that for every $r>0$ 
\begin{align}\label{lem3.5-eq2}
Re\int_{S^{N-1}}
(-\Delta_{S^{N-1}}v)
\overline{v}|v|^{p-2}
\,d\omega\geq 0.
\end{align}
On the other hand, fix $\omega\in S^{N-1}$ and set $w(s)=v(e^{s},\omega)$. Then 
\begin{align*}
\int_{0}^1
\left(r^2\frac{\partial^2 v}{\partial r^2}(r,\omega)+r\frac{\partial v}{\partial r}(r,\omega)\right)
\overline{v(r,\omega)}|v(r,\omega)|^{p-2}
\,\frac{dr}{r}
=
\int_{-\infty}^0
\left(-w''\right)
\overline{w}|w|^{p-2}\,ds.
\end{align*}
Using Hardy's inequality we have 
\begin{align*}
Re\int_{-\infty}^0
\left(-w''\right)
\overline{w}|w|^{p-2}\,ds
&\geq 
\frac{p-1}{p^2}\int_{-\infty}^0s^{-2}|w|^{p}\,ds
=
\frac{p-1}{p^2}\int_{0}^1|\log r|^{-2}|v(r,\omega)|^{p}\,\frac{dr}{r}. 
\end{align*}
Therefore we deduce that 
\begin{align}
\nonumber
&-Re\int_{S^{N-1}}
\int_{0}^1
\left(r^2\frac{\partial^2 v}{\partial r^2}+r\frac{\partial v}{\partial r}\right)
\overline{v}|v|^{p-2}
\,\frac{dr}{r}\,d\omega
\geq
\frac{p-1}{p^2}\int_{S^{N-1}}
\int_{0}^1
|\log r|^{-2}|v(r\omega)|^{p}\,\frac{dr}{r}
\,d\omega
\\
\label{lem3.5-eq3}
&=
\frac{p-1}{p^2}\int_{B_1}
|x|^{-N}\bigl|\log |x|\bigr|^{-2}|v|^{p}
\,dx.
\end{align}
Combining \eqref{lem3.5-eq2} and \eqref{lem3.5-eq3} with \eqref{lem3.5-eq1}, 
we obtain \eqref{log}.
To prove (\ref{log1}), we consider the operator $x^\gamma L$ with $\gamma=(p\theta_0-1)(\alpha-2)$. Then $(N-2+\alpha+\gamma)=(N-2+c)/2=s_0$ and Proposition \ref{gen-pot} applies. In particular, since $f(s_0)=0$, (\ref{zeroreal}) with $v=|x|^{s_0}u$ yields
$$
Re\int_{B_1}
  (-Lu)|x|^{(p\theta_0-1)(\alpha-2)}\overline{u}|u|^{p-2}
=Re\int_{B_1}
  |x|^{2-N}\nabla v \cdot\nabla (\overline{v}|v|^{p-2})\,dx
$$ and (\ref{log1}) follows from (\ref{log}).
\qed

\noindent To prove Theorem \ref{8.1}, as in Section 5, 
we introduce the operator 
\begin{equation*}
L_tu:=Lu - k|x|^{\alpha-2}+k\min\{|x|^{\alpha-2},t \},
\end{equation*}
defined in $D_{p,\alpha}$, see (\ref{domain}),  and, in order to consider separately the singularities at infinity and near the origin, we introduce also the operators 
$L^1_t$ and $L^2_t$ in $B_1$ and $B_{1/2}^c$ defined as $L_t$ on the domains
\begin{equation}\label{domain1}
D^1_{p,\alpha}=\{u\in L^p(B_1) \cap  W^{2,p}(B_{1}\setminus B_\eps)\  \forall\  0<\eps<1,  \ |x|^\alpha D^2u,\ |x|^{\alpha-1} \nabla u,\ |x|^{\alpha-2}u\in L^p(B_1)\}
\end{equation}
and
\begin{equation}\label{domain2}
D^2_{p,\alpha}=\{u\in W^{2,p}(B_{1/2}^c, \ |x|^\alpha D^2u,\ |x|^{\alpha-1} \nabla u,\in L^p(B_{1/2}^c)\},
\end{equation}
respectively. Here $k$ is a fixed constant  large enough, $t>0$ is a parameter and both the operators are endowed with Dirichlet boundary conditions. 
As in Theorem \ref{gen-min} or Proposition \ref{gen-pot} we have, since $k$ is large, the following lemma.
\begin{lem}\label{gen.Lt}
For every $t>0$, 
$L^1_t$
generates an analytic semigroups on $L^p(\R^N)$. 
Moreover, 
%the set 
%$\{u\in C_c^\infty(\overline{B}_1\setminus\{0\})\;;\;u=0\text{ on }\partial B_1\}$
$D_1$ is a core for $L_t^1$.
\end{lem}
Now we state a-priori estimates for $L_t^1$. Due to the presence of the logarithmic term, we cannot directly prove  that (\ref{est-critical}) holds with $\gamma=1$, as in  Lemma \ref{re-int}.
\begin{prop} \label{res-est}
For every $\gamma\in (0,1)$ there are constants $C_\gamma, C'>0$ such that for every $t>0$, $\lambda \in \C_+$ and $u \in D^1_{p,\alpha}$ 
%$u\in C_c^\infty(\ov{B_1}\setminus\{0\})$, 
\begin{equation}\label{est-critical}
\|u\|_{L^p(B_1)}\leq \frac{C_\gamma}{|\lambda|^{\gamma}}
\|\lambda u-L^1_tu\|_{L^p(B_1)}
\end{equation}
and 
\begin{equation}\label{w-est-critical}
\left\|
  |x|^{\theta_0(\alpha-2)}\bigl|\log|x|\bigr|^{-\frac{2}{p}}u
\right\|_{L^p(B_1)}\leq C'\|\lambda u-L^1_tu\|_{L^p(B_1)}.
\end{equation}
If $\theta_0=\frac{1}{p}$, then $\gamma=1$ is allowed.
\end{prop}
{\sc Proof.} We proceed as in Lemma \ref{re-int}, in $B_1$ rather than $\R^N$, and consider first the case when $\theta_0 \ge \frac{1}{p}$. By density we may assume that $u \in D_1$. The estimate 
\begin{equation} \label{weighted1}
|\lambda|\left\||x|^{(\theta_0-\frac{1}{p})(\alpha-2)}u\right\|_{L^p(B_1)}^p
\leq 
C\|\lambda u-L^1_tu\|_{L^p(B_1)}
\left\||x|^{\frac{p\theta_0-1}{p-1}(\alpha-2)}u\right\|_{L^p(B_1)}^{p-1}.
\end{equation}
is identical to (\ref{weighted}) and obtained as in  Lemma \ref{re-int}, recalling that Proposition \ref{gen-pot} holds in the critical case. Observe that if $\theta_0=\frac{1}{p}$ (or 
$\frac{N-2+\alpha}{p}=\frac{N-2+c}{2}$) the above inequality gives (\ref{est-critical}) with $\gamma=1$. 
If $\theta_0>\frac{1}{p}$ we note that
\[
0 < \theta_0-\frac{1}{p}<\frac{p\theta_0-1}{p-1}
=
\theta_0-\frac{1-\theta_0}{p-1}.
\]
We choose $\delta\in ]0,\frac{1-\theta_0}{p-1}]$ and 
we see from Lemma \ref{intgamma} 
\[
\left\||x|^{\frac{p\theta_0-1}{p-1}(\alpha-2)}u\right\|_{L^p(B_1)}
\leq
\left\||x|^{(\theta_0-\delta)(\alpha-2)}u\right\|_{L^p(B_1)}^\tau
\left\||x|^{(\theta_0-\frac{1}{p})(\alpha-2)}u\right\|_{L^p(B_1)}^{1-\tau},
\]
where  
\[
\tau= \frac{p\theta_0-1} {(p-1)(1-p\delta)}. 
\]
Using the above inequality in (\ref{weighted1}) we have
\begin{equation}\label{aux-cri}
|\lambda|
\left\||x|^{(\theta_0-\frac{1}{p})(\alpha-2)}u\right\|_{L^p(B_1)}^{1+(p-1)\tau}
\leq 
C\|\lambda u-L^1_tu\|_{L^p(B_1)}
\left\||x|^{(\theta_0-\delta)(\alpha-2)}u\right\|_{L^p(B_1)}^{(p-1)\tau}.
\end{equation}
On the other hand, by Lemma \ref{power-est} we see that 
\begin{align*}
&
\frac{(p-1)(\alpha-2)^2e^2\delta^2}{4}
\left\||x|^{(\theta_0-\delta)(\alpha-2)}u\right\|_{L^p(B_1)}^{p} \\
&\leq
\|\lambda u-L^1_tu\|_{L^p(B_1)}
\left\||x|^{\frac{p\theta_0-1}{p-1}(\alpha-2)}u\right\|_{L^p(B_1)}^{p-1}
\\
&\leq
\|\lambda u-L^1_tu\|_{L^p(B_1)}
\left\||x|^{(\theta_0-\delta)(\alpha-2)}u\right\|_{L^p(B_1)}^{(p-1)\tau}
\left\||x|^{(\theta_0-\frac{1}{p})(\alpha-2)}u\right\|_{L^p(B_1)}^{(p-1)(1-\tau)}.
\end{align*}
Combining the previous estimate with \eqref{aux-cri}, we have
\begin{align*}
&
\frac{(p-1)(\alpha-2)^2e^2\delta^2}{4}
\left\||x|^{(\theta_0-\delta)(\alpha-2)}u\right\|_{L^p(B_1)}^{p-(p-1)\tau}
\\
&\leq
\|\lambda u-L^1_tu\|_{L^p(B_1)}
\left\||x|^{(\theta_0-\frac{1}{p})(\alpha-2)}u\right\|_{L^p(B_1)}^{(p-1)(1-\tau)}
\\
&\leq
\|\lambda u-L^1_tu\|_{L^p(B_1)}
\left(
\frac{C}{|\lambda|}
\|\lambda u-L^1_tu\|_{L^p(B_1)}
\left\||x|^{(\theta_0-\delta)(\alpha-2)}u\right\|_{L^p(B_1)}^{(p-1)\tau}
\right)^{\frac{(p-1)(1-\tau)}{1+(p-1)\tau}}.
\\
&=
\|\lambda u-L^1_tu\|_{L^p(B_1)}^{\frac{p}{1+(p-1)\tau}}
\left(
\frac{C}{|\lambda|}
\left\||x|^{(\theta_0-\delta)(\alpha-2)}u\right\|_{L^p(B_1)}^{(p-1)\tau}
\right)^{\frac{(p-1)(1-\tau)}{1+(p-1)\tau}}.
\end{align*}
This yields
\begin{align}\label{w-est-pre}
\left\||x|^{(\theta_0-\delta)(\alpha-2)}u\right\|_{L^p(B_1)}
\leq
\left(
\frac{4}{(p-1)(\alpha-2)^2e^2\delta^2}
\right)^{\frac{1+(p-1)\tau}{p}}
\left(
\frac{C}{|\lambda|}
\right)^{\frac{(p-1)(1-\tau)}{p}}
\|\lambda u-L^1_tu\|_{L^p(B_1)}.
\end{align}
Next, we prove \eqref{est-critical}. By Proposition \ref{gen-pot}, 
we have 
\[
Re\left[e^{i\omega}\int_{B_1}(-L^1_tu+k|x|^{\alpha-2}u)\overline {u}|u|^{p-2}\,dx\right]\geq 0,\quad \omega\in [-\frac{\pi}{2}+\omega_0,\frac{\pi}{2}-\omega_0],
\]
for some $\omega_0\in (0,\frac{\pi}{2})$. Thus we see that 
\begin{align*}
\|u\|_{L^p(B_1)}^p
&\leq 
\frac{1}{|\lambda|\cos \omega_0}\|\lambda u - L^1_tu\|_{L^p(B_1)}\|u\|_{L^p(B_1)}^{p-1}
+\frac{k}{|\lambda|\cos \omega_0}\left\||x|^{\frac{\alpha-2}{p}}u\right\|_{L^p(B_1)}^p.
\end{align*}
Fix $\delta<\theta_0-\frac{1}{p}$. Then $p(\theta_0-\delta)>1$ and
H\"older's inequality yields 
\[
\left\||x|^{\frac{\alpha-2}{p}}u\right\|_{L^p(B_1)}
\leq 
\left\|u\right\|_{L^p(B_1)}^{1-\frac{1}{p(\theta_0-\delta)}}
\left\||x|^{(\theta_0-\delta)(\alpha-2)}u\right\|_{L^p(B_1)}^{\frac{1}{p(\theta_0-\delta)}}.
\]
Hence using Young's inequality and \eqref{w-est-pre} we obtain 
\begin{align*}
\|u\|_{L^p(B_1)}^{\frac{1}{\theta_0-\delta}}
&\leq 
\frac{1}{|\lambda|\cos \omega_0}\|\lambda u - L^1_tu\|_{L^p(B_1)}\|u\|_{L^p(B_1)}^{\frac{1}{\theta_0-\delta}-1}
+
\frac{k}{|\lambda|\cos \omega_0}\left\||x|^{(\theta_0-\delta)(\alpha-2)}u\right\|_{L^p(B_1)}^{\frac{1}{\theta_0-\delta}}
\\
&\leq 
(\theta_0-\delta)\left(
\frac{1}{|\lambda|\cos \omega_0}
\|\lambda u - L^1_tu\|_{L^p(B_1)}
\right)^{\frac{1}{\theta_0-\delta}}
+
(1-\theta_0+\delta)
\|u\|_{L^p(B_1)}^{\frac{1}{\theta_0-\delta}}
\\
&\quad+
\frac{kC_\delta^{\frac{1}{\theta_0-\delta}}}{\cos \omega_0}
|\lambda|^{-1-\frac{(p-1)(1-\tau)}{p(\theta_0-\delta)}}\|\lambda u - L^1_tu\|_{L^p(B_1)}^{\frac{1}{\theta_0-\delta}}.
\end{align*}
Consequently, 
setting $\gamma(\delta):=\theta_0-\delta+\frac{(p-1)(1-\tau)}{p}$, we have  \eqref{est-critical}
\[
\|u\|_{L^p(B_1)}
\leq \left(\frac{C_\delta'}{|\lambda|}+\frac{C_{\delta}''}{|\lambda|^{\gamma(\delta)}}\right)
\|\lambda u- L^1_tu\|_{L^p(B_1)}
\]
and we note that $\lim_{\delta\downarrow 0}\gamma(\delta)=1$. Next, we prove \eqref{w-est-critical}. By Lemma \ref{log-est}, since $L_t^1=L-V_t$ with $V_t \ge 0$, we have 
\begin{align*}
\frac{p-1}{p^2}\left\||x|^{\theta_0(\alpha-2)}\bigl|\log|x|\bigr|^{-\frac{2}{p}}u\right\|_{L^p(B_1)}^p
&\leq
Re\int_{B_1}(\lambda u-L^1_tu)|x|^{(p\theta_0-1)(\alpha-2)}\ov{u}|u|^{p-2}\,dx
\\
&
\leq\|\lambda u-L^1_tu\|_{L^p(B_1)}
\left\||x|^{\frac{p\theta_0-1}{p-1}(\alpha-2)}u\right\|_{L^p(B_1)}^{p-1}.
\end{align*}
From estimate \eqref{w-est-pre} with $\delta=\frac{1-\theta_0}{p-1}$, hence $\tau=1$, we obtain 
\begin{align*}
\left\||x|^{\theta_0(\alpha-2)}\bigl|\log|x|\bigr|^{-\frac{2}{p}}u\right\|_{L^p(B_1)}^p
\leq C\|\lambda u-L^1_tu\|_{L^p(B_1)}^p.
\end{align*}
Finally, we consider the case $0 \le \theta_0<\frac{1}{p}$ by using the dual operator $(L_t^1)^*$ in $L^{p'}(B_1)$ and  proceeding as in Lemma \ref{re-int}. 
One verifies that $(L_t^1)^*$ 
satisfies $\frac{N-2+\tilde{c}}{2}=\frac{N}{p'}+\tilde{\theta}(\alpha-2)$ 
with $\tilde{\theta}1-\theta > \frac{1}{p'}$. Therefore applying \eqref{est-critical} to $(L_t^1)^*$ in $L^{p'}(B_1)$, 
we have $\|(\lambda - (L_t^1)^{*})^{-1}\|\leq C(|\lambda|^{-1}+|\lambda|^{-\gamma})$ and 
by duality we obtain \eqref{est-critical} for $L_t^1$ in $L^{p}(B_1)$. Moreover, 
since $(p\theta_0-1)(\alpha-2)\geq 0$, Lemma \ref{log-est} implies that 
\begin{align*}
\frac{p-1}{p^2}\left\||x|^{\theta_0(\alpha-2)}\bigl|\log|x|\bigr|^{-\frac{2}{p}}u\right\|_{L^p(B_1)}^p
&
\leq\|\lambda u-L^1_tu\|_{L^p(B_1)}
\left\||x|^{\frac{p\theta_0-1}{p-1}(\alpha-2)}u\right\|_{L^p(B_1)}^{p-1}
\\
&
\leq\|\lambda u-L^1_tu\|_{L^p(B_1)}
\left\|u\right\|_{L^p(B_1)}^{p-1}
\\
&
\leq C\|\lambda u-L^1_tu\|_{L^p(B_1)}^p.
\end{align*}
This completes the proof.
\qed

\noindent Next we prove a-priori estimates for $L_t$, for large $t$,  by gluing the resolvents of $L_t^1$ and $L_t^2$.

\begin{prop} \label{glob-res-est}
For every $\gamma\in(0,1)$, there are constants $\tau, \rho, C_\gamma, C'>0$ such that if $t \ge \tau$, $\lambda \in \C_+$, $|\lambda | \ge \rho$,  and $u \in D_{p,\alpha}$ 
%$u\in C_c^\infty(\ov{B_1}\setminus\{0\})$, 
\begin{equation}\label{est-critical1}
\|u\|_{L^p}\leq \frac{C_\gamma}{|\lambda|^\gamma}
\|\lambda u-L_tu\|_{L^p}.
\end{equation}
and 
\begin{equation} \label{est-critical2}
\left\|
|x|^{\theta_0(\alpha-2)}\bigl|\log|x|\bigr|^{-\frac{2}{p}}u
\right\|_{L^p(B_{1/2})}\leq C'\|\lambda u-L_tu\|_{L^p}.
\end{equation}
%If $\theta_0=\frac{1}{p}$, then $\gamma=1$ is allowed.
\end{prop}
{\sc Proof.}  Observe that  if $t \ge 2^{2-\alpha}$ then $L^2_t$ coincides with $L$ in $B_{1/2}^c$. Since $\alpha < 2$, the function  $a(x) =|x|^\alpha$ satisfies the inequality $|\nabla a^{1/2}| \le C$ in  $B_{1/2}^c$, therefore, by \cite{for-lor}, the operator $L^2_t$ generates an analytic semigroup and the resolvent estimate 
\begin{equation} \label{resL2}
\|(\lambda-L_t^2)^{-1}\|_p\leq \frac{C}{|\lambda|}
\end{equation}
holds for $\lambda\in \C_+$, $|\lambda |\ge \rho$. By virtue of \cite[Section 5]{met-spi3} 
%and \cite{gluing},
 we can represent the resolvent of $L_t$ by gluing together the resolvents of  $L^1_t$ and $L^2_t$. In order to do this  we need  gradient estimates for $L_t^1$ and $L_t^2$  in an annulus $\Sigma_{r_1,r_2}=B_{r_2}\setminus B_{r_1} \subset  B_1\setminus B_{1/2}$. We fix $\frac 12 <s_1 <r_1 <r_2 <s_2 <1$ and we use
the classical interior estimates for uniformly elliptic operators. Since the coefficients of $L^1_t$  are uniformly bounded with respect to $t$ in the annulus $\Sigma_{s_1,s_2}$, there exists $C>0$, independent of $t>0$ such that
for every $u\in D_{p,\alpha}^1$ 
$$\|\nabla u\|_{p,\Sigma_{r_1,r_2}}\leq C\left (\eps\|L_t^1u\|_{p, \Sigma_{s_1,s_2}}+\frac{1}{\eps}\|u\|_{p,\Sigma_{r_1,r_2}}\right ).$$
Using (\ref{est-critical}) it follows that, for every $\lambda\in \C^+$ 
\begin{align*}
\|\nabla u\|_{p,\Sigma_{r_1,r_2}}&\leq C\left (\eps\|\lambda u- L_t^1u\|_{p, \Sigma_{s_1,s_2}}+\eps\|\lambda u\|_{p\Sigma_{s_1,s_2}}+\frac{1}{\eps}\|u\|_{p,\Sigma_{s_1,s_2}}\right )\\&\leq
 C\left (\eps\|\lambda u- L_t^1u\|_{p}+\eps|\lambda|^{1-\gamma}\|\lambda u-L_t^1 u\|_{p}+\frac{1}{\eps}|\lambda|^{-\gamma}\|\lambda u-L_t^1 u\|_p\right ).
\end{align*}
By choosing $\eps=|\lambda|^{-\delta}$ with $1-\gamma<\delta<\gamma$, we get for $r=\min\{-\delta,\delta-\gamma, 1-\delta-\gamma\}<0$
$$\|\nabla u\|_{p,\Sigma_{r_1,r_2}}\leq C|\lambda|^{r}\|\lambda u- L_t^1u\|_{p}.$$
In a similar way one proves gradient estimates for $L^2_t$.
Following the method of \cite[Section 5]{met-spi3} one constructs an approximate  resolvent $R(\lambda)$ for $L_t$, $$R(\lambda)=\eta_1(\lambda-L^1_t)^{-1}\eta_1+\eta_2(\lambda-L^2_t)^{-1}\eta_2$$ where $\eta_1,\eta_2$ are smooth functions supported in $B_{r_2}$, $B_{r_1}^c$ respectively and such that $\eta_1^2+\eta_2^2=1$. The operator $R(\lambda)$ satisfies $(\lambda-L_t)R(\lambda)=I+S(\lambda)$ and $\|S(\lambda)\| \le 1/2$ for $|\lambda|$ large, because of the gradient estimates. Then, for $\lambda \in \C_+ $, $|\lambda| $ large, we have
$$
(\lambda-L_t)^{-1}=R(\lambda) (I+S(\lambda))^{-1}
$$
and hence (\ref{est-critical1}) follows from (\ref{est-critical}) and (\ref{resL2}).
Estimate (\ref{est-critical2}) follows similarly from (\ref{w-est-critical}), since $\eta_2$ vanishes near $0$.
\qed   

\noindent {\sc Proof} (Theorem \ref{8.1}). We consider the operator $L_{int}$ defined in (\ref{Lint-critical1}). For $\lambda>0$, $\lambda-L_{int}$  is injective, by Proposition  \ref{inj}. As in the proof of Theorem \ref{gen-int} one sees that for $\lambda \ge \rho >0$, $f \in L^p(\R^N)$, $(\lambda-L_n)^{-1}f \to (\lambda-L_{int})^{-1}f$. By \eqref{est-critical1} it follows that  if  $\lambda \in \C_+$, $|\lambda | \ge \rho$,  $\lambda-L_{int}$ is invertible and,  for every $\gamma <1$
$$\|(\lambda-L_{int})^{-1}\| \le \frac{C_\gamma}{ |\lambda|^\gamma}.$$
For $s>0$ let $I_s:L^p \to L^p$ defined by $I_su(x)=u(sx)$. Clearly
$I_s$ is invertible with inverse $I_{s^{-1}}$ and
$\|I_su\|_p=s^{-N/p}\|u\|_p$.
Since $L=s^{2-\alpha}I_sLI_s^{-1}$, if $\lambda \in \C_+$,
$\lambda=r\omega$ with $|\omega|=\rho$  (hence $\omega$ belongs to the resolvent set)  then
the equality 
$$
\lambda-L_{int}=I_sr\left(\omega-\frac{s^{2-\alpha}L_{int}}{r}\right)I_s^{-1}
$$
with $s=r^\frac{1}{2-\alpha}$ shows that $\C_+$ is in the resolvent set and yields the decay
$$
\|(\lambda-L_{int})^{-1}\|_p \le \frac{C}{|\lambda|}\max\{\|(\omega-L_{int})^{-1}\|_p: |\omega|=\rho, \omega \in \C_+\},$$
For $\lambda >0$, positivity and coherence with respect to $p$ of $(\lambda-L_{int})^{-1}$ follow since 
$(\lambda-L_{int})^{-1}f=\lim_{n \to \infty}(\lambda-L_n)^{-1}f$. \qed

\subsection{Positive results for $\alpha>2$} \label{maggiore2}
The generation result proved in the critical case for $\alpha<2$ can  be extended by using similar arguments to the case $\alpha>2$. Recall that $s_0=\frac{N-2+c}{2}$.

\begin{teo}\label{8.1Maggiore}
 Assume that  
\begin{equation}\label{conj-criticalMaggiore}
s_0+2-\alpha\leq\frac{N}{p}\leq s_0
\end{equation}
 and define $L_{int}$ through the domain
\begin{equation} \label{Lint-criticalMaggiore}
D_{int}(L)=\{u\in D_{max}(L)\;;\;|x|^{\theta_0(\alpha-2)}\bigl|\log|x|\bigr|^{-\frac{2}{p}}u\in L^p(B^c_{2})\}, 
\end{equation}
where $\theta_0\in [0,1]$ satisfies $\frac{N}{p}=s_0+\theta_0(2-\alpha)$.
Then
$L_{int}$
generates a positive analytic semigroup in $L^p(\R^N)$  which is coherent with respect to all $p$ satisfying \eqref{conj-criticalMaggiore}. 
\end{teo}
 We only state the main steps.

\begin{prop}\label{injMaggiore}
For $\lambda >0$ the operator $\lambda-L_{int}$ is injective.
\end{prop}

\begin{prop}\label{log-est2}
Set \[
\widetilde{D}_1:=\{u\in C_c^\infty(\R^N\setminus B_1)\;;\;u=0\text{ on }\partial B_1\}.
\]
For every 
$v\in \widetilde{D}_1$
\begin{align}\label{log2}
Re\int_{B_1^c}
  |x|^{2-N}(\langle \nabla v,\nabla (\overline{v}|v|^{p-2})\rangle
\,dx
\geq 
\frac{p-1}{p^2}\int_{B_1^c}
|x|^{-N}\bigl|\log |x|\bigr|^{-2}|v|^{p}
\,dx.
\end{align}
In particular, 
if $u\in \widetilde{D}_1$, 
and $v=|x|^{\frac{N-2+c}{2}}u$, then 
\begin{align}
Re\int_{B_1^c}
  (-Lu)|x|^{(p\theta_0-1)(\alpha-2)}\overline{u}|u|^{p-2}
\,dx
&\geq 
\frac{p-1}{p^2}\int_{B_1^c}
|x|^{p\theta_0(\alpha-2)}\bigl|\log |x|\bigr|^{-2}|u|^{p}
\,dx.
\end{align}
\end{prop}

\begin{prop} %\label{res-est2}
For every $\gamma\in (0,1)$ there are constants $C_\gamma, C'>0$ such that for every $t>0$, $\lambda \in \C_+$ and $u\in \widetilde{D}_1$, 
%$u\in C_c^\infty(\ov{B_1}\setminus\{0\})$, 
\begin{equation}%\label{est-critical2}
\|u\|_{L^p(B_1^c)}\leq \frac{C_\gamma}{|\lambda|^{\gamma}}
\|\lambda u-L_tu\|_{L^p(B_1^c)}
\end{equation}
and 
\begin{equation}%\label{w-est-critical2}
\left\|
  |x|^{\theta_0(\alpha-2)}\bigl|\log|x|\bigr|^{-\frac{2}{p}}u
\right\|_{L^p(B_1^c)}\leq C'\|\lambda u-L_tu\|_{L^p(B_1^c)}, 
\end{equation}
where $L_t=L-k|x|^{\alpha-2}+k\min\{t, |x|^{\alpha-2}\}$.
\end{prop}
By using the propositions stated above, we deduce  Theorem \ref{8.1Maggiore} arguing as for Theorem \ref{8.1}.

\subsection{The equalities $L_{int}=L_{min}$ and $L_{int}=L_{max}$}

Here we investigate when $L_{int}$ coincides with $L_{min}$ or $L_{max}$. 

\begin{prop}\label{Lmax-critical} 
Assume $\alpha <2$ and $s_0 \le \frac{N}{p} \le s_0+2-\alpha$. Then 
\begin{itemize}
\item[(i)] If $\frac{N}{p}=s_0$, then $L_{int}=L_{max}$;
\item[(ii)] If $\frac{N}{p}\neq s_0$, then $N(\lambda-L_{max})\neq \{0\}$,  hence $ L_{int} \neq L_{max}$.
\end{itemize}
\end{prop}
{\sc Proof.} By the definition of $L_{int}$, see (\ref{Lint-critical}), (i) is obvious since $\theta_0=0$. We show (ii), that is, 
$N(\lambda-L_{max})\neq \{0\}$. 
We use Lemma \ref{bhv} with $k=\tilde k=0$ and take $v(x)=u_2(|x|)$.
Since $\frac{N}{p}>s_0$, 
we see from \eqref{u1-behave-res} that $v\in L^{p'}(\R^N)$. 
This means that $N(\lambda-L_{max})\neq \{0\}$.
\qed

\noindent By duality, the following proposition directly follows from Proposition \ref{Lmin-critical1}.

\begin{prop}\label{Lmin-critical1}
Assume $\alpha <2$ and $s_0 \le \frac{N}{p} \le s_0+2-\alpha$. Then 
\begin{itemize}
\item[(i)] If $\frac{N}{p}=s_0+2-\alpha$, then $L_{int}=L_{min}$, that is, $C_c^\infty(\Omega)$ is a core for $L_{int}$;
\item[(ii)] If $\frac{N}{p}\neq s_0+2-\alpha$, then $\overline{R(\lambda-L_{min})}\neq L^p(\R^N)$, hence $ L_{int} \neq L_{min}$.
\end{itemize}
\end{prop}

\noindent The case $\alpha>2$ is similar. 

\begin{prop}\label{Lmax-critical2}
Assume $\alpha >2$ and $s_0+2-\alpha \le \frac{N}{p} \le s_0$. Then
\begin{itemize}
\item[(i)] If $\frac{N}{p}=s_0$, then $L_{int}=L_{max}$;
\item[(ii)] If $\frac{N}{p}\neq s_0$, then $N(\lambda-L_{max})\neq \{0\}$,  hence $ L_{int} \neq L_{max}$.
\end{itemize}
\end{prop}

\begin{prop}\label{Lmin-critical12}
Assume $\alpha >2$ and $s_0+2-\alpha\le \frac{N}{p}\le s_0$. Then 
\begin{itemize}
\item[(i)] If $\frac{N}{p}=s_0+2-\alpha$, then $L_{int}=L_{min}$, that is, $C_c^\infty(\Omega)$ is a core for $L_{int}$;
\item[(ii)] If $\frac{N}{p}\neq s_0+2-\alpha$, then $\overline{R(\lambda-L_{min})}\neq L^p(\R^N)$, hence $ L_{int} \neq L_{min}$ .
\end{itemize}
\end{prop}

\noindent
Integrability of first and second derivatives for $u \in D_{int}(L)$ can be established as in Theorem \ref{gen-int}. For every $\theta <\theta _0$ we set 
$\alpha'=\alpha'(\theta)=\theta(\alpha-2)+2$ and define for $\theta_0>0$,
\begin{equation} \label{Dreg-crit}
D_{reg}(L)=\left\{
\begin{array}{l}
\left\{u\in D_{max}(L)\;; 
\begin{array}{l}
|x|^{\theta_0(\alpha-2)}\bigl|\log|x|\bigr|^{-\frac{2}{p}}u, 
\in L^p(B_{1/2}),
\\[5pt]
|x|^{\alpha'}D^2u, |x|^{\alpha'-1}\nabla u 
\in L^p(B)
\\[5pt]
|x|^{\alpha}D^2u, |x|^{\alpha-1}\nabla u \in L^p(B^c)
\end{array}
\forall \theta\in (0,\theta_0),\right\}
\ \ {\rm if\ } \alpha <2;
\\ \\
\left\{u\in D_{max}(L)\;; 
\begin{array}{l}
|x|^{\theta_0(\alpha-2)}\bigl|\log|x|\bigr|^{-\frac{2}{p}}u, 
\in L^p(B_{2}^c),
\\[5pt]
|x|^{\alpha'}D^2u, |x|^{\alpha'-1}\nabla u 
\in L^p(B^c)
\\[5pt]
|x|^{\alpha}D^2u, |x|^{\alpha-1}\nabla u \in L^p(B)
\end{array}
\forall \theta\in (0,\theta_0),\right\}
\ \ {\rm if\ } \alpha >2,
\end{array}
\right.
\end{equation}
 where $B=B_1$.
\begin{prop} \label{D_reg_crit}
If $\theta_0>0$, that is $N/p \neq s_0$, then the domains $D_{int}(L)$ and $D_{reg}(L)$ coincide.
\end{prop}
{\sc Proof. } Assume $\alpha <2$ and  let $u \in D_{int}(L)$. We  write $u=u_1+u_2$ where $u_1=u\phi$, $u_2=u(1-\phi)$ and $\phi \in C_c^\infty (\R^N)$ with support in $B_2$ and equal to 1 in $B_1$. We introduce the operator $L_2$ on $\R^N$ in this way: the coefficients of $L_2$ coincide with those of $L$ in $B_1^c$ whereas in $B_1$ they take the (constant) value that they have on $\partial B_1$. $L_2$ is therefore uniformly elliptic with Lipschitz coefficients in $B_1$ and satisfies Hypothesis 2.1 of \cite{for-lor}. By construction the function $u_2$ belongs to the maximal domain of $L_2$ and, by \cite[Proposition 2.9]{for-lor}, $|x|^\alpha D^2 u_2, |x|^{\alpha-1}\nabla u_2 \in L^p(B^c)$, that is  $|x|^\alpha D^2 u, |x|^{\alpha-1}\nabla u \in L^p(B^c)$. To treat $u_1$ we consider the operator $L_1=|x|^{\alpha'-\alpha}L$. Since $\alpha<2$, then $\alpha' \ge \alpha$ and then $u_1 \in D_{max}(L_1)$. Since $\alpha'-2=\theta (\alpha-2)>\theta_0(\alpha-2)$,  by the definition of $L_{int}$, $|x|^{\alpha'
-2}u_1 \in L^p(\R^N)$. By Lemma \ref{aux}, $u_1 \in D_{p,\alpha'}$. It follows that 
$|x|^{\alpha'}D^2u_1, |x|^{\alpha'-1}\nabla u_1, |x|^{\alpha'-2}u_1\in L^p(B)$, hence the same holds for $u$.
\qed
  \begin{os} \label{reg0} {\rm
The case $\theta_0=0$ or $N/p=s_0$ is quite special and we recall that $L_{int}=L_{max}$. Integrability of first and second derivatives can be obtained directly using Proposition \ref{L2}. If $\alpha <2$ and $u \in D_{int}(L)$, then $|x|^2 D^2u, |x| \nabla u \in L^p(B)$ and  $|x|^\alpha D^2u, |x|^{\alpha-1} \nabla u \in L^p(B^c)$  and conversely if $\alpha>2$. To see this we proceed as in the proposition above splitting $u=u_1+u_2$ and treating $u_2$ in the same way. Finally we note that $u_1 \in D_{max}(|x|^{2-\alpha}L)$ and then apply Proposition \ref{L2}.}
\end{os}

\noindent As in Section 5 one shows the minimality of $(\lambda-L_{int})^{-1}$, noting that the proof of Lemma \ref{max-princ} extends to the critical case, choosing $\theta_0$ such that $N/p+\theta_0(\alpha-2)=s_0$.
\begin{prop} \label{minimal-sem-crit}
Let $\lambda >0$, $f \ge 0$ and let $0 \le u \in D_{max}(L)$ satisfy $\lambda u -Lu=f$. Then $ (\lambda-L_{int})^{-1}f\le u$.
\end{prop}

\subsection{Negative results}
We show that if $\frac{N}{p}$ falls outside the closed interval (\ref{conjecture1}), then no realization $L_{min} \subset L \subset L_{max}$ generates a semigroup in $L^p(\R^N)$.

\begin{teo}\label{nonex-critical1}
\begin{itemize}
\item[(i)]
If $\alpha<2$ and $\frac{N}{p}> s_0+2-\alpha$, 
or $\alpha>2$ and $\frac{N}{p}< s_0+2-\alpha$. 
Then $N(\lambda-L_{min})\neq \{0\}$. 
\item[(ii)] If $\alpha<2$ and $\frac{N}{p}< s_0$, 
or $\alpha>2$ and $\frac{N}{p}> s_0$. 
Then $\overline{R(\lambda-L_{max})}\neq L^p(\R^N)$.
\end{itemize}
\end{teo}
{\sc Proof.} (i) We give a proof only for $\alpha<2$.
As in  Proposition \ref{counter1} we consider 
radial solutions of the equation 
\[
\lambda v-Lv=0.
\]
We use Lemma \ref{bhv} with $k=0$ and choose  $v=u_2$ so that  $v$ satisfies \eqref{u1-behave-res}. 
Since $s_0+2-\alpha<\frac{N}{p}$ this implies that 
$v, |x|^{\alpha-2}v\in L^p(\R^N)$ and hence, 
by Lemma \ref{aux}, we have 
$v\in D_{min}(L)$ and $\lambda v-Lv=0$. The proof of (ii) follows from (i), by duality.
\qed

%%%%%%%%%%%%%%%%%%%%%%%==7. Some examples==%%%%%%
\section{Examples}
In this section we specialize our results to particular operators.

\begin{esem}\label{sch-sub}
{\rm 
We consider Schr\"odinger operators with inverse square potential $L=\Delta-\frac{b}{|x|^2}$ 
(that is $\alpha=c=0$) assuming $b+\left(\frac{N-2}{2}\right)^2>0$. 
In this case 
\begin{equation*} 
%%%%%%%%%%%%%%%%%==red color==%%%%%%%%%%%%%%%%%%%
s_{1} = \frac{N-2}{2} - \sqrt{b+\left(\frac{N-2}{2}\right)^2}, 
\qquad
s_{2} = \frac{N-2}{2} + \sqrt{b+\left(\frac{N-2}{2}\right)^2}. 
%s_{1,2}=\frac{N-2}{2} \mp\sqrt{b+\left(\frac{N-2}{2}\right)^2}.
\end{equation*} 
Theorem \ref{gen-int} shows that $L_{int}$
endowed with the domain  (\ref{Dreg}) generates a positive analytic semigroup in $L^p(\R^N)$ 
if and only if $$s_1<\frac{N}{p}<s_2+2 \quad {\rm or\  } \quad  \left |\frac{N}{p}-\frac{N}{2}\right |< 1+\sqrt{b+\left(\frac{N-2}{2}\right)^2}.$$ 

\noindent
Observe that this improves the results in \cite{BG} and \cite{BV}. 
We point out that 
although generation results of analytic semigroup for $p$ in the sharp range above have already been proved in \cite[Section 4]{Li-sobol}, 
the description of domain of the generator $D_{int}(L)$ seems to be new.
Let us analyze it in more detail. By Theorem \ref{gen-min} and Proposition \ref{endpoint-Lmin}, if 
$s_1+2\leq \frac{N}{p}<s_2+2 $ 
then $L_{int}$ coincides with $L_{min}$. In particular, 
if $$s_1+2<\frac{N}{p}<s_2+2 \quad {\rm or\  } \quad  \left |\frac{N}{p}-\frac{N}{2}-1\right |< \sqrt{b+\left(\frac{N-2}{2}\right)^2},$$ then the domain  is given by 
\begin{equation*} 
D_{p,\alpha}=\{u\in L^p(\R^N) \cap  W^{2,p}_{loc}(\Omega),\ D^2u,\ |x|^{-1} \nabla u,\ |x|^{-2}u\in L^p(\R^N)\}.
\end{equation*}
We remark that in this case the generation result for $L_{min}$ is also stated in \cite[Section 3]{Sobajima-per}. 
If $s_1<\frac{N}{p} \le s_1+2$, then setting $\alpha'_1=s_1-N/p+2\in [0,2)$  we have from Theorem \ref{gen-int}
\begin{align*}
D_{int}(L)=&\Bigl \{u\in D_{max}(L)\cap W^{2,p}(B^c): 
|x|^{\alpha'}D^2u, |x|^{\alpha'-1}\nabla u, |x|^{\alpha'-2}u\in L^p(B)\ \  \forall \alpha'>\alpha'_1
 \Bigr \}.
\end{align*}
In Example \ref{sa} we show that, when $p=2$, then 
$L_{int}$ coincides with the Friedrich's extension of $L_{min}$. 

\noindent
When $N=1$, $L$ is the so called Calogero Hamiltonian and the above results  are of $L^p$ generalizations 
of the well-known properties of the Calogero operator in $L^2$. 
In fact. specializing to the case $N=1$ (where we recall that $\R$ should be substituted by $[0,\infty[$) and $p=2$, we obtain: 
if $b\geq \frac{3}{4}$, then 
$L_{\min}$ is nonnegative and selfadjoint and coincides with $L_{\max}$; 
if $-\frac{1}{4}\leq b < \frac{3}{4}$, then 
$L_{min}$ and $L_{max}$ are not selfadjoint but 
there exists a selfadjoint extension of $L_{min}$. The critical case $b=-\frac{1}{4}$ is explained in next example.

}
\end{esem}

\begin{esem}
{\rm Consider now the Schr\"odinger operators with inverse square potential $L=\Delta-\frac{b}{|x|^2}$ 
in the critical case $b+\left(\frac{N-2}{2}\right)^2=0$, where $s_0=\frac{N-2}{2}$. If
$s_0\leq\frac{N}{p}\leq s_0+2$ or (for $N \ge 3$)
$$
\frac{2N}{N+2} \le p \le \frac{2N}{N-2}
$$
then, by Theorem \ref{8.1}, there exists an operator $L_{int}$ such that 
$L_{min}\subseteq L_{int} \subseteq L_{max}$ and 
$L_{int}$ generates a positive analytic semigroup in $L^p(\R^N)$ 
which is coherent with respect to all $p$. 
It is worth noticing that the interval of generation is closed, 
in contrast with the case 
$b+(\frac{N-2}{2})^2>0$. 
Observe that this improves the results in \cite{BG}, \cite{BV} and \cite[Section 4]{Li-sobol} where,
although generation results of analytic semigroup for $p$ in the sharp range above have already been proved, 
the description of domain of the generator $D_{int}(L)$ seems to be new.
More precisely, by Proposition \ref{Lmin-critical1}, if 
$\frac{N}{p}=s_0+2 $ 
then $L_{int}$ coincides with $L_{min}$  and, by  Proposition \ref{Lmax-critical}, if $\frac{N}{p}=s_0$ 
then $L_{int}$ coincides with $L_{max}$.
If
$$
s_0\le \frac{N}{p}\leq s_0+2,
$$ then, by Proposition \ref{inj},  the domain  is given by 

\[
D_{int}(L)=\{u\in D_{max}(L)\;;\;|x|^{\theta_0(\alpha-2)}\bigl|\log|x|\bigr|^{-\frac{2}{p}}u\in L^p(B_{1/2})\}, 
\]
where $\theta_0\in[0,1]$ satisfies $\frac{N}{p}=s_0+2\theta_0$. Integrability of first and second derivatives for $u \in L_{int}$ is given in Proposition \ref{D_reg_crit}.

\noindent If $N=1$, then $b=-\frac{1}{4}$ and $L$ is the Calogero operator in $[0,\infty[$.
Then $L_{min}$ and $L_{max}$ are not selfadjoint but there exists a selfadjoint extension of $L_{min}$. 
In this case $L_{int}$ coincides with the Friedrich's extension of $L_{min}$ 
(this fact in a more general context is explained in Example \ref{sa}).}
\end{esem}

\begin{esem}\label{sa}
{\rm 
More generally, we can consider the formally selfadjoint operators 
\[
L=div(|x|^{\alpha}\nabla)-b|x|^{\alpha-2},
\]
which corresponds to $c=\alpha$ in \eqref{defL} and 
we focus our attention to $p=2$. 
If $\alpha=2$,  $L_{min}=L_{max}$ are selfadjoint and their  
domain is already given 
in \cite{met-soba-spi}, see also Proposition \ref{L2}.

\noindent We consider the case $\alpha\neq 2$ and $D_\alpha=b+(\frac{N-2+\alpha}{2})^2\geq 0$. 
Since $c=\alpha$ the function $f$ defined in (\ref{deff}) satisfies 
\begin{align*}
f(s)=b+s(N-2+\alpha-s), \qquad f\left(\frac{N-2+\alpha}{2}\right)=D_\alpha.
\end{align*}
If $D_{\alpha}>0$, then
condition \eqref{conjecture-other}, 
which is equivalent to \eqref{conjecture}, 
is satisfied with $\theta=\frac{1}{2}$. 
If $D_{\alpha}=0$, then the assumption in Theorem \ref{8.1} is satisfied with $\theta_0=\frac{1}{2}$.
Therefore, $L_{int}$ defined by 
\eqref{Lint} if $D_{\alpha}>0$ and \eqref{Lint-critical} if $D_{\alpha}=0$ 
generates an analytic semigroup in $L^2(\R^N)$ 
and its domain is characterized by \eqref{Dreg} id $D_\alpha >0$,  which gives a precise regularity. 
Moreover, $L_{int}$ is the  limit of $L_t$ in the resolvent sense 
(see Subsection 5.1 and Section 6) and each $L_t$ is nonnegative and selfadjoint for every $t>0$, since it coincides with $(L_t)_{min}$. This yields that $L_{int}$ is also selfadjoint.

\noindent 
It is worth noticing that, since $L_{min}$ is symmetric, 
$L_{min}$ is selfadjoint if and only if $L_{max}$ is selfadjoint. This means that the 
conditions on generation by $L_{min}$ and $L_{max}$ given in  Theorems \ref{gen-min}, \ref{gen-max}  coincide. 
This fact can be easily found via the identity
\begin{align}\label{value.sa}
f\left(\frac{N}{2}\right)=f\left(\frac{N}{2}+\alpha-2\right)
=D_{\alpha}-\left(\frac{\alpha-2}{2}\right)^{2}
\end{align}
Moreover, \eqref{value.sa} implies that if $D_{\alpha}\geq(\frac{\alpha-2}{2})^2$, that is, 
\[
b\geq -\left(\frac{N-2+\alpha}{2}\right)^2+\left(\frac{\alpha-2}{2}\right)^2,
\]
then $L_{\min}$ is nonnegative and selfadjoint 
and coincides with $L_{\max}$ (see also Proposition \ref{endpoint-Lmin}), hence with $L_{int}$. 
On the contrary, if $0\leq D_{\alpha}<(\frac{\alpha-2}{2})^2$, that is, 
\[
-\left(\frac{N-2+\alpha}{2}\right)^2
\leq b< -\left(\frac{N-2+\alpha}{2}\right)+\left(\frac{\alpha-2}{2}\right)^2,
\] 
then $L_{min}$ (and $L_{max}$) does not generate a semigroup in $L^2(\R^N)$ but its self-adjoint extension $L_{int}$ does it.  $L_{int}$ is the unique among the infinitely many self-adjoint extensions of $L_{min}$ which has the minimality property with respect to positive solutions, as explained in Propositions \ref{minimal-sem}, \ref{minimal-sem-crit}.

\noindent
Note that the constant $(\frac{\alpha-2}{2})^2$ with $\alpha=0$ coincides with the difference 
of the optimal constants in the usual Hardy and Rellich inequalities $\frac{N(N-4)}{4}-(\frac{N-2}{2})^2=1$.
%%% I think Okazawa want to add this comment.

\noindent
When 
$0\leq D_{\alpha}$ a self-adjoint extension of  
$L_{min}$ can be constructed by closing the nonnegative form 
\[
\mathfrak{a}(u,v)=\int_{\R^N}|x|^{\alpha}\nabla u\cdot\nabla \overline{v}\,dx+b\int_{\R^N}|x|^{\alpha-2}u\overline{v}\,dx,
\quad D(\mathfrak{a})=C_c^\infty(\Omega).
\]
This extension is called the Friedrich's extension $L_F$ of $L_{min}$ and  is 
one of nonnegative selfadjoint extensions of $L_{min}$ (not unique, unless $L_{min}$ itself is self-adjoint).  In general, $D(L_F)$, the domain of $L_F$, is given only formally. However, we can characterize $D(L_F)$ by showing that $L_F$ and $L_{int}$ coincide.
Since both operators generate semigroups, 
it suffices to observe that $D(L_F) \subset D_{int}(L)$ and
we divide the proof accordingly to $D_{\alpha}>0$ and $D_{\alpha}=0$.

\medskip

\noindent{ (The case $D_{\alpha}>0$).} 
We see from \eqref{reale} with $p=2$ and $c=\alpha$ that for every $u\in C_c^\infty(\Omega)$,
\begin{align*}
\mathfrak{a}(u,u)\geq D_\alpha \int_{\R^N}|x|^{\alpha-2}|u|^2\,dx
=D_{\alpha}\bigl\||x|^{\frac{\alpha-2}{2}}u\bigr\|_2^2.
\end{align*}
This implies that $D(\ov{\mathfrak{a}})\subset D(|x|^{\frac{\alpha-2}{2}})$, 
where $\ov{\mathfrak{a}}$ is the closure of the form $\mathfrak{a}$.
Hence, since $\theta=1/2$,  by the definition of $L_{int}$, see \eqref{Lint-critical}, we have
\[
D(L_F)\subset D_{max}(L)\cap D(\ov{\mathfrak{a}})\subset D_{int}(L).
\]

\medskip

\noindent{ (The case $D_{\alpha}=0$).} 
Set 
\[
U_0:=
\begin{cases}
B_{1/2} & \text{if }\alpha<2, 
\\
B_{2}^c & \text{if }\alpha>2, 
\end{cases}
\qquad 
U_1:=
\begin{cases}
B_{1} & \text{if }\alpha<2, 
\\
B_{1}^c & \text{if }\alpha>2
\end{cases}
\]
and let $\eta\in C_c^\infty(\R^N)$ satisfy $0 \le \eta \le 1$, 
$\eta\equiv 1$ in $U_0$ and $\eta\equiv 0$ in $U_1^c$. 
Then using \eqref{log1} if $\alpha<2$ and \eqref{log2} if $\alpha>2$, with $p=2$ and $\theta_0=1/2$,  we see that 
\begin{align*}
\left\|\chi_{U_0}|x|^{\frac{\alpha-2}{2}}\bigl|\log|x|\bigr|^{-1}u\right\|_2^2
\leq 
\int_{U_1}|x|^{\alpha-2}\bigl|\log|x|\bigr|^{-2}|\eta u|^2\,dx
\leq 
4\mathfrak{a}(\eta u, \eta u).
\end{align*}
Thus, we have 
\begin{align*}
\mathfrak{a}(\eta u, \eta u)
&
=
\int_{\R^N}
  \eta^2\left(|x|^{\alpha}|\nabla u|^2+b|x|^{\alpha-2}|u|^2\right)
\,dx
\\
&\quad+
2{\rm Re} \int_{U_0^c}
  |x|^{\alpha}(\eta\nabla \eta)\cdot(\ov{u}\nabla u)
\,dx
+
\int_{U_0^c}
  |x|^{\alpha}|\nabla \eta|^2|u|^2
\,dx
\\
&
\leq 
\mathfrak{a}(u,u)+
2\left\||x|^{\frac{\alpha}{2}}\nabla\eta\right\|_{\infty}\|u\|_2
\left\||x|^{\frac{\alpha}{2}}\nabla u\right\|_{L^2(U_1\setminus U_0)}
+\left\||x|^\frac{\alpha}{2}\nabla\eta\right\|_{\infty}^2\|u\|_2^2.
\end{align*}
Since
\begin{align*}
\left\||x|^{\frac{\alpha}{2}}\nabla u\right\|_{L^2(U_1\setminus U_0)}^2
&\leq
\int_{U_1\setminus U_0}
  \left(|x|^{\alpha}|\nabla u|^2+b|x|^{\alpha-2}|u|^2\right)
\,dx
+
|b|
\int_{U_1\setminus U_0}
  |x|^{\alpha-2}|u|^2
\,dx
\\
&\leq 
\mathfrak{a}(u,u)+2^{|\alpha-2|}|b|\, \|u\|_2^2,
\end{align*} 
we obtain 
\[
\left\|\chi_{U_0}|x|^{\frac{\alpha-2}{2}}\bigl|\log|x|\bigr|^{-1}u\right\|_2^2
\leq 
C\left(\mathfrak{a}(u, u)+\|u\|_2^2\right), 
\]
where $C$ is a constant independent of $u$. This implies that 
\[
D(\overline{\mathfrak{a}})\subset 
D\left(\chi_{U_0}|x|^{\frac{\alpha-2}{2}}\bigl|\log|x|\bigr|^{-1}\right).
\] 
Therefore, from \eqref{Lint-critical} we have 
\[
D(L_F)\subset D_{max}(L)\cap D(\ov{\mathfrak{a}})\subset D_{int}(L).
\]
}
\end{esem}

\begin{esem}
{\rm Let $b=c=0$, that is $L=|x|^\alpha\Delta$ and assume first that $N\neq 2$ so that 
$b+\left(\frac{N-2+c}{2}\right)^2=\left(\frac{N-2}{2}\right)^2>0$. 
Since  $s_1=0$, $s_2=N-2$,   $L_{int}$ endowed with the domain  (\ref{Dreg}) 
generates a positive analytic semigroup in $L^p(\R^N)$ if and only if 
$$\min\{0,2-\alpha\}<\frac{N}{p}<N-2+\max\{0,2-\alpha\}.$$ 
If $\alpha<2$, the condition reads $\frac{N}{p}<N-\alpha$ and, if $\alpha>2$, $\frac{N}{p}<N-2$.\\
By Theorem \ref{gen-min} $L_{int}=L_{min}$ if 
$2-\alpha<N/p<N-\alpha$ and the domain is given by
\begin{equation*} 
D_{p,\alpha}=\{u\in L^p(\R^N) \cap W^{2,p}_{loc}(\Omega),\ |x|^\alpha D^2u,\ |x|^{\alpha-1} \nabla u,\ |x|^{\alpha-2}u\in L^p(\R^N)\}.
\end{equation*}
By Theorem \ref{gen-max}, $L_{int}=L_{max}$ if $N/p<N-2$.
We  observe also that,  when $N=1$ and $\alpha\geq 2$, the interval of admissible $p$ is contained in the negative axis and the  operator is not a generator, as proved in \cite{met-spi2} for the operator $(1+|x|^\alpha)\Delta$. 

\noindent
Observe that this improves the results of \cite[Section 8]{met-spi1}. Indeed here we get a more precise description of the domain of $L$. Moreover we establish here non existence results for semigroups outside the above interval whereas in \cite{met-spi1} only the  non existence  of a positive semigroups is proved.

\noindent Let us consider the critical case $N=2$ where we have $s_0=0$. 
If  $\alpha > 2$, the interval $[s_0+\min\{0, 2-\alpha\}, s_0+\max\{0, 2-\alpha\}]=[2-\alpha,0]$ is contained in the negative real axis and therefore the operator $L$ is not a generator in any $L^p(\R^2)$. We point out that the same result has been obtained in \cite{met-spi2} for the operator $(1+|x|^\alpha)\Delta$. When $\alpha<2$, the operator  $L_{int}$, endowed with the domain 
\[
D_{int}(L)=\{u\in D_{max}(L)\;;\;|x|^{\theta_0(\alpha-2)}\bigl|\log|x|\bigr|^{-\frac{2}{p}}u\in L^p(B_{1/2})\}, 
\]
where $\theta_0\in(0,1]$ satisfies $\frac{2}{p}=\theta_0(2-\alpha)$,
generates an analytic semigroup for $p>\frac{2}{2-\alpha}$. In particular, if $\alpha< 0$, the operator $L_{int}$ generates an analytic semigroup for every $1< p<\infty$.  
}
\end{esem}

\begin{esem}
 {\rm Let $b=0$, that is $L=|x|^\alpha\Delta+c|x|^{\alpha-1}\frac{x}{|x|}\cdot \nabla$. If $N-2+c\neq 0$ then $\left(\frac{N-2+c}{2}\right)^2>0$ and $s_1=0$, $s_2=N-2+c$. 
By Theorem \ref{gen-int}, if $$\min\{0,2-\alpha\}<\frac{N}{p}<N-2+c+\max\{0,2-\alpha\},$$ 
$L_{int}$ endowed with the domain (\ref{Dreg}) generates a positive analytic semigroup in $L^p(\R^N)$. In particular, if $\alpha<2$, the condition reads $\frac{N}{p}<N-\alpha+c$ and, if $\alpha>2$, $\frac{N}{p}<N-2+c$.
By Theorem \ref{gen-min} $L_{int}=L_{min}$ when
$2-\alpha<N/p<N-\alpha+c$
and the domain is given by
\begin{equation*} 
D_{p,\alpha}=\{u\in L^p(\R^N) \cap  W^{2,p}_{loc}(\Omega),\ |x|^\alpha D^2u,\ |x|^{\alpha-1} \nabla u,\ |x|^{\alpha-2}u\in L^p(\R^N)\}.
\end{equation*}
By Theorem \ref{gen-max}, $L_{int}=L_{max}$ if $N/p<N-2+c$.
\\
Observe that this improves the results \cite{met-spi-tac}. Indeed here the degeneracy near the origin is also allowed and the domain description is more precise.

\noindent The critical case $b+\left(\frac{N-2+c}{2}\right)^2=0$ occurs for $c=2-N$ and we have $s_0=0$. 
As in the previous example, if  $\alpha > 2$, the interval $[s_0+\min\{0, 2-\alpha\}, s_0+\max\{0, 2-\alpha\}]=[2-\alpha, 0]$ is contained in the negative real axis, therefore the operator $L$ is not a generator in any $L^p(\R^N)$. The same phenomenon  has been already proved in  \cite{met-spi1} for the operator $(1+|x|^\alpha)\Delta+c|x|^{\alpha-1}\frac{x}{|x|}\nabla$. When $\alpha<2$, the operator $L_{int}$, endowed with the domain 
\[
D_{int}(L)=\{u\in D_{max}(L)\;;\;|x|^{\theta_0(\alpha-2)}\bigl|\log|x|\bigr|^{-\frac{2}{p}}u\in L^p(B_{1/2})\}, 
\]
where $\theta_0\in(0,1]$ satisfies $\frac{2}{p}=\theta_0(2-\alpha)$,
generates an analytic semigroup for $p>\frac{2}{2-\alpha}$. In particular, if $\alpha<0$, the operator $L_{int}$ generates an analytic semigroup for every $1 < p<\infty$.}
\end{esem}

\begin{esem}
{\rm For certain choices of the parameters $\alpha,\ b$ and $c$, $L_{int}$ can generate an analytic semigroup in $L^p(\R^N)$ even though $L$ is not dissipative for any $1<q<\infty$ and  $L_{max}$ does not generate for any $1<q<\infty$. Similarly, $L_{int}$ can generate an analytic semigroup in $L^p(\R^N)$ even though $L_{min}$ and  $L_{max}$ do not generate for any $1<q<\infty$.                                

\begin{itemize}\item[(a)]                                                        
\noindent  Assume that $b= 0$ and $N-2+c< 0$. It follows that $s_1=N-2+c< s_2= 0$. Therefore the operator $L_{max}$ does not generate an analytic semigroup for any $1<q<\infty$. If, in addition, we assume that $N\geq 2$ and $0\leq \alpha<2$, the dissipativity condition 
 $$s_1\leq \frac{N+\alpha-2}{p}\leq s_2$$ is never satisfied  but the generation condition for $L_{int}$ is valid for some $p$ since $s_2+2-\alpha> 0$.

\item[(b)]
\noindent We keep the conditions  $b= 0$ and $N-2+c < 0$  so that $s_1< s_2= 0$ and the operator $L_{max}$ does not generate for any $1<q<\infty$ but we assume $\alpha < 2$ and $N\leq s_1+2-\alpha$, that is $\alpha \le c$. 
It follows that $L_{min}$ never generates an analytic semigroup.
Finally observe that, since $s_1\leq 0$ and $s_2+2-\alpha>N$, the operator $L_{int}$ generates an analytic semigroup for every $1<p<\infty$.

\item[(c)]
\noindent If $L_{int}$ generates for some $p$ one can always find a $1<q<\infty$ such that $L$ is dissipative in $L^q$ or $L_{min}$ generates or $L_{max}$ generates in $L^q$.

\smallskip
\noindent In fact, assume that $\alpha < 2$ and that the generation condition for $L_{int}$ is true for some $1<p<\infty$ that is $$s_1<\frac{N}{p}<s_2+2-\alpha.$$
In order to violate the generation condition for $L_{max}$ for every $1<q<\infty$ we should have $s_1<s_2 \le 0$. Indeed, if $s_1<0<s_2$, then we can find some $q$ such that $s_1<\frac{N}{q}<s_2$. If $s_1$ and $s_2$ are positive, the generation condition for $L_{max}$ is violated only if $s_1 \ge N$ but this is not possible since  $s_1<\frac{N}{p}$.\\
\noindent Therefore we have: $s_1 <s_2 \le 0$ and $s_2+2-\alpha>0$. If $L_{min}$ does not generate in any $L^q$, then $s_1+2-\alpha \ge N$. If we choose $q$ such that
$$
\frac{N+\alpha-2}{s_1} <q <\frac{N+\alpha-2}{s_2}
$$
then $1<q<\infty$, $s_1 <(N+\alpha-2)/q <s_2$ and $L$ is dissipative in $L^q$.
\end{itemize}
}

\end{esem}

\begin{esem} \rm{
Let $L_1=|x|^\alpha \Delta+c|x|^{\alpha-1}\frac{x}{|x|}\nabla-b|x|^{\alpha-2}$ in the ball $B_R$, with Dirichlet boundary conditions. We define the domain of $L_1$  and deduce generation results for this operator in the ball by those in the whole space.

\begin{defi} 
\begin{align*}
D_{max}(L_1)= \{u \in L^p(B_R)\cap W^{2,p}\left (B_R \setminus B_\eps\right)\ \forall \eps>0: u(x)=0\ {\rm if\ } |x|=R, \  L_1u \in L^p(B_R)\}.
\end{align*}
\end{defi}
Observe that  the Dirichlet boundary condition $u(x)=0$ for $|x|=R$ makes sense, since $u$ has second derivatives in $L^p$ in a neighborhood of the boundary of $B_R$. By elliptic regularity $L_1$ is closed on its maximal domain.
If $\alpha \ge 2$ the function $a(x) =|x|^\alpha$ satisfies the inequality $|\nabla a^{1/2}| \le C$ in the ball $B_R$ even though not globally in $\R^N$ when $\alpha >2$. In analogy with \cite{for-lor} we define the domain of $L_1$ as follows.

\begin{defi} 
If $\alpha \ge 2$ we set
\begin{align*}
D_p(L_1)= \{u & \in L^p(B_R)\cap W^{2,p}\left (B_R \setminus B_\eps\right)\ \forall \eps>0: u(x)=0\ {\rm if\ } |x|=R,  \\ & |x|^{\alpha/2}\nabla u, |x|^\alpha D^2u \in L^p(B_R)\}.
\end{align*}
\end{defi}
\noindent
By the results in \cite{for-lor} and in \cite{met-spi1}, we immediately get generation for every $1<p<\infty$ when $\alpha\geq 2$.
\begin{prop} \label{easy-case}
If $\alpha \ge 2$, then $D_{max}(L_1)=D_p(L_1)$, the operator $L_1$ is closed on its domain and generates an analytic semigroup in $L^p(B_R)$ for every $1<p<\infty$.
\end{prop}
Concerning the case $\alpha<2$, the result proved in $\R^N$ is still true. It can be deduced by Theorem \ref{gen-int} by using the methods of \cite[Proposition 5.7]{met-spi1}.

\begin{prop} \label{main-palla}
Let $1<p<\infty$, $\alpha < 2$. If $b+(N-2+c)^2/4 > 0$, then a suitable realization of $L_{1,min} \subset L_{1,int} \subset L_{1,max}$ generates a semigroup in $L^p(B_R)$ if and only if $s_1<N/p<s_2+2-\alpha$. The semigroup is analytic and positive.
Moreover, setting  $\alpha'=\theta(\alpha-2)+2$,  $D_{int}(L_{1})$ is given by 
\begin{equation*}
\Bigl \{u\in D_{p,max}(L_1)\;; 
|x|^{\alpha'}D^2u, |x|^{\alpha'-1}\nabla u, |x|^{\alpha'-2}u\in L^p(B_R)\ \  {\rm for\ every\ } \theta \in I \Bigr \}
\end{equation*}
where $I$ is the interval of all $\theta\in [0,1]$ such that $f\left(\frac{N}{p}+\theta(\alpha-2)\right)>0$.
If $b+(N-2+c)^2/4=0$, set $s_0=\frac{N-2+c}{2}$,  a suitable realization of $L_{1,min} \subset L_{1,int} \subset L_{1,max}$ generates a semigroup in $L^p(B_R)$ if and only if $s_0 \le N/p \le s_0+2-\alpha$. The semigroup is analytic and positive.
Moreover, 
\begin{equation*} 
D_{int}(L_1)=\{u\in D_{max}(L_1)\;;\;|x|^{\theta_0(\alpha-2)}\bigl|\log|x|\bigr|^{-\frac{2}{p}}u\in L^p(B_{1/2})\}, 
\end{equation*}
where $\theta_0\in [0,1]$ satisfies $\frac{N}{p}=s_0+\theta_0(2-\alpha)$.
\end{prop}

\noindent Let us show the compactness of the resolvent for $\alpha<2$.
\begin{prop} 
Let $\alpha < 2$  and the assumptions of Proposition \ref{main-palla} be satisfied. Then the resolvent of $L_{1,int}$ is compact.
\end{prop}
Let us prove that  $D_{int}(L_{1})$ is compactly embedded
into $L^p(B_R)$. Consider the case $b+(N-2+c)^2/4>0$.  Let $\mathcal{U}$ be a bounded subset of
$D(L_{1,int})$. By the domain characterization, we obtain $\int_{B_R}|x|^{p(\alpha'-2)}|u|^p\leq M$ for
some positive $M$ and for every $u\in\mathcal{U}$. Since $\alpha' <2$, given
$\eps>0$, there exists $0<r<R$ such that $\int_{|x|<r}|u|^p<\eps^p$
for every $u\in \mathcal{U}$. Let $\mathcal{U}'$ be the set of the restrictions of the functions in $\mathcal{U}$ to $B_R\setminus B_r$.
Since the embedding of $W^{2,p}(B_R\setminus B_r)$ into $L^p(B_R\setminus B_r)$ is compact, the set $\mathcal{U}'$ which is bounded in $W^{2,p}(B_R\setminus B_r)$ is totally bounded in $L^p(B_R\setminus B_r)$. Therefore there exist $n\in\N$, $f_1,\ldots,f_n\in L^p(B_R\setminus B_r)$ such that
$$\mathcal{U}'\subseteq \bigcup_{i=1}^n\{f\in L^p(B_R\setminus B_r):\ \|f-f_i\|_{L^p(B_R\setminus B_r)}<\eps\}.$$
Set $\ov{f}_i=f_i$ in $B_R\setminus B_r$ and $\ov{f}_i=0$ in $B_r$. Then $\ov{f}_i\in L^p(B_R)$ and
$$\mathcal{U}\subseteq \bigcup_{i=1}^n\{f\in L^p(B_R):\ \|f-\ov{f}_i\|_{L^p(B_R)}<2\eps\}.$$
It follows that $\mathcal{U}$ is relatively compact in $L^p(B_R)$.\\
\noindent If $b+(N-2+c)^2/4=0$ and $\frac{N}{p}\neq s_0$ (that is $\theta_0 \neq 0$), the proof follows in similar way since the weight $|x|^{\theta_0(\alpha-2)}\bigl|\log|x|\bigr|^{-\frac{2}{p}}$ tends to $\infty$ as $x \to 0$.
If $\frac{N}{p}= s_0$, we consider the adjoint operator $\tilde {L}_{1,int}$ in $L^{p'}(B_R)$ whose domain is given as above with $\tilde {\theta}_0=1-\theta_0$=1, see the proof of Lemma \ref{re-int}. Then the resolvent of $\tilde {L}_{1,int}$ is compact in $L^{p'}(B_R)$, hence the semigroup, since it is analytic. By duality the semigroup generated by $L_{1,int}$ is compact in $L^p(B_R)$, hence the resolvent.
 \qed
}
\end{esem}

\begin{esem} {\rm
Let $L_2=|x|^\alpha \Delta+c|x|^{\alpha-1}\frac{x}{|x|}\nabla-b|x|^{\alpha-2}$ in the exterior domain $B_R^c$, with Dirichlet boundary conditions. We proceed as in the previous Example. 

\begin{defi} 
\begin{align*}
D_{max}(L_2)=\{u \in L^p(B_R^c)\cap W^{2,p}\left (B_R^c \cap B_r\right)\ \forall r>0: u(x)=0\ {\rm if\ } |x|=R, \  L_2u \in L^p(B_R^c)\}.
\end{align*}
\end{defi}
As before,  the Dirichlet boundary condition $u(x)=0$ for $|x|=R$ makes sense, since $u$ has second derivatives in $L^p$ in a neighborhood of the boundary of $B_R^c$. By local elliptic regularity, $L_2$ is closed on its maximal domain.
Observe that, when $\alpha \le 2$, the function  $a(x) =|x|^\alpha$ satisfies the inequality $|\nabla a^{1/2}| \le C$ in the exterior domain $B_R^c$. By following \cite{for-lor}, we can also define the domain $D_p(L_2$) as follows.

\begin{defi} 
If $\alpha \le 2$ we set
\begin{align*}
D_p(L_2)= \{u & \in L^p(B_R^c)\cap W^{2,p}\left (B_R^c \cap B_r\right)\ \forall r>0: u(x)=0\ {\rm if\ } |x|=R,   \\&
|x|^{\alpha/2}\nabla u,\ |x|^\alpha D^2u \in L^p(B_R^c)\}.
\end{align*}
\end{defi}
\noindent
As before 
the 
first generation results immediately follows 
from 
\cite{for-lor} and the results in \cite{met-spi1}. 
\begin{prop} \label{easy-case2}
If $\alpha \le 2$, then $D_{max}(L_2)=D_p(L_2)$, the operator $L_2$ is closed on its domain and, for every $1<p<\infty$ generates an analytic semigroup in $L^p(B_R^c)$.
\end{prop}
\noindent
In  the case $\alpha>2$, the following result can be deduced from Theorem \ref{gen-int} through  \cite[Proposition 5.6]{met-spi1}.

\begin{prop} \label{main-fuori-palla}
Let $1<p<\infty$, $\alpha >2$. If $b+(N-2+c)^2/4 > 0$, a suitable realization of $L_{2,min} \subset L_{2,int} \subset L_{2,max}$ generates a semigroup in $L^p(B_R^c)$ if and only if $s_1+2-\alpha<N/p<s_2$. The semigroup is analytic and positive.
Moreover, setting $\alpha'=\theta(\alpha-2)+2$,  $D_{int}(L_{2})$ is given by 
\begin{equation*}
\Bigl \{u\in D(L_{max})\;; 
|x|^{\alpha'}D^2u, |x|^{\alpha'-1}\nabla u, |x|^{\alpha'-2}u\in L^p(B_R^c)\ \  {\rm for\ every\ } \theta \in I \Bigr \}
\end{equation*}
where $I$ is the interval of all $\theta\in [0,1]$ such that $f\left(\frac{N}{p}+\theta(\alpha-2)\right)>0$.
If $b+(N-2+c)^2/4=0$, set $s_0=\frac{N-2+c}{2}$,  a suitable realization of $L_{2,min} \subset L_{2,int} \subset L_{2,max}$ generates a semigroup in $L^p(B_R)$ if and only if $s_0+2-\alpha \le N/p \le s_0$. The semigroup is analytic and positive.
Moreover, 
\begin{equation*} 
D_{int}(L_2)=\{u\in D_{max}(L_2)\;;\;|x|^{\theta_0(\alpha-2)}\bigl|\log|x|\bigr|^{-\frac{2}{p}}u\in L^p(B^c_{R})\}, 
\end{equation*}
where $\theta_0\in [0,1]$ satisfies $\frac{N}{p}=s_0+\theta_0(2-\alpha)$.
\end{prop}

\begin{prop} \label{comp1}
Let $\alpha>2$ and the conditions of Theorem \ref{main-fuori-palla} be satisfied.  Then the resolvent of
$L_2$ is compact in $L^p(B_R^c)$.\end{prop}

}
\end{esem}

\section*{Appendix}
\renewcommand{\thesection}{A}

\setcounter{lem}{0}

\begin{teo}[{\cite[Theorem 1.1]{Sobajima-per}}]\label{perturbation}
Let $A$ and $B$ be densely defined operators in Banach space $X$.
Assume that 
\begin{itemize}
\item[i)]
$-A$ generates a bounded analytic semigroup on $X$\/$;$ 
\item[ii)]
$D(A)\subset D(B)$ and there exists $\beta_0>0$ such that for every $u\in D(A)$, 
\begin{align}
\label{a1}
{\rm Re\,}(Au, F(Bu))_{X,X'}&\,\geq \beta_0\|Bu\|^2, 
\\
\label{a2}
(u, F(Bu))_{X,X'}&\,\geq 0, 
\end{align}
\end{itemize}
where $F$ is single-valued duality map from $X$ to the dual space $X'$.
Then for every $k\in \mathbb{C}$ satisfying ${\rm Re}\,k>-\beta_0$, 
$-(A+kB)$ with domain $D(A)$ also generates a 
bounded analytic semigroup on $X$. 
Moreover, if $A$ has a compact resolvent, then $A+kB$ also has a compact resolvent. 

\noindent
Assume further that $X$ is a Banach lattice, 
$A$ has positive resolvent and $B$ is also positive. Then 
$A+kB$ is also positive resolvent for every $k\in (-\beta_0,0)$.
\end{teo}

%\newpage


\begin{thebibliography}{99}

\bibitem{AGG06}%%%%%==
\refer{W. Arendt, J.A.\ Goldstein, G.R.\ Goldstein:} 
      {Outgrowths of Hardy's inequality, 
      Recent advances in differential equations 
      and mathematical physics, 51--68, 
      Contemp. Math.}{\bf 412}, 
      {Amer. Math. Soc. Providence, RI, 2006}. 

\bibitem{BG}
\refer{P. Baras, J.A. Goldstein:} 
{The heat equation with a singular potential,} 
{Trans. Amer. Math. Soc.}{284}{(1984), 121--139}.

\bibitem{BW}
\refbook{R. Beals, R. Wong:} 
{Special Functions,}{Cambridge Studies in Advances Mathematics,} 
{126} {(2010)}. 

\bibitem{BV}
\refer{H. Brezis, J.L. Vazquez: } 
{Blow-up solutions of some nonlinear elliptic problems,} 
{Rev. Mat. Univ. Compl. Madrid}{10}{(1997), 443--469}.


\bibitem{CM}
\refer{X. Cabre, Y. Martel: } 
{Existence versus explosion istantanee pour des equationes 
de la chaleur lineaires avec potentiel singulaires} 
{C.R. Acad. Sci. Paris}{329}{(1999), 973--978}.

\bibitem{costa}
\refer{D. G. Costa:}{Some new and short proofs for a class of Caffarelli-Kohn-Nirenberg type inequalities}%
{J. Math. Anal. Appl. }{337}{ (2008) 311--317}.


%\bibitem{engel-nagel}
%\refbook{K.J. Engel, R. Nagel:} {One parameter semigroups for linear evolutions equations,} {Springer-Verlag, Berlin, (2000)}.

\bibitem{for-lor}
\refer{S. Fornaro, L. Lorenzi:} 
{Generation results for elliptic operators with unbounded 
diffusion coefficients in $L^p$ and $C_b$-spaces,} 
{Discrete and Continuous Dynamical Systems} 
{A18} {(2007), 747--772}.

\bibitem{Gitman}
\refbook{ D. M. Gitman, I. V. Tyutin and B. L. Voronov:} 
{Self-adjoint extensions and spectral analysis in the
Calogero problem,}{J. Phys. A: Math. Theor.} 
{43} {(2010), 145--205}. 

\bibitem{Goldstein}
\refbook{J.A. Goldstein:} 
{Semigroups of Linear Operators and Applications}, 
{Oxford Mathematical Monographs, 
Oxford University Press, New York (1985)}. 

\bibitem{Kato}
\refbook{T. Kato:} 
{Perturbation Theory for Linear Operators,} 
{Grundlehren der math. Wissenschaten {\bf 132}, 
Springer-Verlag, Berlin and New York (1976)}. 

\bibitem{Li-sobol}
\refer{V. Liskevich, Z. Sobol, H. Vogt: }{On the $L^p$-theory 
for $C_0$-semigroups associated with second-order 
elliptic operators II,}{Journal of Functional Analysis} 
{193}{(2002), 55--76}. 

\bibitem{met-soba-spi}
\refer{G. Metafune, M. Sobajima, C. Spina:} 
{Weighted Calder\'on-Zygmund and Rellich inequalities 
in $L^p$,} 
{preprint (2013).}{}{}

\bibitem{gluing}
\refer{G. Metafune, M. Sobajima, C. Spina:}{Spectral properties of operators obtained by localization methods,}{preprint (2014).}{}


\bibitem{met-spi}
\refer{G. Metafune, C. Spina: } 
{An integration by parts formula in Sobolev spaces,} 
{Mediterranean Journal of Mathematics}{5}{(2008), 359--371}. 

\bibitem{met-spi2}
\refer{G. Metafune, C. Spina: } 
{Elliptic operators with unbounded 
diffusion 
coefficients 
in $L^p$ spaces,} 
{Annali Scuola Normale Superiore di Pisa Cl. Sc. (5),}
{11}{(2012), 
303--340
}.

\bibitem{met-spi1}
\refer{G. Metafune, C. Spina: } 
{A degenerate elliptic operators with unbounded coefficients,} 
{ Atti Accad. Naz. Lincei Cl. Sci. Fis. Mat. Natur. Rend. Lincei,} {}{to appear}.

\bibitem{met-spi3}
\refer{G. Metafune, C. Spina:} 
{Kernel estimates for some elliptic operators 
with unbounded coefficients,} 
{Discrete and Continuous Dynamical Systems} 
{A32 (6)}{(2012), 2285--2299}.


\bibitem{met-spi-tac}
\refer{G. Metafune, C. Spina, C. Tacelli: }{ Elliptic operators with unbounded diffusion and drift
coefficients in $L^p$ spaces,}{ Advances Diff. Equat., }
{}{ (to appear)}.

\bibitem{met-spi-tac1}
\refer{G. Metafune, C. Spina, C. Tacelli: }{On a class of elliptic operators with unbounded diffusion coefficients,}{ }
{}{preprint (2014)}.

%\bibitem{olver}
%\refbook{F.W.J. Olver:} 
%{Asymptotic and Special Functions,} 
%{A K Peters, (1997)}.

\bibitem{okazawa}
\refer{N. Okazawa:} {$L^p$-theory of Schr\"odinger operators 
with strongly singular potentials,} 
{Japan. J. Math.}{22}{(1996), 199--239}.

\bibitem{Oka-Soba}
\refer{N. Okazawa, M. Sobajima,} 
{$L^p$-theory for Schr\"odinger operators perturbed 
by singular drift terms,} 
{Proceedings of  the conference ``Differential Equations, Inverse Problems and Control Theory (2013)'', }{}{(to appear)}. 

\bibitem{OSY}
\refer{N. Okazawa, M. Sobajima, T. Yokota:} 
{Existence of solutions to heat equations 
with singular lower order terms,} 
{J. Differential Equations,}{256}{(2014), 3568-3593}. 

\bibitem{Sobajima}
\refer{M. Sobajima:} 
{$L^p$-theory for second-order elliptic operators 
with unbounded coefficients,} 
{J. Evol. Equ.}{12}{(2012), 957-971}.

\bibitem{Sobajima-per}
\refer{M. Sobajima:} 
{A class of relatively bounded perturbations 
for generators of bounded analytic semigroups 
in Banach spaces,} 
{J. Math. Anal. Appl.}{(416)}{(2014), 855-861}. 
,
\bibitem{soba-spin}
\refer{M. Sobajima, C. Spina: } 
{Generation results for some elliptic second order 
operators with unbounded diffusion coefficients,} 
{}{}{preprint (2014)}.

\bibitem{SW}
\refer{M. Sobajima, S. Watanabe:} 
{Landau-Lifschitz conjecture about the motion 
of a quantum mechanical particle 
under the inverse square potential,} 
{}{}{preprint (2013)}. 


\bibitem{spina}
\refer{C. Spina:} 
{Kernel estimates for some elliptic elliptic 
operators with unbounded diffusion coefficients 
in the one-dimensional and bi-dimensional cases,} 
{Semigroup Forum}{86 (1)}{(2013), 67--82}.

\bibitem{VZ}
\refer{J.L. Vazquez, E. Zuazua: } 
{The Hardy Inequality and the Asymptotic Behaviour 
of the Heat Equation with an Inverse Square Potential,} 
{Journal of Functional Analysis}{173}{(2000), 103--153}. 

\end{thebibliography}
\end{document}